\documentclass[11pt]{article} 

\usepackage[utf8]{inputenc}

\usepackage{graphicx}
\usepackage[ruled,vlined]{algorithm2e}
\usepackage{algorithmic}
\usepackage{epsfig}
\usepackage{epstopdf}
\usepackage{amssymb,amsmath}
\usepackage{amsmath}
\usepackage{amsthm}
\usepackage{amsfonts}
\usepackage{psfrag}
\usepackage{rotating}
\usepackage{latexsym}
\usepackage{dsfont}
\usepackage{stmaryrd}
\usepackage{amsthm}
\usepackage{amssymb}
\usepackage{amsmath}
\usepackage{mathtools}
\usepackage{mathrsfs}    
\usepackage{siunitx}
\usepackage{pgfpages}
\usepackage[Euler]{upgreek} 
\usepackage{isomath}
\usepackage{physics}
\usepackage{pifont} 
\usepackage{bbm}
\usepackage{empheq}
\usepackage[thicklines]{cancel}
\usepackage{subcaption}
\usepackage[authoryear,semicolon]{natbib}

\usepackage[force,almostfull]{textcomp}
\usepackage{fancyvrb}
\usepackage{textcomp}
\usepackage{url}
\usepackage{ifdraft}
\usepackage{url}
\usepackage{multirow}
\usepackage[colorlinks=true,linkcolor=blue,citecolor=green,urlcolor=red]{hyperref}
\usepackage{rotating}
\usepackage{xspace}
\usepackage{array}
\usepackage{multido}
\usepackage{enumerate}
\usepackage[utf8]{inputenc}
\newlength\figureheight 
\newlength\figurewidth 
\usepackage{graphics}
\usepackage{pgfplots}
\pgfplotsset{compat=newest}
\pgfplotsset{plot coordinates/math parser=false}
\allowdisplaybreaks

\usepackage{scalefnt}

\newtheoremstyle{specialcasestyle}{1mm}{1mm}{\upshape}{}{\bfseries\upshape}{.}{0mm}{}
\theoremstyle{specialcasestyle}

\newtheorem{assump}{Assumption}

\newtheorem{corr}{Corollary}

\newtheorem{rem}{Remark}
\newtheorem{thm}{Theorem}

\theoremstyle{remark}               

\DeclareMathOperator{\Law}{Law}

\newcommand{\E}{\mathbb{E}}

\newcommand{\Ptwo}{\mathcal{P}_2}

\usepackage{tikz}
\usetikzlibrary{positioning}
\usetikzlibrary{arrows.meta}

\theoremstyle{definition}

\usepackage{makecell}

\usepackage{enumitem}

\usepackage{cleveref}
\crefname{rem}{remark}{remark}
\Crefname{rem}{Remark}{Remark}
\crefname{assump}{assumption}{assumption}
\Crefname{assump}{Assumption}{Assumption}
\crefname{thm}{theorem}{theorem}
\Crefname{thm}{Theorem}{Theorem}
\crefname{prop}{proposition}{proposition}
\Crefname{prop}{Proposition}{Proposition}
\crefname{corr}{corollary}{corollary}
\Crefname{corr}{Corollary}{Corollary}

\usepackage{xcolor}
\usepackage[most]{tcolorbox}
\newtcolorbox{assumptionsA}[1][]{%
   enhanced,breakable,
   colback=blue!5!white,
   colframe=blue!75!black,
   title={Assumption (A)},
   nameref={Assumption (A)},
   #1
}
\newtcolorbox{assumptionsB}[1][]{%
   enhanced,breakable,
   colback=blue!5!white,
   colframe=blue!75!black,
   title={Assumption (B)},
   nameref={Assumption (B)},
   #1
}

\usepackage{placeins} 

\begin{document} 

\title{Computing expectations associated with McKean--Vlasov stochastic differential equations via probability-flow ordinary differential equations}
\author{\underline{Mathias Horres} \thanks{Corresponding author, Institute for Geometry and Practical Mathematics, RWTH Aachen University, Aachen, Germany ({\tt horres@igpm.rwth-aachen.de})},
 Shyam Mohan Subbiah Pillai \thanks{	Weierstrass Institute for Applied Analysis and Stochastics, Anton--Wilhelm--Amo Str. 39, 10117 Berlin, Germany ({\tt pillai@wias-berlin.de}).},
  and Ra\'ul Tempone \thanks{Computer, Electrical and Mathematical Sciences \& Engineering Division (CEMSE), King Abdullah University of Science and Technology (KAUST), Thuwal, Saudi Arabia ({\tt raul.tempone@kaust.edu.sa}).}  \thanks{This work was supported by the KAUST Office of Sponsored Research (OSR) under Award No. URF/1/2584-01-01 and the Alexander von Humboldt Foundation. For open access purposes, the authors have applied a Creative Commons Attribution (CC BY) license to any author-accepted manuscript version arising from this submission.}
}
\date{}
\maketitle
\thispagestyle{empty}
\begin{abstract}
This work considers the numerical approximation of expectations of functionals of $d$-dimensional McKean--Vlasov stochastic differential equations (MV-SDEs) with random initial conditions. Standard numerical schemes employ a stochastic interacting particle system of $P$ coupled $d$-dimensional SDEs, whose empirical measure converges to the MV-SDE law. However, naive Euler--Maruyama Monte Carlo (MC) estimators based on particle systems achieve a computational complexity of $\mathcal{O}\left( \text{TOL}^{-4} \right)$ to satisfy a prescribed error tolerance TOL. To reduce this complexity, this work reformulates the associated nonlinear Fokker--Planck equation as a Liouville equation describing the evolution of the law of a deterministic mean-field ordinary differential equation (ODE) with random initial conditions, sharing the same marginal distributions as the original MV-SDE. The resulting mean-field ODE eliminates the need to approximate It\^{o} integrals and permits the use of higher-order integration schemes in  time and stochastic space. However, simulating trajectories requires the score function of the MV-SDE's marginal law, which is approximated using a statistical estimator based on a system of $P$ coupled $d$-dimensional ODEs. The proposed estimator combines randomized quasi-MC in the stochastic variables and Runge--Kutta integrators in time, with efficient score estimation. Under assumptions on weak-error and variance convergence, the resulting  estimator achieves a computational complexity of $\mathcal{O}\left(\text{TOL}^{-\left(\max\{\beta/\alpha,2\} + 1/\gamma \right)}\right)$, where $\alpha$ and $\beta$ denote, respectively, the weak-error and variance convergence rates with respect to $P$,  and $\gamma$ is the order of the time-integration scheme.
\end{abstract}

\section{Introduction}
\label{sec:Introduction}
This work considers the approximation of $\mathbb{E}[g(X(T))]$, where $X$ solves a $d$-dimensional McKean--Vlasov stochastic differential equation (MV-SDE) \citep{mckean1966class} with random initial conditions, and $g:\mathbb{R}^d \rightarrow \mathbb{R}$ is a given scalar observable. MV-SDEs have drift and diffusion coefficients that depend on the law of the solution. \citet{mckean1966class} introduced MV-SDEs to model the collective behavior of large interacting particle systems. They arise across fields such as mathematical finance \citep{1100796777}, collective animal behavior \citep{RadekErban2012KineticandRelatedModels}, and oscillator systems \citep{acebron2005kuramoto}. The dependence on the law introduces significant complexity in the theoretical analysis and numerical approximation, on which a substantial body of literature exists \citep{mishura2020existence,Propagation_of_chaos_I,Propagation_of_chaos_II}.

The existence and uniqueness of solutions to MV-SDEs have been established under suitable assumptions by \citet{mckean1966class}, \citet{Sznitman1991TopicsIP}, and \citet{crisan2010approximate}, and more recently  \citet{hammersley2021weak} and \citet{mishura2020existence}. The Kolmogorov partial differential equations (PDEs) \citep{buckdahn2017mean}, \citep{crisan2018smoothing}, are of particular interest, especially the Kolmogorov forward equation, also called the Fokker--Planck equation, whose solution is the law of the MV-SDE. Unlike the linear case for standard SDEs, this equation is nonlinear. 

A direct approximation of the $d$-dimensional nonlinear Fokker--Planck equation becomes numerically infeasible in large dimensions. Hence, MV-SDEs are approximated using an interacting $P$-particle system, a system of $P$ coupled $d$-dimensional SDEs \citep{Sznitman1991TopicsIP}. Under suitable conditions, the interacting particle system converges in the mean-field limit as $P\to\infty$, which is the propagation-of-chaos phenomenon \citep{Sznitman1991TopicsIP}. The solution to the interacting particle system is approximated numerically using a suitable time-integration scheme with $N$ time steps. Monte Carlo (MC) estimators using this approach have a complexity of $\mathcal{O}\left(\text{TOL}^{-4}\right)$ \citep{MLMC} to achieve a prescribed error tolerance TOL. Multi-index Euler--Maruyama particle-system-based MC methods \citep{MLMC} further reduces the complexity to $\mathcal{O}\left(\text{TOL}^{-2} \log \left( \text{TOL}^{-1} \right)^2\right)$. 

This work aims to reduce this computational complexity based on the seminal work of \citet{GradientEstimation}. The Fokker--Planck equation for MV-SDEs is rewritten as a Liouville equation \citep{liouville1838note} of a $d$-dimensional probability-flow mean-field ODE. The Liouville equation describes the dynamics of a mean-field ODE with random initial conditions, that has the same marginal distribution as the original MV-SDE. With a smooth ODE flow, higher-order quadrature methods in time and stochastic space are employed to avoid lower-order It\^{o} integral discretization schemes. However, this ODE requires  the score function of the marginal law of the original MV-SDE. Such an estimator was first introduced by \citet{GradientEstimation} to approximate solutions to Fokker--Planck equations. A similar score-estimation problem also arises in score-based generative modeling \citep{song2021scorebasedgenerativemodelingstochastic}. Similar deterministic particle estimators have been explored in Bayesian inference: for example in \citep{pathiraja2019discretegradientscomputationalbayesian}, \citep{Reich_2019}, and \citep{carrillo2019blobmethoddiffusion} to approximate PDEs of Wasserstein gradient-flow type. The score estimator proposed by \citet{ScoreMatching} is employed in this work, resulting in a system of $P$-coupled $d$-dimensional ODEs. The estimator combines this score estimator with higher-order integrators in time and in stochastic space. Under suitable regularity assumptions, this achieves a lower complexity compared to the standard numerical approaches. 

Compared with \citet{GradientEstimation}, this work advances from deterministic particle solvers for Fokker--Planck equations to expectation estimation for MV-SDEs with a tolerance-driven work model. \citet{LU2024112859} focused on a Kullback--Leibler-based transport analysis, whereas the current work focuses on randomized quasi-MC (rQMC) sampling of the initial law and a tolerance-adaptive estimator. In this setting, rQMC is employed to sample only the initial law, whereas low-discrepancy sequences are propagated through a stochastic particle system with coupled noise in the work by \citet{rached2024forwardpropagationlowdiscrepancy}. \citet{song2021scorebasedgenerativemodelingstochastic} also considered probability-flow ODEs in the context of score-based generative modeling.


The contributions of this work are as follows:
\begin{itemize}
\item[(i)] Analytic framework: This work presents the Liouville reformulation for MV-SDEs and the score-based mean-field ODE estimator.
\item[(ii)] Conditional theory: Under the stated assumptions, this work derives a score-approximation error bound, rQMC weak error bound, and cost model.
\item[(iii)] Algorithmic development: An adaptive parameter-balancing strategy for a given tolerance TOL, revealing the predicted behavior.
\item[(iv)] Numerics: Numerical experiments are conducted on moderate-dimensional Gaussian and Kuramoto MV-SDEs.
\end{itemize}

This work extends the probability-flow ODE concept from \citet{GradientEstimation} to MV-SDEs, combines it with reproducing kernel Hilbert space (RKHS) score estimation, and derives rQMC-driven complexity bounds and a tolerance-adaptive algorithm. \Cref{sec:Mean-Field ODE Estimator} quantifies the four error sources and obtains bounds for each under sufficient conditions. \Cref{sec:Work Complexity} analyzes the work complexity of the estimator, and \Cref{sec:Adaptive Algorithm} provides the adaptive algorithm. Finally, \Cref{sec:Numerical Experiments} verifies the theoretical findings via numerical experiments on five-dimensional Gaussian and Kuramoto MV-SDEs. 

\section{Mean-field ordinary differential equation estimator}
\label{sec:Mean-Field ODE Estimator}

\subsection{McKean--Vlasov SDEs}
\label{ssec:McKean-Vlasov SDEs}

Let $\{W(t) \mid t\in[0,T]\}$ be a standard $d$-dimensional Wiener process on the filtered probability space $(\Omega,\mathcal{F},\{\mathcal{F}_t\}_{t\in[0,T]},\mathbb{P})$, and $X(0)$ be $\mathcal{F}_0$-measurable with $X(0) \sim \mu$ and independent of $\{W(t) \mid t\in[0,T]\}$. The filtration $\{\mathcal{F}_t\}_{t\in[0,T]}$ is the usual augmentation of $\sigma(X(0), W(s) \mid 0 \leq s \leq t)$. This work considers the stochastic process $X:[0,T]\times\Omega\to\mathbb{R}^d$ 
\begin{equation}\label{McKean-Vlasov SDE}
\begin{split}
dX(t) &= b(X(t),\mu_t) dt + \sigma(X(t),\mu_t) dW(t),\quad t\in (0,T]\\
X(0) &\sim \mu,
\end{split}
\end{equation}
where $\mathcal{P}(\mathbb{R}^d)$ denotes the set of probability measures on $\mathbb{R}^d$ and $\mu_t:= \text{Law}(X(t)) \in \mathcal{P}(\mathbb{R}^d)$. When $\mu_t$ admits a density with respect to the Lebesgue measure, it is denoted by $\rho_t$, i.e., $\mu_t(dx) = \rho_t(x) dx$, and $\nabla \log  \rho_t$ represents the corresponding score function. The drift and diffusion functions are given by $b:\mathbb{R}^d\times\mathcal{P}(\mathbb{R}^d) \to \mathbb{R}^d$, and $\sigma:\mathbb{R}^d\times\mathcal{P}(\mathbb{R}^d) \to \mathbb{R}^{d\times d}$. The initial measure $\mu \in\mathcal{P}(\mathbb{R}^d)$ is considered given. The associated Fokker--Planck equation of the MV-SDE \eqref{McKean-Vlasov SDE} is given by:
\begin{equation}  \label{Fokker Planck}
\begin{split}
\partial_t \rho_t(x) &=  -\sum_{i=1}^d \partial_{x_i} \left( b_i(x,\mu_t) \rho_t(x) \right) + \sum_{i,j=1}^d \partial_{x_i} \partial_{x_j}\left( D_{ij}(x,\mu_t) \rho_t(x) \right)  ,\quad t>0, x\in\mathbb{R}^d\\
\rho_0(x) &= \rho(x)
\end{split}
\end{equation}
where $D(x,\mu_t) = \frac{1}{2}\sigma(x,\mu_t)\sigma(x,\mu_t)^T$ \citep{Ambrosio2008GradientFlows,Bogachev2015FPK}. For $p\geq 1$, the general Wasserstein-$p$ distance between the probability measures $\mu, \nu \in\mathcal{P}_p(\mathbb{R}^d)$, the space of probability measures with a bounded moment of order $p$, is defined by
\begin{align*}
W_{p}(\mu,\nu) := \inf_{\pi \in \Pi(\mu,\nu)} \left( \int_{\mathbb{R}^d\times \mathbb{R}^d} \|x - y\|^p \pi(dx,dy) \right)^{1/p} = \inf_{ \substack{X\sim\mu \\ Y\sim\nu} } \mathbb{E}[\|X-Y\|^p]^{1/p}
\end{align*}
where $\Pi(\mu,\nu)$ denotes the set of all couplings of $\mu$ and $\nu$, that is, the set of probability measures on $\mathbb{R}^d\times \mathbb{R}^d$ with the first and second marginals respectively equal to $\mu$ and $\nu$ \citep{Propagation_of_chaos_I}. Based on the work by \citet{CarmonaDelarue2018MFG1}, and \citet{CarmonaDelarue2018MFG2}, this work assumes that $C>0$ exists such that, for all $x,y\in\mathbb{R}^d$ and all $\mu,\nu\in\mathcal{P}_2(\mathbb{R}^d)$, the following holds:
\begin{align}\label{Lipschitz MVSDE}
\|b(x,\mu) - b(y,\nu)\|_2 + \|\sigma(x,\mu) - \sigma(y,\nu)\|_2 \leq C\left( \|x-y\|_2 + W_2(\mu,\nu) \right).
\end{align}
Let $\mu\in \mathcal{P}_2(\mathbb{R}^d)$. Then, for any $T>0$, the MV-SDE \eqref{McKean-Vlasov SDE} has a strong solution on $[0,T]$ and its probability density is the weak solution to the nonlinear Fokker--Planck equation \eqref{Fokker Planck}. Let $g:\mathbb{R}^d\to\mathbb{R}$ be a given scalar-valued quantity of interest (QoI). This work is interested in the following expected value at final time $T$:
\begin{align}\label{QoI}
\mathbb{E}_{\mu_T}\left[ g(X(T)) \mid \mu_0 = \mu \right].
\end{align}
Directly solving \eqref{McKean-Vlasov SDE} is infeasible due to the dependence of $X(t)$ on its measure $\mu_t$. Solving the $d$-dimensional nonlinear parabolic Fokker--Planck equation \eqref{Fokker Planck} is computationally cumbersome. Instead, strong approximations to the MV-SDE \eqref{McKean-Vlasov SDE} can be employed. For $P$ particles in $\mathbb{R}^d$, where $\mathbf{X}(t) = (X_1^P(t),\ldots,X_P^P(t))^T$, the following empirical measure is defined:
\begin{align*}
\bar{\mu}^{\mathbf{X}(t)} = \frac{1}{P} \sum_{p=1}^P \delta(\cdot - X_p^P(t)),
\end{align*}
where $\delta(\cdot - x)$ denotes the delta distribution centered at $x$. Then, the MV-SDE \eqref{McKean-Vlasov SDE} is approximated using an interacting stochastic particle system \citep{Sznitman1991TopicsIP}, where $\bar{\mu}^{\mathbf{X}(t)}$ approximates $\mu_t$. For $p=1,\ldots,P$, $X_p^P(t)$ solves the coupled SDEs
\begin{equation}\label{Interacting SDE}
\begin{split}
dX_p^P(t) &= b(X_p^P(t),\bar{\mu}^{\mathbf{X}(t)}) dt + \sigma(X_p^P(t),\bar{\mu}^{\mathbf{X}(t)}) dW_p(t),\quad t\in (0,T]\\
X_p^P(0) &\sim \mu,
\end{split}
\end{equation}
where $\{W_p\}_{p=1}^P$ represents mutually independent $d$-dimensional Wiener processes, independent of the starting points $\{X_p^P(0)\}_{p=1}^P$. Define
\begin{equation}\label{PD drift diffusion}
\begin{split}
&\mathbf{b}:\mathbb{R}^{Pd} \to \mathbb{R}^{Pd},~ \mathbf{X}\mapsto ( b(X_1,\bar{\mu}^{\mathbf{X}}),\ldots, b(X_P,\bar{\mu}^{\mathbf{X}})  )^T\\
&\boldsymbol{\sigma}:\mathbb{R}^{Pd} \to \mathbb{R}^{Pd \times Pd},~ \mathbf{X}\mapsto \text{diag}( \sigma(X_1,\bar{\mu}^{\mathbf{X}}),\ldots, \sigma(X_P,\bar{\mu}^{\mathbf{X}}  )),\\
&\mathbf{D}:\mathbb{R}^{Pd} \to \mathbb{R}^{Pd \times Pd},~ \mathbf{X} \mapsto \frac{1}{2}\boldsymbol{\sigma}(\mathbf{X}) \boldsymbol{\sigma}(\mathbf{X})^T,
\end{split}
\end{equation}
where $\mathbf{X} = (X_1,\ldots,X_P)^T$  is constructed using $P$ $d$-dimensional components. Then, the $Pd$-dimensional linear Fokker--Planck equation associated with the interacting particle system \eqref{Interacting SDE} is given as follows:
\begin{equation} \label{Fokker Planck interacting}
\begin{split}
\partial_t \boldsymbol{\rho}_t(\mathbf{x}|\mathbf{y}) &=  - \nabla_{\mathbf{x}} \cdot\left( \mathbf{b}(\mathbf{x})\boldsymbol{\rho}_t(\mathbf{x}|\mathbf{y}) - \nabla_{\mathbf{x}} \cdot\left( \mathbf{D}(\mathbf{x})\boldsymbol{\rho}_t(\mathbf{x}|\mathbf{y}) \right) \right),\quad t>0, \mathbf{x}\in\mathbb{R}^{Pd}\\
\boldsymbol{\rho}_0(\mathbf{x}|\mathbf{y}) &= \boldsymbol{\delta}(\mathbf{x} - \mathbf{y} ),
\end{split}
\end{equation}
where $\boldsymbol{\rho}_t(\cdot|\mathbf{y})$ represents the joint conditional density of $\mathbf{X}^P(t) = (X_1^P(t),\ldots,X_P^P(t))^T$, given that $\mathbf{X}^P(0) = \mathbf{y}\in\mathbb{R}^{Pd}$, which is distributed according to $\boldsymbol{\mu}$. Solving the Fokker--Planck equation \eqref{Fokker Planck interacting} numerically up to a desired tolerance quickly becomes infeasible due to the large dimension $Pd$ of the system.

In general, propagation-of-chaos estimates are available in strong or Wasserstein metrics at an order of $\mathcal{O}\left( P^{-1/2}\right)$, whereas weak rates for smooth observables can improve under additional regularity \citep{Sznitman1991TopicsIP}. Under additional assumptions on drift and diffusion, a weak-convergence rate of $\mathcal{O}\left( P^{-1}\right)$ has been shown \citep{haji2021simple}. In this case, MC methods achieve a complexity of $\mathcal{O}\left( \text{TOL}^{-4} \right)$ \citep{MLMC}. 

\subsection{Mean-field ordinary differential equations}
\label{ssec:Mean-Field ODE}
This work aims to reduce the $\mathcal{O}\left(\text{TOL}^{-4}\right)$ computational complexity using probability-flow ODEs. Instead of using \eqref{Interacting SDE} to approximate \eqref{QoI}, this study follows the work by \citet{GradientEstimation} and constructs a mean-field ODE, that captures the time-marginals of the MV-SDE \eqref{McKean-Vlasov SDE} and avoids It\^{o} integrals, allowing the use of higher-order Runge--Kutta (RK) methods. All randomness of the resulting ODE stems from the initial measure; hence rQMC methods can be applied by exploiting the smoothness of the ODE flow. The mean-field ODE is obtained by rewriting the Fokker--Planck equation \eqref{Fokker Planck} in the form of a Liouville equation:
\begin{equation} \label{sec. exp: Liouville equation}
\begin{split}
\partial_t \rho_t(x) &= - \nabla_x \cdot \left( \hat{b}(x,\mu_t)\rho_t(x) \right),\quad t>0, x\in\mathbb{R}^d\\
\rho_0(x) &= \rho(x),
\end{split}
\end{equation}
with velocity
\begin{align}\label{mean-field-ode-velocity}
\hat{b}:\mathbb{R}^d\times\mathcal{P}(\mathbb{R}^d) \to\mathbb{R}^d,~ (x,\mu_t)\mapsto b(x,\mu_t) - \nabla_x \cdot D(x,\mu_t) - D(x,\mu_t) \nabla_x \log \rho_t(x),
\end{align}
where $(\nabla_x \cdot D(x,\mu_t))_i = \sum_{j=1}^d \partial_{x_j} D_{ij}(x,\mu_t)$, $i=1,\ldots,d$. The Liouville equation \eqref{sec. exp: Liouville equation} can be interpreted as the PDE satisfied by the probability density function of the process $\hat{X}:[0,T]\times \Omega\to \mathbb{R}^d$ that follows the mean-field ODE:
\begin{equation} \label{sec. exp: mean field ODE true gradlog}
\begin{split}
\frac{d\hat{X}(t)}{dt} &= \hat{b}\left(\hat{X}(t),\mu_t\right),\quad t\in (0,T]\\
\hat{X}(0) &\sim \mu.
\end{split}
\end{equation}

By itself, this approach is not yet useful due to the score function $D(x,\mu_t)\nabla_x \log \rho_t(x)$ in the mean-field ODE \eqref{sec. exp: mean field ODE true gradlog}. One might be tempted to approximate the measure $\mu$ in \eqref{sec. exp: mean field ODE true gradlog} with an empirical measure of $P$ particles $\mathbf{x} = (x_1,\ldots,x_P)^T$
\begin{align*}
\bar{\mu}^\mathbf{x} = \frac{1}{P}\sum_{p=1}^P \delta(\cdot - x_p).
\end{align*}
However, the score of an empirical measure is not well-defined. Therefore, the next subsection focuses on an efficient estimation of the score.

\begin{rem}
~The existence and uniqueness of the \eqref{sec. exp: mean field ODE true gradlog} are only assumed. Due to the score dependence, the velocity is irregular, a topic still under investigation \citep{banos2015,Bayer2024,bauer2018strongsolutionsmeanfieldstochastic,REN2023291}.
However, $\hat{X}(t)$ is not directly approximated. Instead, the score is approximated, yielding a more regular velocity and a second mean-field ODE $\hat{X}^\infty(t)$. This second mean-field ODE is approximated using an interacting particle system, and the existence and uniqueness of $\hat{X}^\infty(t)$ hold in the limit $P\to\infty$ under \nameref{ass:A}, introduced below.
\end{rem}

\subsection{Score estimation}
\label{ssec:Score Estimation}
Hyvärinen's score matching \citep{ScoreMatching} is employed to approximate the unknown $\nabla\log\rho_t$ in some approximation space $\mathcal{V}$, as in the work by \citet{GradientEstimation}. This work restricts itself to the Gaussian RKHS $\mathcal{V}=\mathcal{H}$, defined by the kernel functions
\begin{align*}
K:\mathbb{R}^d\times\mathbb{R}^d \to\mathbb{R},~ (x,y)\mapsto \exp\left( -\frac{\|x-y\|_2^2}{2\ell^2} \right)
\end{align*}
for a length scale $\ell\in\mathbb{R}^+$. Each component of $\nabla\log\rho_t$ is approximated via the minimization problem:
\begin{equation}\label{argmin}
\begin{split}
\partial_{x_i} \log \rho_t \approx r_i + \arg\min_{\phi\in\mathcal{H}} \left( J_i^r[\phi,\mu_t] + \lambda \|\phi\|_{\mathcal{H}}^2 \right),\quad i = 1,\ldots,d,
\end{split}
\end{equation}
where
\begin{align*}
&J_i^r[\phi,\mu_t] 
:= \int_{\mathbb{R}^d} \left( \phi^2(x) + 2r_i(x)\phi(x) + 2\partial_{x_i} \phi(x) \right) \ \mu_t(dx),
\end{align*}
and $\lambda>0$ is a regularization parameter. The functions $r_i \in \mathcal{H}$ can be freely chosen and are useful when further information about the score is known, e.g., the limit measure $\mu_\infty$, see \citep{GradientEstimation}. Other minimization spaces $\mathcal{V}$ are also possible, such as other RKHSs, polynomial spaces, or even other minimization approaches. The solution of the minimization problem \eqref{argmin} is written as follows:
\begin{equation}\label{explicit argmin}
r_i + \arg\min_{\phi\in\mathcal{H}}  \left( J_i^r[\phi,\mu_t] + \lambda \|\phi\|_{\mathcal{H}}^2 \right) = r_i + \phi_{i,\mu_t}^*,\quad i=1,\ldots,d,
\end{equation}
and is given by
\begin{align*}
 \left(\lambda I + \mathcal{C}_\mu\right)\phi_{i,\mu}^* =& - \left( \mathcal{C}_\mu r_i + \xi_{i,\mu} \right),
\end{align*}
where $I:\mathcal{H}\to\mathcal{H}$ represents the identity and 
\begin{equation*}
\mathcal{C}_\mu:\mathcal{H} \to \mathcal{H},\quad  \mathcal{C}_\mu f := \int_{\mathbb{R}^d} K(\cdot,y) f(y) \ \mu(dy), \quad \xi_{i,\mu} := \int_{\mathbb{R}^d} \partial_{x_i} K(x,\cdot) \ \mu(dx)\in\mathcal{H}, \quad i = 1,\ldots,d.
\end{equation*}
See \Cref{sec:Score Estimator} for further details. Let us denote $r = \left( r_1,\ldots,r_d\right)^T$, and $\phi_{\mu}^* :=\allowbreak \Big( \phi_{1,\mu}^*,\ldots, \phi_{d,\mu}^*\Big)^T$. For the RKHS $\mathcal{H}$, this study defines a new mean-field ODE $\hat{X}^\infty:[0,T]\times\Omega\to\mathbb{R}^d$ given by
\begin{equation}\label{mean-field ODE approximation}
\begin{split}
\frac{d \hat{X}^\infty(t)}{dt} &= b^\mathcal{H}(\hat{X}^\infty(t),\mu_t^\infty),\quad t\in (0,T]\\
\hat{X}^\infty(0) &\sim \mu,
\end{split}
\end{equation}
with velocity
\begin{align}\label{mean-field-ode-score-veclocity}
b^\mathcal{H}:\mathbb{R}^d\times\mathcal{P}(\mathbb{R}^d) \to \mathbb{R}^d,~ (x,\mu)\mapsto  b(x,\mu) - \nabla \cdot D(x,\mu) - D(x,\mu)\left( r(x) +  \phi_{\mu}^*(x)\right).
\end{align}
The notation $b^\mathcal{H}$ indicates that the velocity field of \eqref{mean-field ODE approximation} uses an approximation of $\nabla \log \rho_t$ in $\mathcal{H}$, where $\mu_t^\infty\in\mathcal{P}(\mathbb{R}^d)$ denotes the law of $\hat{X}^\infty(t)$. To summarize, $\mu_t = \text{Law}(X(t))$, and $\rho_t = \frac{d\mu_t}{dx}$, just as $\mu_t^\infty = \text{Law}(\hat{X}^\infty(t))$, and $\rho_t^\infty = \frac{d\mu_t^\infty}{dx}$. Before proceeding, a set of standing \nameref{ass:A} is stated.

\begin{assumptionsA}\label{ass:A}
\textbf{Regularity of measures:} This work assumes that $\mu_t$ and $\mu_t^\infty$ admit strictly positive densities $\rho_t$ and $\rho_t^\infty$ on the relevant domain for $t\in[0,T]$, with sufficient spatial regularity for the Fokker--Planck and Liouville equations to hold in the stated sense. This work assumes periodic/no-flux boundary conditions, compact-domain boundary conditions compatible with the flux, or sufficient decay at infinity to ensure that all integrations by parts in the Liouville reformulation and in score matching have no boundary contribution.\\
\textbf{Boundedness:} For all $x\in\mathbb{R}^d$, and $\mu\in\mathcal{P}(\mathbb{R}^d)$,
\begin{equation}\label{D bound}
\begin{split}
\|D(x,\mu)\|_2 &\leq B,\\
\|\nabla\cdot D(x,\mu)\|_2 &\leq B.
\end{split}
\end{equation}
Similarly, for all $x,y\in\mathbb{R}^d$, and $i=1,\ldots,d$,
\begin{equation}\label{A kernel bound}
\begin{split}
&K(x,y) \leq B,\\
&|\partial_{x_i} K(x,y)| \leq B,\\
&\|\partial_{x_i} K(x,\cdot)\|_{\mathcal{H}} \leq B,\\
&|r_i(x)| \leq B.
\end{split}
\end{equation}

\textbf{Lipschitz continuity:} For all $x,y\in\mathbb{R}^d$ and $\mu,\nu\in\mathcal{P}(\mathbb{R}^d)$, this work assumes
\begin{equation}\label{Lipschitz_divergence}
\begin{split}
\|\nabla \cdot D(x,\mu) - \nabla \cdot  D(y,\nu)\|_2 &\leq C(\|x-y\|_2 + W_2(\mu,\nu))
\end{split}
\end{equation}
for a positive constant $C$. For all $x,y,z\in\mathbb{R}^d$, $i=1,\ldots,d$,
\begin{equation}
\begin{split}
\|r_i(x) - r_i(y)\|_2&\leq C\|x-y\|_2,\\
|K(x,z) - K(y,z)|&\leq C\|x-y\|_2, \\
|\partial_{x_i} K(x,z) - \partial_{x_i} K(y,z)| &\leq C\|x-y\|_2,\\
\|g(x) - g(y)\|_2&\leq \text{Lip}(g) \|x-y\|_2,
\end{split}
\end{equation}
with positive constants $C$ and $\text{Lip}(g)$.\\
\textbf{Growth condition:} For the multi-index $\kappa = (\kappa_1,\ldots,\kappa_d)$ of order $|\kappa| = \sum_{i=1}^d \kappa_i$, the multi-index partial derivative operator is defined as
\begin{align*}
D_x^\kappa = \frac{\partial^{\kappa_1 + \ldots + \kappa_d}}{\partial x_1^{\kappa_1} \ldots \partial x_d^{\kappa_d}}.
\end{align*}
The QoI $g:\mathbb{R}^d\to\mathbb{R}$ is assumed to satisfy
\begin{align}
|D_x^\kappa g(x)| \leq C(1+\|x\|_2^{R_g}),\quad|\kappa|  \leq 2,
\end{align}
where $C$, and $R_g\leq 2$ are positive constants.

For $\mathbf{x}\in\mathbb{R}^{Pd}$ with $d$-dimensional components $x_p$, $p=1,\ldots,P$ and $d$-dimensional empirical measure $\bar{\mu}^\mathbf{x} := \frac{1}{P}\sum_{p=1}^P \delta(\cdot - x_p)$ the  following function is defined:
\begin{align*}
\nabla\cdot \mathbf{D}(\mathbf{x}) & = \left( \nabla\cdot D(x_1,\bar{\mu}^\mathbf{x}),\ldots, \nabla\cdot D(x_P,\bar{\mu}^\mathbf{x})  \right)^T\in \mathbb{R}^{Pd}.
\end{align*}
For all $\mathbf{x}\in\mathbb{R}^{Pd}$, similar growth rates are assumed:
\begin{equation}
\begin{split}
\|\mathbf{b}(\mathbf{x})\|_2 & \leq C(1+\|\mathbf{x}\|_2),\\
\|\nabla \cdot \mathbf{D}(\mathbf{x})\|_2& \leq C(1+\|\mathbf{x}\|_2),\\
\|\mathbf{D}(\mathbf{x})\|_2 &\leq C(1+\|\mathbf{x}\|_2)
\end{split}
\end{equation}
with a positive constant $C$. Moreover, the following defines the growth for all derivatives up to and including those of order $2$:
\begin{equation}
\begin{split}
\|D_\mathbf{x}^\kappa \mathbf{b}(\mathbf{x})\|_2& \leq C'(1+\|\mathbf{x}\|_2^R),\\
\| D_\mathbf{x}^\kappa \nabla \cdot \mathbf{D}(\mathbf{x})\|_2& \leq C'(1+\|\mathbf{x}\|_2^R),\\
\|D_\mathbf{x}^\kappa \mathbf{D}(\mathbf{x})\|_2& \leq C'(1+\|\mathbf{x}\|_2^R),\quad  |\kappa| \leq 2
\end{split}
\end{equation}
with the positive constants $C'$, and $R$. This study also assumes that:
\begin{equation}
\begin{split}
\|D_x^\kappa r_i(x)\|_2& \leq C'(1+\|x\|_2^R),\quad i=1,\ldots,d\\
\| D_x^\kappa K(x,y)\|_2& \leq C'(1+\|x\|_2^R),\quad  |\kappa| \leq 2
\end{split}
\end{equation}
for all $x,y\in\mathbb{R}^d$, with positive constants $C'$, and $R$.\\
\textbf{Initial measure:} This work imposes $\mu \in \mathcal{P}_2(\mathbb{R}^d) \cap \mathcal{P}_q(\mathbb{R}^d)$ for $q\geq 2R+2$ and assumes absolute continuity with respect to the Lebesgue measure.\\
\textbf{RKHS:} The kernel $K:\mathbb{R}^d\times\mathbb{R}^d\to\mathbb{R}$, defining the RKHS $\mathcal{H}$, is symmetric and positive-definite. \\
\textbf{Reproducing property for derivatives:} The kernel satisfies $K\in C^2(\mathbb{R}^d\times\mathbb{R}^d)$ \citet{supportvectormachines}, so that for all $f\in\mathcal{H}$ and all $i=1,\ldots,d$,
\begin{align}\label{Reproducing Property for Derivatives}
\partial_{x_i} f(x) = \langle f, \partial_{x_i} K(x,\cdot) \rangle_\mathcal{H}.
\end{align}
\end{assumptionsA}

\begin{rem}\label{remark: absolutely continuous initial density}
~The approximation of the score in \eqref{explicit argmin} becomes
\begin{equation}\label{kernel regression particle system}
\begin{split}
\partial_{x_i} \log \rho(x) &\approx r_i(x) + \frac{1}{P} \sum_{p=1}^P a_{p\,i} K(x,x_p)
\\
a_{p\,i} &:=  - \sum_{k=1}^P \left( (\mathbf{K}^2+ \lambda \mathbf{K})^{-1} \right)_{pk} \sum_{l=1}^P \left( \partial_{x_i} K(x_l, x_k)  + K(x_l,x_k) r_i(x_l) \right),
\end{split}
\end{equation}
for any empirical measure
\begin{align*}
\bar{\mu}^\mathbf{x} = \frac{1}{P}\sum_{p=1}^P \delta(\cdot - x_p)
\end{align*}
of particles $\mathbf{x}=(x_1,\ldots,x_P)^T$, using the matrix $\mathbf{K}_{ij} := K(x_i,x_j)$. For the particle system, the parameter $\lambda$ is normalized by the system size $P$, solving $r_i + \arg\min_{\phi\in\mathcal{H}}  \left( J_i^r[\phi,\mu_t] + \lambda/P \|\phi\|_{\mathcal{H}}^2 \right)$, $i=1,\ldots,d$, instead. See Section 5.2 in the work by \citet{GradientEstimation} for further details. Thus, $\lambda>0$ serves as a Tikhonov regularization parameter, improving the stability of the minimizer from \eqref{explicit argmin}.
\end{rem}


The score in the mean-field ODE \eqref{sec. exp: mean field ODE true gradlog} is approximated using Hyvärinen's score matching \citep{ScoreMatching}, resulting in the mean-field ODE \eqref{mean-field ODE approximation}. To quantify the resulting error, the approximation $r + \phi_{\mu_t^\infty}^*$ of $\nabla \log \rho_t$ in a general approximation space $\mathcal{V}$ is assumed to satisfy
\begin{align}\label{approximation assumption}
\exists  \epsilon(t) \geq 0~\forall t\in[0,T]~:~ \| \left(\nabla \log \rho_t - \left( r(x) +  \phi_{\mu_t^\infty}^* \right) \right)_i\|_{L_{\mu_t}^2(\mathbb{R}^d)} \leq \epsilon(t),
\end{align}
for all $i=1,\ldots,d$, with the weighted $L^2$-norm $\|f\|_{L_{\mu_t}^2(\mathbb{R}^d)} = \left(\int_{\mathbb{R}^d} f(x)^2  \ \mu_t(dx)\right)^{1/2}$. \Cref{approximation assumption} is an abstract approximation property of the score surrogate in the mean-field ODE \eqref{mean-field ODE approximation}. 

\begin{thm}{~(Vanishing Wasserstein stability under a score residual)}\label{thm:vanishing-wasserstein-stability}
\\Let \(T>0\), let \(\mu_0=\mu_0^\infty\in \Ptwo(\mathbb{R}^d)\). Let \(\hat{X}(t),\hat{X}^\infty(t)\) be absolutely continuous processes on a common probability space such that
\[
\hat{X}(0)=\hat{X}^\infty(0) \quad\text{a.s.}, \qquad \Law\left(\hat{X}(t)\right)=\mu_t,\qquad \Law\left(\hat{X}^\infty(t)\right)=\mu_t^\infty,
\]
and, for almost every (a.e.) \(t\in[0,T]\),
\[
\frac{d\hat{X}(t)}{dt}=\hat{b}(\hat{X}(t),\mu_t),\qquad
\frac{d\hat{X}^\infty(t)}{dt}=b^\mathcal{H}(\hat{X}^\infty(t),\mu_t^\infty),
\]
with $\hat{b}$ and $b^\mathcal{H}$ defined in \eqref{mean-field-ode-velocity}, and \eqref{mean-field-ode-score-veclocity} respectively. This work assumes that
\begin{align*}
\exists  \epsilon(t) \geq 0~\forall t\in[0,T]~:~ \| \left(\nabla \log \rho_t - \left( r +  \phi_{\mu_t^\infty}^* \right) \right)_i\|_{L_{\mu_t}^2(\mathbb{R}^d)} \leq \epsilon(t),
\end{align*}
for all components $i=1,\ldots,d$, where $\mu_t(dx) = \rho_t(x) \ dx$ and $\mu_t^\infty(dx) = \rho_t^\infty(x) \ dx$ and assumes that, in addition to \nameref{ass:A} 
\begin{align}\label{lipschitz of velocity general}
\|b^\mathcal{H}(x,\mu) - b^\mathcal{H}(y,\nu)\|_2 \leq C(\|x-y\|_2 + W_2(\mu,\nu)),
\end{align}
for all $x,y\in\mathbb{R}^d$ and all $\mu,\nu\in\Ptwo(\mathbb{R}^d)$.
Then, for every \(t\in[0,T]\),
\[
W_2(\mu_t,\mu_t^\infty)^2
\le
B^2d\int_0^t e^{(4C+2\tilde{C}+1)(t-s)}\,\epsilon(s)^2 \ ds,
\]
with $\tilde{C} = C\left( 3 + \sqrt{2Bd}\left( \sqrt{B}\frac{B^2+B}{\lambda} + B \right) \right)$. In particular, if \(\epsilon\equiv 0\), then
\[
W_2(\mu_t,\mu_t^\infty)=0 \qquad\text{for all }t\in[0,T].
\]
Hence, \(\mu_t=\mu_t^\infty\) for all \(t\). If $g$ is globally Lipschitz, then
\begin{align*}
&\Bigg| \mathbb{E}\left[g\left(\hat{X}(T)\right) ~\Bigg|~  \mu_0 = \mu \right] - \mathbb{E}\left[g\left( \hat{X}^{\infty} (T) \right) ~\Bigg|~ \hat{X}^{\infty} (0) \sim \mu\right] \Bigg|
\\
 \leq  &\text{Lip}(g) B \sqrt{d} \left( \int_0^T e^{(4C+2\tilde{C}+1)(T-s)}\,\epsilon(s)^2 \ ds \right)^{1/2}.
\end{align*}
\end{thm}
\begin{proof}
See Appendix~\ref{sec:Wasserstein}.
\end{proof}

\begin{rem}
~Since $g$ is globally Lipschitz by \nameref{ass:A}, differences of expectations are controlled by the $W_2$-distance. Let $\pi\in\Pi(\mu_t,\mu_t^\infty)$ be any coupling. Then, by the Cauchy--Schwarz inequality,
\begin{align*}
&\Bigg| \int_{\mathbb{R}^d} g \ d\mu_t - \int_{\mathbb{R}^d} g \ d\mu_t^\infty \Bigg|
= \Bigg| \int_{\mathbb{R}^d\times\mathbb{R}^d} \left( g(x) - g(y) \right) \pi(dx,dy) \Bigg|
\\
\leq &\text{Lip}(g) \left( \int_{\mathbb{R}^d\times\mathbb{R}^d} \|x-y\|_2^2 \ \pi(dx,dy) \right)^{1/2}.
\end{align*}
All integrals are finite, because $g$ grows at most linearly and $\mu_t,\mu_t^\infty\in\mathcal{P}_2(\mathbb{R}^d)$. Taking the infimum over all couplings $\pi$ yields
\begin{align*}
&\Bigg| \mathbb{E}\left[g\left(\hat{X}(t)\right) ~\Bigg|~  \mu_0 = \mu \right] - \mathbb{E}\left[g\left( \hat{X}^{\infty} (t) \right) ~\Bigg|~ \hat{X}^{\infty} (0) \sim \mu\right] \Bigg| \leq \text{Lip}(g)\, W_1(\mu_t,\mu_t^\infty)
\\
 \leq &\text{Lip}(g)\, W_2(\mu_t,\mu_t^\infty),
\end{align*}
which, combined with \Cref{thm:vanishing-wasserstein-stability}, gives the stated bound on the difference of expectations.
\end{rem}
Direct access to $\mu_t^\infty$ does not exist in the mean-field ODE \eqref{mean-field ODE approximation}; thus, $\mu_t^\infty$ is approximated using an empirical measure, yielding the corresponding particle system. For a set of initial points $\boldsymbol{\zeta} = (\zeta_1,\ldots,\zeta_P)^T$ at time $t\geq 0$ for $p=1,\ldots,P$, this work defines $\mathbf{\bar{X}}^{P|t|\boldsymbol{\zeta}}(s) = (\bar{X}^{P|t|\zeta_1}(s),\ldots,\bar{X}^{P|t|\zeta_P}(s))^T$ with $\bar{X}^{P|t|\zeta_p}:[t,T]\times\Omega\to\mathbb{R}^d$, $p=1,\ldots,P$. The particle system solves the following:
\begin{equation}\label{interacting particle system}
\begin{split}
\frac{d\bar{X}^{P|t|\zeta_p}(s)}{ds} =& b^\mathcal{H}\left(\bar{X}^{P|t|\zeta_p}(s),\bar{\mu}^{\mathbf{\bar{X}}^{P|t|\boldsymbol{\zeta}}(s)}\right),\quad s\in (t,T],~ p=1,\ldots,P,\\
\bar{X}^{P|t|\zeta_p}(t) = \zeta_p,
\end{split}
\end{equation}
using the empirical measure
\begin{align}\label{empirical measure interacting particle system}
\bar{\mu}^{\mathbf{\bar{X}}^{P|t|\boldsymbol{\zeta}}(s)} = \frac{1}{P}\sum_{p=1}^P \delta\left(\cdot - \bar{X}^{P|t|\zeta_p}(s)\right).
\end{align}

\begin{figure}[htbp]
  \centering
  \begin{subfigure}[t]{0.45\textwidth}
    \centering
	\includegraphics[width=\textwidth]{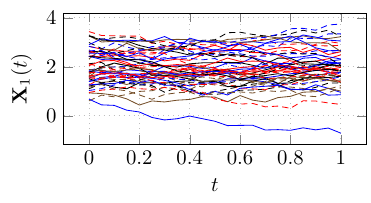}
\caption{Stochastic differential equation sample paths in the first dimension $\mathbf{X}_1(t)$.}

  \end{subfigure}%
  \hfill
  \begin{subfigure}[t]{0.45\textwidth}
\centering

\includegraphics[width=\textwidth]{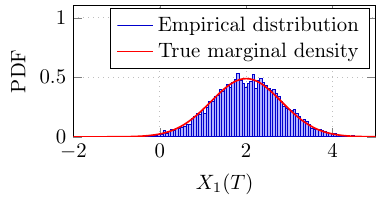}
\caption{Stochastic differential equation histogram and true marginal density in the first dimension $\mathbf{X}_1(T)$.}

\end{subfigure}

  \vspace{0.5em} 

  \begin{subfigure}[t]{0.45\textwidth}
    \centering

\includegraphics[width=\textwidth]{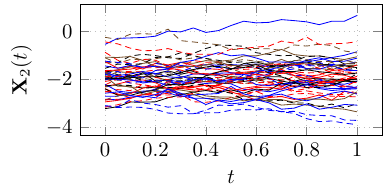}
\caption{Stochastic differential equation sample paths in the second dimension $\mathbf{X}_2(t)$.}

  \end{subfigure}%
  \hfill
  \begin{subfigure}[t]{0.45\textwidth}
\centering

\includegraphics[width=\textwidth]{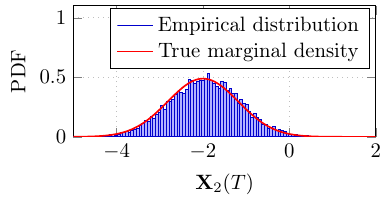}
\caption{Stochastic differential equation histogram and true marginal density in the second dimension $\mathbf{X}_2(T)$.}

\end{subfigure}

  \caption{Fifty sample paths and histograms at $T=1$ of the 2D Gaussian model in Appendix~\ref{sec:Gaussian Model} using an interacting particle system of size $P=10000$. The interacting particle systems were created for the original McKean--Vlasov stochastic differential equation (MV-SDE).}
  \label{fig:Sample path histogram SDE}
\end{figure}

\begin{figure}[htbp]

  \begin{subfigure}[t]{0.45\textwidth}
    \centering

\includegraphics[width=\textwidth]{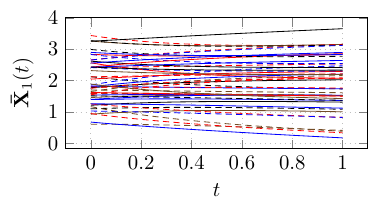}
\caption{Ordinary differential equation sample paths in the first dimension $\mathbf{\bar{X}}_1(t)$.}

  \end{subfigure}%
  \hfill
  \begin{subfigure}[t]{0.45\textwidth}
    \centering
    
\includegraphics[width=\textwidth]{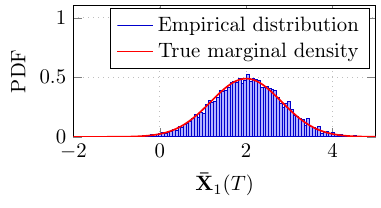}
\caption{Ordinary differential equation histogram and true marginal density in the first dimension $\mathbf{\bar{X}}_1(T)$.}

\end{subfigure}

  \vspace{0.5em} 

  \begin{subfigure}[t]{0.45\textwidth}
    \centering

\includegraphics[width=\textwidth]{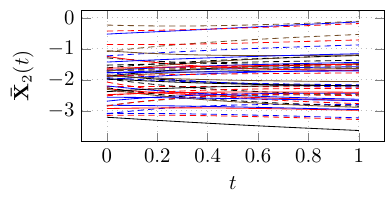}
\caption{Ordinary differential equation sample paths in the second dimension $\mathbf{\bar{X}}_2(t)$.}

  \end{subfigure}%
  \hfill
  \begin{subfigure}[t]{0.45\textwidth}
\centering

\includegraphics[width=\textwidth]{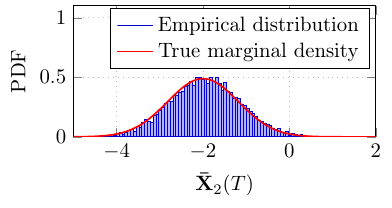}
\caption{Ordinary differential equation histogram and true marginal density in the second dimension $\mathbf{\bar{X}}_2(T)$.}

\end{subfigure}

  \caption{Fifty sample paths and histograms at $T=1$ of the 2D Gaussian model in Appendix~\ref{sec:Gaussian Model} using an interacting particle system of size $P=10000$. The interacting particle systems were created for the corresponding mean-field ordinary differential equation (ODE).}
  \label{fig:Sample path histogram ODE}
\end{figure}

Figures~\ref{fig:Sample path histogram SDE} and \ref{fig:Sample path histogram ODE} demonstrate the effect of this approach, plotting paths approximating an MV-SDE and its corresponding mean-field ODE, and their histograms at time $t=T$. Both systems were simulated with the same starting points. Figures~\ref{fig:Sample path histogram SDE} and \ref{fig:Sample path histogram ODE} correspond to the two-dimensional (2D) Gaussian model defined in \eqref{app Gaussian model} in Appendix~\ref{sec:Gaussian Model}, using
\begin{align*}
m = \begin{pmatrix}
1\\
-1
\end{pmatrix},\quad \Gamma = C_{\text{diff}} = 0.5 \begin{pmatrix}
1 & -0.5\\
-0.5 & 1
\end{pmatrix}
\end{align*} and Gaussian kernels for the score approximation.

This work aims to improve the convergence rates of the interacting particle system \eqref{interacting particle system} to the mean-field ODE \eqref{mean-field ODE approximation} by using rQMC sampling for the initial measure. All randomness in the interacting particle system \eqref{interacting particle system} comes from the initial measure $\mu$. Therefore, the aim is to be efficient when sampling $P$ starting points for the $P$ interacting particles. Although plain MC sampling is possible, rQMC sampling is employed to exploit the underlying smoothness of the problem to improve convergence rates.

\subsection{Randomized quasi-Monte Carlo sampling}
\label{sssec:Randomized Quasi-Monte Carlo Sampling}

\citet{a46f288a95ff44a589103486f365ea79}, \citet{Caflisch_1998}, \citet{randomnumbergeneration} and \citet{dick2010digital} have provided thorough reviews of rQMC sampling. QMC methods serve to approximate integrals over the unit cube. To recast the problem, let $\boldsymbol{\mu}_t^P$, where $t\in[0,T]$, denote the joint probability distribution of $\bar{\mathbf{X}}^P(t)$, the solution of the interacting particle system \eqref{interacting particle system}, and let $\boldsymbol{\mathcal{T}}_t^P:[0,1]^{Pd}\to\mathbb{R}^{Pd}$ be a transport, such that
\begin{align}\label{transport interacting particle system}
\boldsymbol{\mathcal{T}}_t^P \# \boldsymbol{\mathcal{U}}[0,1]^{Pd} = \boldsymbol{\mu}_t^P,
\end{align}
where $\mathcal{U}[0,1]^{Pd}$ represents the uniform distribution on the unit cube $[0,1]^{Pd}$, and such that
\begin{align*}
\boldsymbol{\mathcal{T}}_t^P( \boldsymbol{\xi} ) = \bar{\mathbf{X}}^P(t),
\end{align*}
for $\boldsymbol{\xi}\in [0,1]^{Pd}$ and $\bar{\mathbf{X}}^P(t) = (\bar{X}_1^P(t),\ldots,\bar{X}_P^P(t))^T$ being a solution to the interacting particle system \eqref{interacting particle system}. This work defines the following:
\begin{align*}
\mathbf{g}_P:\mathbb{R}^{Pd} \to \mathbb{R},~ \mathbf{X} \mapsto \frac{1}{P}\sum_{p=1}^P g(\mathbf{X}_p).
\end{align*}
As $\mathbf{g}_P$ averages $P$ identically distributed marginals, its expectation equals the single-marginal expectation under the corresponding law. Therefore,
\begin{equation}\label{integration unitcube}
\begin{split}
\mathbb{E}\left[  \frac{1}{P}\sum_{p=1}^P g(\bar{X}_p^P(T))  \right] &= \int_{\mathbb{R}^{Pd}} \mathbf{g}_P(\mathbf{x}) \ \boldsymbol{\mu}_T^P(\mathbf{dx})\\
&= \int_{[0,1]^{Pd}} \mathbf{g}_P(\boldsymbol{\mathcal{T}}_T^P( \boldsymbol{\xi} )) \ \boldsymbol{\mathcal{U}}[0,1]^{Pd}(\mathbf{d}\boldsymbol{\xi}).
\end{split}
\end{equation}
Typically, the Koksma--Hlawka inequality \citep{dick2010digital} quantifies the approximation of the integral \eqref{integration unitcube} with low-discrepancy points, which requires the total variation of $\mathbf{g}_P(\boldsymbol{\mathcal{T}}_t^P( \boldsymbol{\xi} ))$ to be bounded in the sense of Hardy and Krause. In this case, the Hardy--Krause assumption is unreasonable, because the integral over $\mathbb{R}^{Pd}$ is converted into an integral over $[0,1]^{Pd}$. The framework established by Owen's growth condition \citep{OwenGrowth} extends the Koksma--Hlawka inequality to integrands with infinite total variation but at most polynomial growth up to the boundary. 
This work defines $P$ independently evolved solutions $\hat{\mathbf{X}}^{\infty|t|\boldsymbol{\zeta}}(s)=(\hat{X}^{\infty|t|\zeta_1}(s),\ldots,\hat{X}^{\infty|t|\zeta_P}(s))^T$, constructed analogously to the interacting particle system \eqref{interacting particle system}. Each particle $\hat{X}^{\infty|t|\zeta_p}:[t,T] \times\Omega\to\mathbb{R}^d$, $p=1,\ldots,P$ evolves independently according to
\begin{equation}\label{rQMC: independent limit system transport}
\begin{split}
\frac{d\hat{X}^{\infty|t|\zeta_p}(s)}{ds} &=  b^\mathcal{H}\left( \hat{X}^{\infty|t|\zeta_p}(s),  \mu_s^\infty \right),\quad s\in (t,T], p = 1,\ldots,P\\
\hat{X}^{\infty|t|\zeta_p}(t) &= \zeta_p.
\end{split}
\end{equation}
The corresponding Kolmogorov backward equation (KBE) for $u(t,\boldsymbol{\zeta}) = \mathbf{g}_P\left( \bar{\mathbf{X}}^{P|t|\boldsymbol{\zeta}}(T) \right)$ is:
\begin{equation}\label{rQMC: constant cost to go}
 \begin{split}
 \frac{\partial u}{\partial t}(t,\boldsymbol{\zeta})  + \sum_{p=1}^{P}   b^\mathcal{H}(\boldsymbol{\zeta}_p,\bar{\mu}^{\boldsymbol{\zeta}}) \cdot \nabla_{\boldsymbol{\zeta}_p} u(t,\boldsymbol{\zeta})  &= 0,\quad (t,\boldsymbol{\zeta})\in [0,T) \times \mathbb{R}^{Pd}\\
u(T,\boldsymbol{\zeta}) &= \mathbf{g}_P(\boldsymbol{\zeta}),
 \end{split}
 \end{equation}
The term $\bar{\mu}^{\boldsymbol{\zeta}}$ represents the empirical measure given by $\boldsymbol{\zeta}$
\begin{align*}
\bar{\mu}^{\boldsymbol{\zeta}} = \frac{1}{P}\sum_{p=1}^P \delta(\cdot - \zeta_p).
\end{align*}
To quantify and bound the rQMC weak-error, additional assumptions (\nameref{ass:B}) are required.

\begin{assumptionsB}\label{ass:B}
\textbf{Linearization:} This work assumes that $b(x,\mu)$, $D(x,\mu)$ and $\nabla \cdot D(x,\mu)$ have a linear and bounded linearization, yielding the following kernel functions:
\begin{align*}
A_{b}&:\mathbb{R}^d\times\mathbb{R}^d\times \mathcal{P}(\mathbb{R}^d) \to \mathbb{R}^d,~(x,y;\mu) \mapsto A_b(x,y;\mu),
\\
A_{D}&:\mathbb{R}^d\times\mathbb{R}^d\times \mathcal{P}(\mathbb{R}^{d}) \to \mathbb{R}^{d\times d},~(x,y;\mu) \mapsto A_D(x,y;\mu),
\\
A_{\nabla\cdot D}&:\mathbb{R}^d\times\mathbb{R}^d\times \mathcal{P}(\mathbb{R}^d) \to \mathbb{R}^d,~(x,y;\mu) \mapsto A_{\nabla\cdot D}(x,y;\mu),
\end{align*}
which can be employed to linearize the following, for any $\mu,\nu\in\mathcal{P}(\mathbb{R}^d)$, and any $x\in\mathbb{R}^d$
\begin{equation}\label{linearization drift diffusion}
\begin{split}
b(x,\mu) - b(x,\nu)&= \int_{\mathbb{R}^d} A_{b}(x,y;\mu) \ (\mu - \nu)(dy) + R_{b}(x;\mu,\nu),
\\
D(x,\mu) - D(x,\nu) &= \int_{\mathbb{R}^d} A_{D}(x,y;\mu) \ (\mu - \nu)(dy) + R_{D}(x;\mu,\nu),
\\
\nabla \cdot D(x,\mu) - \nabla \cdot D(x,\nu) &= \int_{\mathbb{R}^d} A_{\nabla \cdot D}(x,y;\mu) \ (\mu - \nu)(dy) + R_{\nabla\cdot D}(x;\mu,\nu),
\end{split}
\end{equation}
where the remainders satisfy
\begin{align*}
& \mathbb{E}_{\boldsymbol{\xi}} \Bigg[ \sup_{t\in[0,T]} \sup_{x\in\mathbb{R}^{d}} \max \Bigg\{ \Big\|R_{b}\left(x;\mu_t^\infty,\bar{\mu}^{\hat{\mathbf{X}}^{\infty|0|\boldsymbol{\mathcal{T}}(\boldsymbol{\xi}(\mathbf{x}))}(t)}\right) \Big\|_2, \\ 
&\hspace*{1cm} \Big\|R_{D}\left(x;\mu_t^\infty,\bar{\mu}^{\hat{\mathbf{X}}^{\infty|0|\boldsymbol{\mathcal{T}}(\boldsymbol{\xi}(\mathbf{x}))}(t)}\right)\Big\|_2,  \Big\|R_{\nabla\cdot D}\left(x;\mu_t^\infty,\bar{\mu}^{\hat{\mathbf{X}}^{\infty|0|\boldsymbol{\mathcal{T}}(\boldsymbol{\xi}(\mathbf{x}))}(t)}\right)\Big\|_2 \Bigg\}\Bigg] \\
\leq & C_R P^{-\alpha_R} ,
\end{align*}
where $\hat{\mathbf{X}}^{\infty|0|\boldsymbol{\mathcal{T}}(\boldsymbol{\xi}(\mathbf{x}))}(t)=(\hat{X}^{\infty|0|\mathcal{T}_0(\xi(x_1))}(t),\ldots,\hat{X}^{\infty|0|\mathcal{T}_0(\xi(x_P))}(t))^T$ solves the mean-field ODE \eqref{rQMC: independent limit system transport} and $\alpha_R>1$. Moreover, $\boldsymbol{\mathcal{T}}(\boldsymbol{\xi}(\mathbf{x})) = (\mathcal{T}_0(\xi(x_1)),\ldots,\mathcal{T}_0(\xi(x_P)))^T$ represents a set of initial conditions, where $\xi$ denotes one randomization of the low-discrepancy points $x_1,\ldots,x_P$ and $\mathcal{T}_0:[0,1]^d\to\mathbb{R}^d$ indicates a transport such that $(\mathcal{T}_0)\#\mathcal{U}[0,1]^d = \mu$. $C_R$ is a constant independent of $t$, $x$, and $P$.

\textbf{Transport:} For every $t\in[0,T]$, a measurable transport map $\mathcal{T}_t:[0,1]^d\to\mathbb{R}^d$ is assumed to exist such that
\begin{align}\label{rQMC: transport}
\mathcal{T}_t\#\mathcal{U}[0,1]^d = \mu_t^\infty
\end{align}
and
\begin{align}\label{rQMC: transport particle}
\mathcal{T}_t(\xi(x_p)) = \hat{X}^{\infty|0|\zeta_p}(t)
\end{align}
where $\hat{X}^{\infty|0|\zeta_p}(t)$ denotes solutions to the mean-field ODE \eqref{rQMC: independent limit system transport} with the initial conditions $\zeta_p\sim\mu$, $p=1,\ldots,P$.\\
\textbf{Uniform randomized quadrature bound:}
Define
  \begin{align*}
  \mathcal{A}_t^b &= \left\lbrace z \mapsto [A_b(x,z;\mu_t^\infty)]_{i} \mid x\in\mathbb{R}^d,~1\leq i \leq d \right\rbrace
  \\
  \mathcal{A}_t^{\text{div}} &= \left\lbrace z \mapsto [A_{\nabla\cdot D}(x,z;\mu_t^\infty)]_{i} \mid x\in\mathbb{R}^d,~1\leq i \leq d \right\rbrace
  \\
  \mathcal{A}_t^D &= \left\lbrace z \mapsto [A_D(x,z;\mu_t^\infty)]_{ij} \mid x\in\mathbb{R}^d,~1\leq i,j \leq d \right\rbrace
  \\
  \mathcal{S}_t^r &= \left\lbrace z \mapsto \left[ (\lambda I + \mathcal{C}_{\mu_t^\infty})^{-1}K(x,\cdot) \right](z) r_i(z) \mid x\in\mathbb{R}^d,~1\leq i \leq d \right\rbrace
  \\
    \mathcal{S}_t^{\partial} &= \left\lbrace z \mapsto  \left[ \partial
    _{x_i} \left( (\lambda I + \mathcal{C}_{\mu_t^\infty})^{-1} K(x,\cdot) \right) \right](z) \mid x\in\mathbb{R}^d,~1\leq i \leq d \right\rbrace
  \end{align*}
  and
    \begin{align*}
  \mathcal{F}= \{ z\mapsto h_1(z) h_2(z) \mid h_1\in\mathcal{H}, h_2\in\mathcal{H}, \|h_1\|_\mathcal{H}\leq 1, \|h_2\|_\mathcal{H}\leq 1 \}
  \end{align*} as in \eqref{new function space} from Appendix \ref{sec:Score Estimator}. Then let
  \begin{align*}
  \mathcal{G}_t = \mathcal{A}_t^b \cup \mathcal{A}_t^{\text{div}} \cup \mathcal{A}_t^D \cup \mathcal{S}_t^r \cup \mathcal{S}_t^{\partial} \cup \mathcal{F}.
  \end{align*}
  This work assumes that the supremum below is measurable and that there exists a constant $C_{\mathcal{G}}<\infty$, independent of $P$, such that
  \begin{align}\label{measure convergence}
&\mathbb{E}_{\boldsymbol{\xi}} \left[ \sup_{t\in[0,T]} \sup_{f\in \mathcal{G}_t}\Bigg|\int_{\mathbb{R}^d} f(z) \ (\mu_t^\infty - \bar{\mu}^{\hat{\mathbf{X}}^{\infty|0|\boldsymbol{\mathcal{T}}(\boldsymbol{\xi}(\mathbf{x}))}(t)})(dz) \Bigg|\right] \leq C_{\mathcal{G}} P^{-\alpha_Q},
\end{align}
where $\alpha_Q\in (0,1]$.

\textbf{Kolmogorov backward equation:} This work assumes that, for each $P$, the flow generated by the vector field of the interacting particle system \eqref{interacting particle system} is continuously differentiable with respect to the initial condition $\boldsymbol{\zeta}$, and that the value function $u(t,\boldsymbol{\zeta}) = \mathbf{g}_P(\mathbf{\bar{X}}^{P|t|\boldsymbol{\zeta}}(T))$ is a classical solution of the KBE \eqref{rQMC: constant cost to go}. Furthermore, for all particle systems of size $P\in\mathbb{N}$, a constant $C_u$ independent of $P$ and $t$ exists, such that
\begin{align}\label{bound nabla u}
\sup_{t\in[0,T]} \sup_{P\in\mathbb{N}} \sup_{\boldsymbol{\zeta}\in\mathbb{R}^{Pd}} \sum_{p=1}^P \Big\| \nabla_{\boldsymbol{\zeta}_p} u\left(t,\hat{\mathbf{X}}^{\infty|0|\boldsymbol{\zeta}}(t)\right) \Big\|_2 \leq C_u.
\end{align}
\end{assumptionsB}

\begin{rem}
~The existence of the transport \eqref{rQMC: transport} follows, from the Knothe--Rosenblatt rearrangement; see \citet{transport}. The convergence shown in \Cref{thm:rQMC-transport} relies on \eqref{measure convergence}, which is not derived from standard rQMC theory: The class $\mathcal{G}_t$ depends on the mean-field flow and the score surrogate. Establishing verifiable sufficient conditions is left for future work. Measurability for \eqref{measure convergence} may, for example, be ensured by requiring the transported class
\begin{align*}
\{f\circ\mathcal{T}_t \mid t\in[0,T],~ f\in\mathcal{G}_t\}
\end{align*}
to admit a countable subclass determining the quadrature-error supremum.
\end{rem}

\begin{thm}{~(rQMC transport)}\label{thm:rQMC-transport}
\\ The set of initial conditions is denoted by $\boldsymbol{\mathcal{T}}(\boldsymbol{\xi}(\mathbf{x})) = (\mathcal{T}_0(\xi(x_1)),\ldots,\mathcal{T}_0(\xi(x_P)))^T$, where $\xi$ denotes one randomization of the low-discrepancy points $x_1,\ldots,x_P$ and $\mathcal{T}_0:[0,1]^d\to\mathbb{R}^d$ represents a transport, such that $(\mathcal{T}_0)\#\mathcal{U}[0,1]^d = \mu$. The interacting particle system $\bar{\mathbf{X}}^{P|0|\boldsymbol{\mathcal{T}}(\boldsymbol{\xi}(\mathbf{x}))}(t)=(\bar{X}^{P|0|\mathcal{T}_0(\xi(x_1))}(t),\ldots,\bar{X}^{P|0|\mathcal{T}_0(\xi(x_P))}(t))^T$ is defined based on \eqref{interacting particle system} and the independently evolved mean-field ODEs $\hat{\mathbf{X}}^{\infty|0|\boldsymbol{\mathcal{T}}(\boldsymbol{\xi}(\mathbf{x}))}(t)=(\hat{X}^{\infty|0|\mathcal{T}_0(\xi(x_1))}(t),\ldots,\hat{X}^{\infty|0|\mathcal{T}_0(\xi(x_P))}(t))^T$ are defined in \eqref{rQMC: independent limit system transport}. Moreover, define $\mathbf{g}_P:\mathbb{R}^{Pd}\to\mathbb{R}, (x_1,\ldots,x_P) \mapsto \frac{1}{P}\sum_{p=1}^P g(x_p)$ for the function $g:\mathbb{R}^d\to\mathbb{R}$. Then, under \nameref{ass:A} and \nameref{ass:B},
\begin{align}\label{rQMC convergence rate}
\Big|\mathbb{E}_{\boldsymbol{\xi}}\left[  \mathbf{g}_P(\bar{\mathbf{X}}^{P|0|\boldsymbol{\mathcal{T}}(\boldsymbol{\xi}(\mathbf{x}))}(T))\right]  - \mathbb{E}_{\boldsymbol{\xi}}\left[ \mathbf{g}_P(\hat{\mathbf{X}}^{\infty|0|\boldsymbol{\mathcal{T}}(\boldsymbol{\xi}(\mathbf{x}))}(T))  \right]\Big| = \mathcal{O}\left( P^{-\alpha}\right),
\end{align}
with $\alpha = \min\{\alpha_Q,\alpha_R\}$.
\end{thm}
\begin{proof}
Appendix~\ref{sec:RQMCTransport} presents the proof.
\end{proof}

Finally, the classical Runge--Kutta scheme of order $\gamma = 4$ is employed for the time discretization. The numerical approximation of ODEs is a classic topic \citep{teschl2024ordinary,ODEsI,ODEsII} and this work does not claim that the classical Runge--Kutta scheme is optimal to approximate the interacting particle system \eqref{interacting particle system}. Instead, as a widely known higher-order scheme, it is selected to discretize the ODE. This paper defines the approximations to \Cref{thm:rQMC-transport} with $N$ time steps using an additional superscript as follows:
\begin{align}\label{numerical discretization}
\mathbf{\bar{X}}^{P|N|0|\boldsymbol{\mathcal{T}}(\boldsymbol{\xi}(\mathbf{x}))}(t_n) = (\bar{X}^{P|N|0|\mathcal{T}_0(\xi(x_1))}(t_n),\ldots,\bar{X}^{P|N|0|\mathcal{T}_0(\xi(x_P))}(t_n))^T,\quad n=0,\ldots,N.
\end{align}
Shorthand notation is introduced for the estimator. With a low-discrepancy sequence $\mathbf{x} = (x_1,\ldots,x_P)$ of $P$ points, and one randomization $\boldsymbol{\xi}$, this work defines
\begin{align*}
\boldsymbol{\psi}^{P|N}\left(\boldsymbol{\xi}\left(\mathbf{x}\right)\right) := \frac{1}{P}\sum_{p=1}^P g\left(\bar{X}^{P|N|0|\mathcal{T}_0(\xi(x_p))}(t_N)\right).
\end{align*}
For $M$ independent randomizations $\boldsymbol{\xi}^{(m)}$, $m = 1,\ldots,M$, this approach obtains the rQMC estimator of $P$ particles, $N$ time steps, and $M$ randomizations
\begin{align}\label{estimator definition}
\mathcal{A}_{\text{rQMC}}^{P|N|M}\left(\mathbf{x}\right) = \frac{1}{M}\sum_{m=1}^M \boldsymbol{\psi}^{P|N}\left(\boldsymbol{\xi}^{(m)}\left(\mathbf{x}\right)\right).
\end{align}
The overall error of the estimator is bounded as follows:
\begin{align}
&\Bigg| \mathbb{E}\left[g\left(X(T)\right) \mid \mu_0 = \mu \right ] - \mathcal{A}_{\text{rQMC}}^{P|N|M}(\mathbf{x})\Bigg|
\notag\\
\leq &\underbrace{ \Bigg| \mathbb{E}\left[g\left(X(T)\right) \mid \mu_0 = \mu \right] - \mathbb{E}\left[g\left(\hat{X}(T)\right) ~\Bigg|~ \mu_0 = \mu \right] \Bigg|}_{=0}
\notag\\
+& \underbrace{\Bigg| \mathbb{E}\left[g\left(\hat{X}(T)\right) ~\Bigg|~  \mu_0 = \mu \right] - \mathbb{E}\left[\mathbf{g}_P\left( \mathbf{\hat{X}}^{\infty|0|\boldsymbol{\mathcal{T}}(\boldsymbol{\xi}(\mathbf{x}))} (T) \right) ~\Bigg|~ \mathcal{T}_0(\xi(x_p)) \sim \mu, ~p=1,\ldots,P \right] \Bigg|}_{=e_1: \text{ Theorem }\ref{thm:vanishing-wasserstein-stability}}
\notag\\
+& \Bigg| \mathbb{E}\left[\mathbf{g}_P\left( \mathbf{\hat{X}}^{\infty|0|\boldsymbol{\mathcal{T}}(\boldsymbol{\xi}(\mathbf{x}))} (T) \right) ~\Bigg|~ \mathcal{T}_0(\xi(x_p)) \sim \mu, ~p=1,\ldots,P \right]
\notag\\ 
&\underbrace{\hspace{3cm} - \mathbb{E}\left[\mathbf{g}_P\left( \mathbf{\bar{X}}^{P|0|\boldsymbol{\mathcal{T}}(\boldsymbol{\xi}(\mathbf{x}))} (T) \right) ~\Bigg|~ \mathcal{T}_0(\xi(x_p)) \sim \mu, ~p=1,\ldots,P \right] \Bigg|}_{=e_2: \text{ Theorem }\ref{thm:rQMC-transport}}
\notag\\
+& \Bigg| \mathbb{E}\left[\mathbf{g}_P\left( \mathbf{\bar{X}}^{P|0|\boldsymbol{\mathcal{T}}(\boldsymbol{\xi}(\mathbf{x}))} (T) \right) ~\Bigg|~ \mathcal{T}_0(\xi(x_p)) \sim \mu, ~p=1,\ldots,P \right]
\notag\\ 
&\underbrace{\hspace{3cm} - \mathbb{E}\left[\mathbf{g}_P\left( \mathbf{\bar{X}}^{P|N|0|\boldsymbol{\mathcal{T}}(\boldsymbol{\xi}(\mathbf{x}))} (t_N) \right) ~\Bigg|~ \mathcal{T}_0(\xi(x_p)) \sim \mu, ~p=1,\ldots,P \right] \Bigg|}_{=e_3: \text{ Time Discretization} }
\notag\\
+& \underbrace{\Bigg| \mathbb{E}\left[\mathbf{g}_P\left( \mathbf{\bar{X}}^{P|N|0|\boldsymbol{\mathcal{T}}(\boldsymbol{\xi}(\mathbf{x}))} (t_N) \right) ~\Bigg|~ \mathcal{T}_0(\xi(x_p)) \sim \mu, ~p=1,\ldots,P \right]  - \mathcal{A}_{\text{rQMC}}^{P|N|M}(\mathbf{x})\Bigg|}_{=e_4: \text{ Sampling Error} }.\label{all errors}
\end{align}
\Cref{tab:summary of Assumptioons} summarizes the stated assumptions and their applications to proving Theorems \ref{thm:vanishing-wasserstein-stability} and \ref{thm:rQMC-transport}. \Cref{tab:summary of notation} summarizes the notation for all the introduced stochastic processes in this section. 

\begin{table}[htbp]
\begin{center}

\begin{tabular}{c|c} 
Assumptions & Applied in
\\
\hline
\nameref{ass:A}: Regularity 
\\ 
and boundedness & existence, Theorems~\ref{thm:vanishing-wasserstein-stability} and \ref{thm:rQMC-transport}
\\
\nameref{ass:B}: Linearization,   & \Cref{thm:rQMC-transport}
\\
Transport, Quadrature and KBE Bound & 
\end{tabular}
\caption{Assumptions and their applications.}
\label{tab:summary of Assumptioons}
\end{center}
\end{table} 

\begin{table}[htbp]
\begin{center}
\begin{tabular}{c|c} 
Notation & Definition \\
\hline
$X(t)\rule{0pt}{2.4ex}$ & Original MV-SDE \eqref{McKean-Vlasov SDE} at time $t$.
\\
\hline
$\hat{X}(t)\rule{0pt}{2.4ex}$ & Mean-field ODE \eqref{sec. exp: mean field ODE true gradlog} with original score at time $t$.
\\
\hline
$\hat{X}^\infty(t)\rule{0pt}{2.4ex}$ & Mean-field ODE \eqref{mean-field ODE approximation} with approximated score at time $t$.
\\
\hline
$\hat{\mathbf{X}}^{\infty|t|\boldsymbol{\zeta}}(s)\rule{0pt}{2.4ex}$ & $P$ independently evolved trajectories of the mean-field ODE \eqref{rQMC: independent limit system transport} at
\\
& time $s$, using score approximation, with starting points $\boldsymbol{\zeta}$ at time $t\leq s$.
\\
\hline
$\bar{\mathbf{X}}^{P|t|\boldsymbol{\zeta}}(s)\rule{0pt}{2.4ex}$ & Interacting particle system \eqref{interacting particle system} of size $P$ at time $s$,
\\
& with starting points $\boldsymbol{\zeta}$ at time $t\leq s$.
\\
\hline
$\bar{\mathbf{X}}^{P|N|t|\boldsymbol{\zeta}}(t_n)\rule{0pt}{2.4ex}$ & Interacting particle system \eqref{numerical discretization} of size $P$ at time $t_n$, with starting points
\\
&  $\boldsymbol{\zeta}$ at time $t$, using $N$ equidistant time steps on the grid
\\
&  $0 = t_0 < t_1 < \ldots < t_N = T$.
\\
\hline
$\mathcal{A}_{\text{rQMC}}^{P|N|M}\left(\mathbf{x}\right)\rule{0pt}{3.4ex}$ & rQMC estimator \eqref{estimator definition} with $P$ particles, $N$ time steps,
\\
& $M$ randomizations, and the low-discrepancy sequence $\mathbf{x}$.
\end{tabular}
\caption{Notation for the various processes and estimators introduced in this work. A hat denotes the mean-field ordinary differential equation (ODE), a bar is reserved for interacting particle systems. Superscript $\infty$: population score surrogate; $P$: number of particles; $N$: number of time steps; $M$: number of randomizations. Bold letters indicate $Pd$-dimensional quantities. The superscript $t|\boldsymbol{\zeta}$ denotes that the system starts at time $t$ from the initial condition $\boldsymbol{\zeta}$.}
\label{tab:summary of notation}
\end{center}
\end{table}

\section{Computational complexity}
\label{sec:Work Complexity}

\begin{figure}[htbp]
\includegraphics[width=\textwidth]{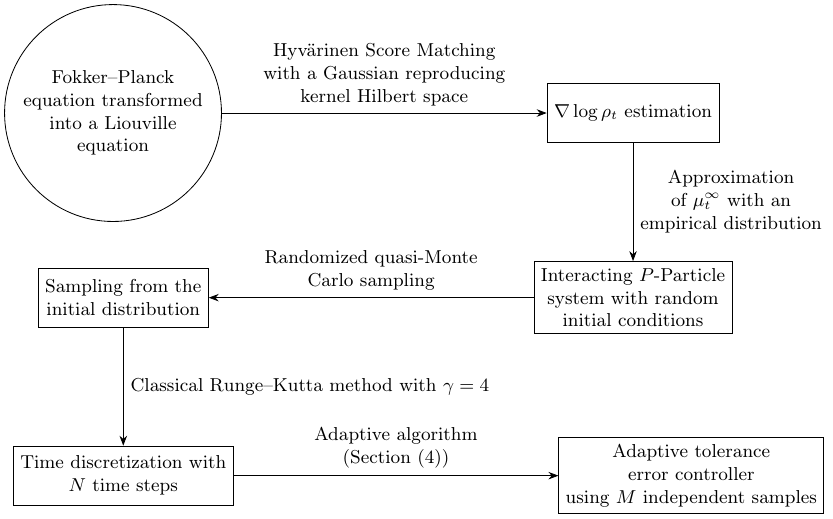}
\caption{Computational architecture of the mean-field ordinary differential equation (ODE) estimator \eqref{estimator definition}.}
\label{fig:architecture}
\end{figure}

Given the error analysis in \Cref{sec:Mean-Field ODE Estimator}, this study now analyzes the work of the estimator.
\begin{itemize}
\item The score is approximated at every time step. This involves solving a linear system of size $P\times P$. The computational work of this operation is $\mathcal{O}\left( P^{\delta_r}\right)$, with rate $\delta_r>0$.
\item The velocity of the interacting mean-field ODE is evaluated for all $P$ trajectories at all time steps, resulting in a computational work of $\mathcal{O}(P^{\delta_v})$, with rate $\delta_v>0$.
\item A Runge--Kutta scheme with $s$ stages and $N$ time steps is employed, thus having $sN$ evaluations of the score surrogate and velocity functions.
\item $M$ independent sets of initial conditions (either MC or rQMC) are generated.
\end{itemize}
In total, the computational work is given by
\begin{align}\label{work}
\mathcal{W} = \mathcal{O}\left( MNP^\delta\right),
\end{align}
where
\begin{align}\label{delta definition}
\delta := \max\{ \delta_r, \delta_v \}.
\end{align}
The current implementation solves a $P\times P$ linear system for the score approximation at each time step using direct factorization. For the modest $P \leq 5000$, used throughout the numerics in \Cref{sec:Numerical Experiments}, $\delta\approx 2$. Asymptotically, this work assumes $\delta_r=3$ for the estimator. For example, a conjugate-gradient solver could be applied to the regularized normal equations (matrix-vector products at $\mathcal{O}\left(P^2\right)$ each) with effective preconditioning, which often converges in a few to tens of iterations with reduced $\delta_r$. Low-rank methods (Nyström/FALKON-style) as well as fast-kernel summation (e.g., fast multipole method or random Fourier features) could be employed to reduce matrix-vector products toward $\mathcal{O}(P)$, which could shift $\delta_r$ toward $1-1.5$, \citep{rudi2017falkon, 174bfb93628e44f5a008a32e7af7da5a,srinivasan2014preconditioned,DARVE2000195,NIPS2015_d14220ee}. For the naive implementation, $\delta_v=2$, which can be reduced with more sophisticated implementations to $\delta_v = 1$ \citep{MLMC}.

The aim is to bound the error $\Big| \mathbb{E}\left[g\left(X(T)\right) \mid \mu_0 = \mu \right ] - \mathcal{A}_{\text{rQMC}}^{P|N|M}(\mathbf{x})\Big|$ by an absolute tolerance TOL by bounding the individual errors listed in \eqref{all errors}. The projection error $e_1$ is fixed a priori by suitable approximation space selection for the score and is independent of the parameters $P$, $N$, and $M$. This work assumes that $e_1$ is negligible, satisfies $\text{TOL}-e_1>0$ and that $e_1$ is set as a fixed approximation error external to the optimization in $(M,N,P)$. The following bias constraint is imposed:
\begin{align}\label{bias constraint}
e_2 + e_3 \leq \theta (\text{TOL} - e_1),
\end{align}
for an error split parameter $\theta \in (0,1)$. Each rQMC randomization $\xi^{(m)}$, $m=1,\ldots,M$, is independent; hence the CLT applies for sufficiently large $M$. For the statistical error $e_4$, the following is imposed:
\begin{align}\label{stat constraint}
e_4  \approx c_\varepsilon \sqrt{ \text{Var}\left( \mathcal{A}_{\text{rQMC}}^{P|N|M}\left(\mathbf{x} \right) \right) } \leq (1-\theta) (\text{TOL} - e_1),
\end{align}
where $c_\varepsilon$ corresponds to the $\left( \frac{1+ \varepsilon}{2} \right)$-quantile of the standard normal distribution. Hence, the following variance constraint can be imposed:
\begin{align*}
\text{Var}\left( \mathcal{A}_{\text{rQMC}}^{P|N|M}\left(\mathbf{x}\right) \right)  \leq \left( \frac{(1-\theta)(\text{TOL} - e_1)}{c_\varepsilon} \right)^2.
\end{align*}
The estimator's bias is assumed to be bounded by 
\begin{align}\label{bias rate}
e_2 + e_3 \leq \frac{B_1}{P^\alpha} + \frac{B_2}{N^\gamma},
\end{align}
with rates $\alpha,\gamma>0$ and positive constants $B_1, B_2 >0$, both independent of $N$ and $P$. The rate $N^{-\gamma}$ and constant $B_2$ are uniform-in-$P$ over the parameter range used in the numerics. This work assumes the following:
\begin{assump}\label{Variance}
~This work assumes the existence of a $V$ independent of $P$ and $N$ and a rate $\beta>0$, so that the variance of the estimator satisfies
\begin{align}\label{rQMC: Variance rate}
\text{Var}\left( \frac{1}{P}\sum_{p=1}^P g\left(\bar{X}_p^P(T)\right) \right) \leq V P^{-\beta}.
\end{align}
\end{assump}
with $\beta>0$ denoting the convergence rate of the estimator variance and a positive constant $V>0$, independent of $N$ and $P$. Assumption \ref{Variance} implies the following bound on the estimator's statistical error:
\begin{align}\label{stat error rate}
e_4  \leq c_\varepsilon \sqrt{\frac{V}{M P^\beta}},
\end{align}
For MC sampling, $\beta = 1$ is expected, and by \Cref{Variance}, $\beta = 2$ is expected for rQMC sampling. For a computationally efficient estimator, the following optimization problem is solved:
\begin{align}\label{eq1}
\left\{
\begin{aligned}
\underset{M, N, P \in \mathbb{N}}{\mathrm{arg\,min}} &\quad M N P^\delta\\
\text{s.t. } &\quad \frac{B_1}{P^\alpha} + \frac{B_2}{N^\gamma} \leq \theta (\text{TOL} - e_1)\\
\text{and } &\quad \frac{c_\varepsilon^2 V}{MP^{ \beta}} \leq (1-\theta)^2 (\text{TOL} - e_1)^2.
\end{aligned}
\right.
\end{align}
The minimization in \eqref{eq1} is relaxed to positive real values of $M$, $N$, and $P$, resulting in the following Lagrangian:
\begin{align}\label{Lagrangian ODE}
\mathcal{L}_{\text{ODE}} =& MNP^\delta + \lambda_1 \left( \frac{B_1}{P^\alpha} + \frac{B_2}{N^\gamma} - \theta (\text{TOL} - e_1) \right) + \lambda_2\left( \frac{c_\varepsilon^2 V}{MP^{ \beta}} - (1-\theta)^2 (\text{TOL} - e_1)^2\right),
\end{align}
with Lagrange multipliers $\lambda_1,\lambda_2 \in \mathbb{R}^+$.

\begin{thm}{~(Conditional complexity under stated bias, variance and work models)}\label{complexity: optimal values}
\\This work considers the mean-field ODE estimator $\mathcal{A}_{\text{rQMC}}^{P|N|M}\left(\mathbf{x}\right) $ in \eqref{estimator definition}. Let $\alpha$, $\beta$, $\gamma$, $\delta$, $B_1$, $B_2$, $V$, $c_\varepsilon>0$ be positive constants, defined in the bias model \eqref{bias rate}, the variance model \eqref{stat error rate}, and the work model \eqref{work}, and let $\delta > \beta$. For any tolerance $\text{TOL} > 0$, the Lagrangian $\mathcal{L}_{\text{ODE}}$ \eqref{Lagrangian ODE} has stationary points given by
\begin{equation}
\begin{split} \label{optimal NMP}
M &= \frac{c_\varepsilon^2V}{(1-\theta)^2} \left( \frac{(\delta-\beta)\theta}{B_1(\delta-\beta + \frac{\alpha}{\gamma})} \right)^{\beta/\alpha} (\text{TOL} - e_1)^{\beta/\alpha - 2}  = \mathcal{O}(\text{TOL} ^{\beta/\alpha - 2}),\\
N &= \left( B_2 \left( \frac{\gamma(\delta-\beta)}{\alpha} + 1 \right) \frac{1}{\theta (\text{TOL} - e_1)} \right)^{1/\gamma} = \mathcal{O}(\text{TOL}^{-1/\gamma}), \\
P& = \left( \frac{B_1(\delta-\beta + \frac{\alpha}{\gamma} )}{\delta-\beta} \frac{1}{\theta (\text{TOL} - e_1)} \right)^{1/\alpha}= \mathcal{O}(\text{TOL}^{-1/\alpha}).
\end{split}
\end{equation}
Let $M_0\geq 1$ be fixed minimum number of rQMC randomisations (or iid MC samples). If $\beta \leq 2\alpha$, then 
\begin{align*}
M' = \max\left\lbrace M_0, \left\lceil \frac{c_\varepsilon^2V}{(1-\theta)^2} \left( \frac{(\delta-\beta)\theta}{B_1(\delta-\beta + \frac{\alpha}{\gamma})} \right)^{\beta/\alpha} (\text{TOL} - e_1)^{\beta/\alpha - 2} \right\rceil  \right\rbrace
\end{align*}
and the resulting computational complexity of this estimator is
\begin{align*}
\mathcal{W} = M'NP^\delta =  \mathcal{O}\left( \text{TOL}^{-\left(2  + \frac{\delta - \beta}{\alpha} + \frac{1}{\gamma} \right)} \right).
\end{align*}
For $\beta> 2\alpha$, the variance constraint is asymptotically inactive at the bias-selected $P$; one may take $M' = M_0$ and thus
\begin{align*}
\mathcal{W} = M'NP^\delta =  \mathcal{O}\left( \text{TOL}^{-\left( \frac{\delta }{\alpha} + \frac{1}{\gamma} \right)} \right).
\end{align*}
\end{thm}
\begin{proof}
See Appendix~\ref{sec:Optimal Work Complexity} for details.
\end{proof} 
When $\delta = \beta$, $MP^\beta$ is a common term in the total work \eqref{work} and in the statistical error constraint \eqref{stat constraint}. Thus, the relaxed optimization in $(M,P)$ at a fixed $MP^\beta$ develops a flat direction. Therefore, in practice, $M \geq M_0$ is fixed (e.g., $M_0 = 20$ for the CLT to hold), and $(P,N)$ is set based on the bias and statistical constraints. This heuristic justifies $M$ as a constant and directly leads to \Cref{corollary optimatl}, by the proof of \Cref{complexity: optimal values} with a constant $M$.
\begin{corr}{~(Computational complexity, $\delta=\beta$ and fixed $M$)}\label{corollary optimatl}
\\Consider the same setting as Theorem \ref{complexity: optimal values}. Furthermore, let $\beta = \delta$, and $M$ be a constant. For any tolerance $\text{TOL} > 0$, the computational cost scaling satisfies
\begin{equation}\label{corollary optimal P N}
\begin{split}
P &=  \mathcal{O}\left( \text{TOL}^{-\max\{1/\alpha,2/\beta\} } \right),
\\
N &=  \mathcal{O}\left( \text{TOL}^{-1/\gamma} \right).
\end{split}
\end{equation}
The resulting computational complexity of this estimator is
\begin{align}\label{corollary optimal work rate}
\mathcal{W} = MNP^\delta =  \mathcal{O}\left( \text{TOL}^{-\left(\max\{\beta/\alpha,2\} +\frac{1}{\gamma} \right)} \right).
\end{align}
\end{corr}

\begin{rem}
~In the numerical regime where $\alpha\approx 1$ and $\beta \approx 2$, \eqref{corollary optimal work rate} reveals that the total computational complexity is given by
\begin{align*}
\mathcal{W} = \mathcal{O}\left( \text{TOL}^{-\left( 2 + \frac{1}{\gamma}  \right)} \right).
\end{align*}
In this implementation, the empirically observed work scaling in the range of the tested particle sizes in Section \ref{sec:Numerical Experiments} is consistent with $\delta= 2$, $\alpha =1$, $\beta = 2$, and $\gamma = 4$, yielding a computational complexity of $\mathcal{O}\left(\text{TOL}^{-2.25}\right)$. For asymptotically large $P$, a dense direct factorization of the $P\times P$ kernel system yields $\delta_r=3$; hence, $\delta=3$. With the same $\alpha$, $\beta$, and $\gamma$, this yields a complexity of $\mathcal{O}\left(\text{TOL}^{-3.25}\right)$. Guaranteeing $\delta = 2$, and thus the rate $\mathcal{O}\left(\text{TOL}^{-2.25}\right)$, requires an $\mathcal{O}(P^2)$ kernel solver (e.g., iterative, low-rank, or fast-kernel methods).
\end{rem}

\begin{rem}
~For brevity, \Cref{complexity: optimal values} and \Cref{corollary optimatl} were formulated for the rQMC estimator $\mathcal{A}_{\text{rQMC}}^{P|N|M}\left(\mathbf{x}\right) $ of \eqref{estimator definition}. However, the same rates hold for various sampling methods, such as standard MC, if the assumptions that lead to the minimization problem \eqref{eq1} apply.
\end{rem}

\FloatBarrier
\section{Adaptive algorithm}
\label{sec:Adaptive Algorithm}
This section presents an algorithm that adaptively chooses values of $M$, $N$, and $P$  for a prescribed tolerance, because the constants and rates in \Cref{complexity: optimal values} and \Cref{corollary optimatl} are not known a priori. 

\subsection{Pilot estimation of rates and constants}
\Cref{adaptive alg. 2} requires inputs $B_1$, $B_2$, and $V$ to compute starting values of $M$, $N$, and $P$, using their optimal values in \eqref{optimal NMP} of \Cref{complexity: optimal values}. The algorithm automatically estimates the values of $B_1$, $B_2$, and $V$, while estimating the rates $\alpha$, and $\gamma$, and $\beta$. These initial estimates of $B_1$, $B_2$, and $V$ need not be precise and can be obtained cheaply. In fact, Algorithm~\ref{adaptive alg. 2} can also be adapted to not require $B_1$, $B_2$ and $V$ and instead use any starting values of $M$, $N$, and $P$.

In the practically relevant regime $\delta \approx 2 \approx \beta $, \Cref{complexity: optimal values} is employed only as scaling rules for the starting values rather than as exact optimal parameters for the implemented solver. \Cref{optimal NMP} can be rewritten as follows:
\begin{align*}
M &= C_1(\alpha,\beta,\gamma,\delta) \text{TOL}_{\text{eff}}^{\beta/\alpha - 2},\\
N &= C_2(\alpha,\beta,\gamma,\delta) \text{TOL}_{\text{eff}}^{-1/\gamma},\\
P &= C_3(\alpha,\beta,\gamma,\delta) \text{TOL}_{\text{eff}}^{-1/\alpha}
\end{align*}
with factors $C_i(\alpha,\beta,\gamma,\delta)$, $i=1,2,3$ independent of TOL. Hence, the dependence on $\text{TOL}_{\text{eff}}$ is governed by $\alpha$, $\beta$, and $\gamma$ but not $\delta$. Thus, if the factors $C_i(\alpha,\beta,\gamma,\delta)$, $i=1,2,3$ are evaluated for an inaccurate value of $\delta$, the algorithm's complexity rate does not change, just the constant associated with its runtime. The $\delta$-independence circumvents the observed problem $\delta \approx \beta$, resulting in $0 < N \ll 1$ and $P \gg 1$. In numerical experiments, $\delta = 3$ was set to obtain reasonable starting values for $M$, $N$, and $P$. However, the inaccuracy in the starting values of $M$, $N$, and $P$ also necessitates the adaptive loop in the proposed algorithm.

\subsection{Bias and statistical error estimators}
The estimator's bias is computed using Richardson extrapolation by generating two low-discrepancy sequences. The first, $\mathbf{x}^{2P}$, is of size $2P$, and the second, $\mathbf{x}^{P}$, is of size $P$. With these sequences, $\mathcal{A}_{\text{rQMC}}^{2P|\left\lceil 2^{\alpha/\gamma}N \right\rceil|M}\left(\mathbf{x}^{2P}\right)$ and $\mathcal{A}_{\text{rQMC}}^{P|N|M}\left(\mathbf{x}^P\right)$ are computed. Thus, the bias is estimated as follows:
\begin{align}\label{bias estimation}
\Big| \mathbb{E}\left[g(X(T)) ~\Big|~ \mu_0 = \mu\right] - \mathcal{A}_{\text{rQMC}}^{P|N|M}\left(\mathbf{x}^P\right) \Big| \approx \frac{1}{1-2^{-\alpha}} \Big| \mathcal{A}_{\text{rQMC}}^{2P|\left\lceil 2^{\alpha/\gamma}N \right\rceil|M}\left(\mathbf{x}^{2P}\right) - \mathcal{A}_{\text{rQMC}}^{P|N|M}(\mathbf{x}^P) \Big|.
\end{align} 
The statistical error is estimated using the CLT:
\begin{equation}\label{statistical error estimation}
\begin{split}
\text{stat}_{\text{error}} &= c_\varepsilon\sqrt{\text{Var}\left( \mathcal{A}_{\text{rQMC}}^{P|N|M}\left(\mathbf{x}^P\right) \right)} \quad \text{with}\\
\text{Var}\left( \mathcal{A}_{\text{rQMC}}^{P|N|M}\left(\mathbf{x}^P\right) \right)& \approx \frac{1}{M(M-1)}\sum_{m=1}^M \left( \mathcal{A}_{\text{rQMC}}^{P|N|M}\left(\mathbf{x}^P\right) - \boldsymbol{\psi}^{P|N}\left(\boldsymbol{\xi}^{(m)}\left(\mathbf{x}^P\right) \right) \right)^2.
\end{split}
\end{equation}

\subsection{Adaptive loop}
The primary structure of Algorithm \ref{adaptive alg. 2} comprises two inner loops to monitor the bias constraint \eqref{bias constraint} and statistical error constraint \eqref{stat constraint}. Accurately estimating the projection error $e_1$ requires a reference solution of the original MV-SDE. The proposed adaptive algorithm does not include such a reference. Instead, this study assumes that the score-estimation space is well chosen so that $e_1$ can be neglected, and
\begin{align*}
\text{TOL}_{\text{eff}} := \text{TOL} - e_1
\end{align*}
is the input of the proposed algorithm. The Gaussian MV-SDE in \Cref{sec:Numerical Experiments} confirms that $e_1$ is negligible, and the mean-field ODE result is compared with its analytic solution. The adaptive loop controls $e_2$, $e_3$, and $e_4$ conditional on the chosen score approximation space, and does not detect a non-negligible score-surrogate bias $e_1$. In models without analytic reference values, this bias must be assessed using score-space refinement, alternative score bases, or problem-specific diagnostics.

The first inner loop in Algorithm \ref{adaptive alg. 2} checks if the bias estimator \eqref{bias estimation} satisfies \eqref{bias constraint}, whereas the second enforces the statistical error constraint \eqref{stat constraint}. This iteration repeats until the bias and statistical constraints are satisfied. By \Cref{Variance} and \Cref{corollary optimatl}, $M = \mathcal{O}\left( 1 \right)$; thus the algorithm simply ensures that $M$ is large enough to make the CLT applicable, but it is kept constant in \Cref{adaptive alg. 2}.

In the first inner loop of Algorithm~\ref{adaptive alg. 2}, the number of time steps $N$ is multiplied by $2^{\frac{\alpha}{\gamma}}$, rounding up to the next integer, and the number of particles $P$ is doubled. Afterwards, in the inner if-clause, the number of particles $P$ is again doubled, and the statistical error is estimated using the CLT by approximation \eqref{statistical error estimation}.
\begin{algorithm}[!ht]
\caption{Adaptive Randomized Quasi-Monte Carlo Algorithm}\label{adaptive alg. 2}
\textbf{Inputs:} $\quad T, \text{TOL}_{\text{eff}}, \theta, B_1,B_2,V, \alpha, \beta, \gamma,c_\varepsilon$\;
Determine $M$, $N$, and $P$ using \eqref{optimal NMP}, rounding up to the next integer with $M\geq 20$\;
Generate one low-discrepancy Sobol sequence $\mathbf{x}^{2P}$ of size $2P$ and obtain $\boldsymbol{\xi}^{(m)}\left(\mathbf{x}^{2P}\right)$ for $m=1,\ldots,M$ by scrambling it\;
Define the nested coarse nets $\boldsymbol{\xi}^{(m)}\left(\mathbf{x}^P\right)$ using the first $P$ points of $\boldsymbol{\xi}^{(m)}\left(\mathbf{x}^{2P}\right)$, for $m=1,\ldots,M$\;
\For{$m = 1,\ldots,M$}{
	Compute $\boldsymbol{\psi}^{2P|\left\lceil 2^{\alpha/\gamma}N \right\rceil}\left( \boldsymbol{\xi}^{(m)}\left(\mathbf{x}^{2P}\right) \right)$ and $\boldsymbol{\psi}^{P|N}\left( \boldsymbol{\xi}^{(m)}\left(\mathbf{x}^{P}\right) \right)$\;
}
Compute $\mathcal{A}_{\text{rQMC}}^{2P|\left\lceil 2^{\alpha/\gamma}N \right\rceil|M}\left(\mathbf{x}^{2P}\right)$ and $\mathcal{A}_{\text{rQMC}}^{P|N|M}\left(\mathbf{x}^P\right)$\;
Estimate the bias by approximation \eqref{bias estimation}\;
Estimate the statistical error by approximation \eqref{statistical error estimation}\;
\While{\textup{$\text{bias} > \theta\text{TOL}_{\text{eff}}$ or $\text{stat}_{\text{error}} > (1-\theta)\text{TOL}_{\text{eff}}$}}{
\If{\textup{$\text{bias} > \theta\text{TOL}_{\text{eff}}$}}{
	$N \leftarrow \left\lceil 2^{\frac{\alpha}{\gamma}}N \right\rceil$\;
	$P \leftarrow 2P$\;
Generate one low-discrepancy Sobol sequence $\mathbf{x}^{2P}$ of size $2P$ and obtain $\boldsymbol{\xi}^{(m)}\left(\mathbf{x}^{2P}\right)$ for $m=1,\ldots,M$ by scrambling it\;
Define the nested coarse nets $\boldsymbol{\xi}^{(m)}\left(\mathbf{x}^P\right)$ using the first $P$ points of $\boldsymbol{\xi}^{(m)}\left(\mathbf{x}^{2P}\right)$, for $m=1,\ldots,M$\;
\For{$m = 1,\ldots,M$}{
	Compute $\boldsymbol{\psi}^{2P|\left\lceil 2^{\alpha/\gamma}N \right\rceil}\left( \boldsymbol{\xi}^{(m)}\left(\mathbf{x}^{2P}\right) \right)$ and $\boldsymbol{\psi}^{P|N}\left( \boldsymbol{\xi}^{(m)}\left(\mathbf{x}^{P}\right) \right)$\;
}
Compute $\mathcal{A}_{\text{rQMC}}^{2P|\left\lceil 2^{\alpha/\gamma}N \right\rceil|M}\left(\mathbf{x}^{2P}\right)$ and $\mathcal{A}_{\text{rQMC}}^{P|N|M}\left(\mathbf{x}^P\right)$\;
Estimate the bias by approximation \eqref{bias estimation}\;
Estimate the statistical error by approximation \eqref{statistical error estimation}\;
}
\If{\textup{$\text{stat}_{\text{error}} > (1-\theta)\text{TOL}_{\text{eff}}$}}{
	$P \leftarrow 2P$\;
Generate one low-discrepancy Sobol sequence $\mathbf{x}^{2P}$ of size $2P$ and obtain $\boldsymbol{\xi}^{(m)}\left(\mathbf{x}^{2P}\right)$ for $m=1,\ldots,M$ by scrambling it\;
Define the nested coarse nets $\boldsymbol{\xi}^{(m)}\left(\mathbf{x}^P\right)$ using the first $P$ points of $\boldsymbol{\xi}^{(m)}\left(\mathbf{x}^{2P}\right)$, for $m=1,\ldots,M$\;
\For{$m = 1,\ldots,M$}{
	Compute $\boldsymbol{\psi}^{2P|\left\lceil 2^{\alpha/\gamma}N \right\rceil}\left( \boldsymbol{\xi}^{(m)}\left(\mathbf{x}^{2P}\right) \right)$ and $\boldsymbol{\psi}^{P|N}\left( \boldsymbol{\xi}^{(m)}\left(\mathbf{x}^{P}\right) \right)$\;
}
Compute $\mathcal{A}_{\text{rQMC}}^{2P|\left\lceil 2^{\alpha/\gamma}N \right\rceil|M}\left(\mathbf{x}^{2P}\right)$ and $\mathcal{A}_{\text{rQMC}}^{P|N|M}\left(\mathbf{x}^P\right)$\;
Estimate the bias by approximation \eqref{bias estimation}\;
Estimate the statistical error by approximation \eqref{statistical error estimation}\;
}
}
\textbf{Return:}  $\mathcal{A}_{\text{rQMC}}^{P|N|M}\left(\mathbf{x}^P\right)$\;
\end{algorithm}

\FloatBarrier
\section{Numerical experiments}
\label{sec:Numerical Experiments}
This work performs numerical experiments on the Gaussian model from Appendix~\ref{sec:Gaussian Model}, which admits an analytic solution, and on the Kuramoto model from Appendix~\ref{sec:Kuramoto Model}. Both models were simulated up to $T=1$. The Gaussian model primarily verifies the rates, due to its known solution, whereas the Kuramoto experiments exercise all score estimation, time integration, and rQMC sampling. The analysis is conditional on a fixed kernel satisfying \nameref{ass:A}.

The drift matrices $A$ and $B$ of the Gaussian model \eqref{app Gaussian model} are given by $A=-\mathbf{1}$, where $\mathbf{1}\in\mathbb{R}^{d\times d}$ denotes a square matrix with every entry equal to $1$, and $B=\mathbf{1}$. The diffusion matrix is $C_{\text{diff}}=\sigma_0 \left( (1+c)I - \frac{c}{d}\mathbf{1}  \right)$, where $c = 0.4$, and $\sigma_0 = 0.5$. For the initial Gaussian density of the Gaussian model, this work sets $m = \mathbf{0} \in \mathbb{R}^d$, and $\Gamma = C_{\text{diff}}$.

The natural frequencies of the Kuramoto model \eqref{app Kuramoto model particle system} are sampled as $\nu \sim \mathcal{U}(-0.2,0.2)^d$. The initial distribution is $\mathcal{N}(m,\Gamma)$ with $m = \mathbf{0}\in\mathbb{R}^d$, $\Gamma = \sigma_0 \left( (1+c)I - \frac{c}{d}\mathbf{1}  \right)\in\mathbb{R}^{d\times d}$, and $c = 0.4$, $\sigma_0 = 0.4$. The diffusion of the Kuramoto model is given by $\sigma = \Gamma$.

Unless stated otherwise, this work approximates the score with Gaussian kernels. The choice of kernel regularization $\lambda$ has little effect on the accuracy for simple symmetric densities, such as Gaussian ones, but is critical for more complex, skewed densities (e.g., log-normal densities). This paper discusses the influence of the hyperparameter $\lambda$ in Appendices~\ref{Influence of the Hyperparameter Gaussian} and \ref{Influence of the Hyperparameter LogNormal}. For the numerical experiments, this work sets $\lambda = 10^{-3}$ to balance accuracy and numerical stability.


To avoid introducing prior structure, $r_i(x) = 0$ for all $i=1,\ldots,d$. For the kernel length scale $\ell$, twice the standard deviation of the current particle system was selected. This choice was made based on the results presented by \citet{GradientEstimation} and makes the kernel random and dependent on time and $P$. Including this construction in the assumptions was beyond the scope of this work.

Regarding rQMC sampling, this work generates a Sobol sequence using the MATLAB sobolset($\cdot$) function, applies scrambling via scramble($\cdot$, ``MatousekAffineOwen''), and then applies net($\cdot$) to obtain sequences of sizes $P$ and $2P$, if necessary. All experiments were conducted in MATLAB R2025b on the Apple M4 platform, with the random number generator initialized using rng($1$).

\subsection{Computational cost scaling of score function approximation}
\label{ssec:Work Rate}

\begin{figure}[htbp]
\centering
\includegraphics[width=.5\textwidth]{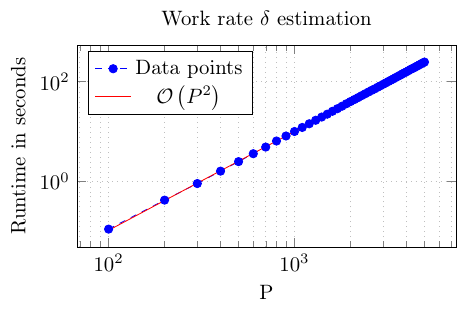}
\caption{Quantity of interest, given by $g(x)$, for the one-dimensional Gaussian model in Appendix~\ref{sec:Gaussian Model} is computed for varying numbers of particles $P$ using randomized quasi-Monte Carlo (rQMC) sampling. Parameters: $M=1$ and $N=30$.}
\label{fig:P vs work}
\end{figure}

\Cref{fig:P vs work} presents the runtime of the mean-field ODE's particle system approximation after $N=30$ forward Euler time steps and one realization $M=1$ with rQMC sampling, plotted against $P\in[100,5000]$. Gaussian kernels were applied to approximate the score at each time point. The computations were performed on the one-dimensional (1D) Gaussian model \eqref{app Gaussian model} with an initial mean of $m = 0$ and initial variance of $\Gamma = 0.5$. This study uses MATLAB's built-in backslash operator to compute the coefficients of the score approximation. Measuring the runtime for the score estimator computation yields a work rate of $\delta = 2$, a finite-range empirical scaling. The asymptotic dense direct solution is $\mathcal{O}(P^3)$ in floating point operations (FLOPs) and $\mathcal{O}(P^2)$ in memory. \Cref{fig:P vs work} only isolates the empirical particle-scaling of the score computation and is not intended as a runtime measurement of the full estimator.

\subsection{Convergence rates}
\label{ssec:Convergence Rates}
\begin{figure}[htbp]
  \centering
  \begin{subfigure}[t]{0.45\textwidth}
    \centering

\includegraphics[width=\textwidth]{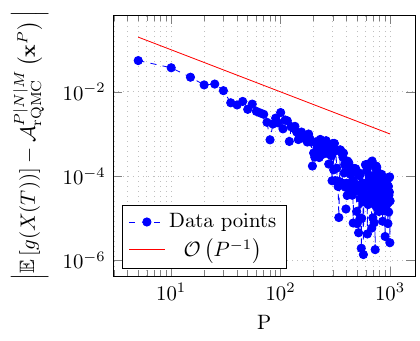}
\caption{Weak error convergence for the Gaussian model in Appendix~\ref{sec:Gaussian Model}, using randomized quasi-Monte Carlo (rQMC) sampling.}
\label{fig:P vs weak error gaussian}

  \end{subfigure}%
  \hfill
  \begin{subfigure}[t]{0.45\textwidth}
\centering

\includegraphics[width=\textwidth]{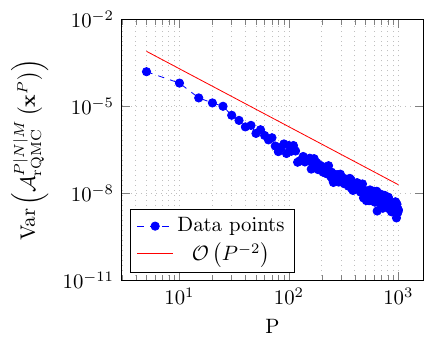}
\caption{Variance convergence for the Gaussian model in Appendix~\ref{sec:Gaussian Model}, using randomized quasi-Monte Carlo (rQMC) sampling.}
\label{fig:P vs estimator variance gaussian}

\end{subfigure}

\vspace{0.1em}

  \begin{subfigure}[t]{0.45\textwidth}
    \centering

\includegraphics[width=\textwidth]{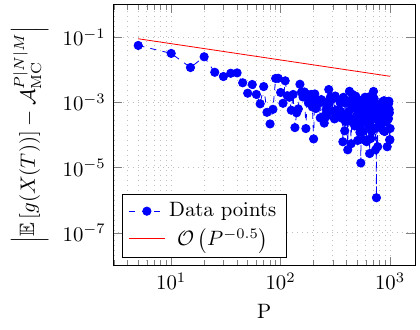}
\caption{Weak error convergence for the Gaussian model in Appendix~\ref{sec:Gaussian Model}, using Monte Carlo (MC) sampling.}
\label{fig:1d P vs weak error gaussian}

  \end{subfigure}%
  \hfill
  \begin{subfigure}[t]{0.45\textwidth}
\centering

\includegraphics[width=\textwidth]{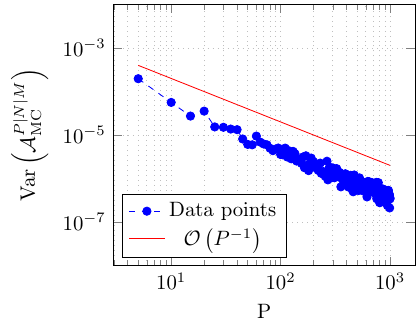}
\caption{Variance convergence for the Gaussian model in Appendix~\ref{sec:Gaussian Model}, using Monte Carlo (MC) sampling. }
\label{fig:1d P vs estimator variance gaussian}

\end{subfigure}

  \caption{Convergence rates for $g(x) = \frac{1}{d} \sum_{i=1}^d \cos(x_i)$ vs. $P$ for the five-dimensional Gaussian model in Appendix~\ref{sec:Gaussian Model}. Fourth-order Runge--Kutta (RK4) and Gaussian kernel regression were used, with both rQMC and MC sampling of initial states. Parameters: $M=40$, and $N=10$.}
  \label{fig:5 dim convergence rate gaussian}
\end{figure}

\begin{figure}[htbp]
  \centering

  \begin{subfigure}[t]{0.45\textwidth}
    \centering

\includegraphics[width=\textwidth]{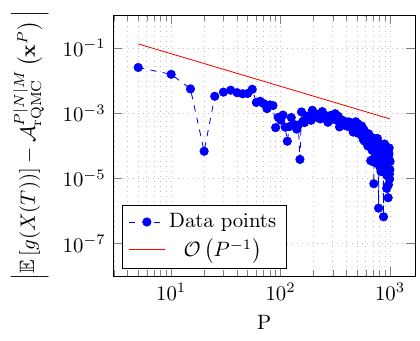}
\caption{Weak error convergence for the Kuramoto model in Appendix~\ref{sec:Kuramoto Model}, using randomized quasi-Monte Carlo (rQMC) sampling.}
\label{fig:P vs weak error kuramoto}

  \end{subfigure}%
  \hfill
  \begin{subfigure}[t]{0.45\textwidth}
\centering

\includegraphics[width=\textwidth]{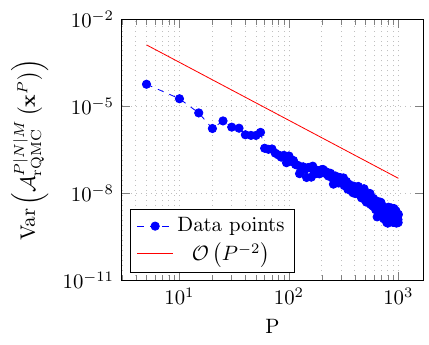}
\caption{Variance convergence for the Kuramoto model in Appendix~\ref{sec:Kuramoto Model}, using randomized quasi-Monte Carlo (rQMC) sampling. }
\label{fig:P vs estimator variance kuramoto}

\end{subfigure}

  \vspace{0.1em} 

  \begin{subfigure}[t]{0.45\textwidth}
    \centering

\includegraphics[width=\textwidth]{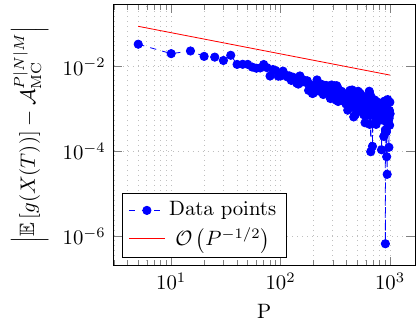}
\caption{Weak error convergence for the Kuramoto model in Appendix~\ref{sec:Kuramoto Model}, using Monte Carlo (MC) sampling.}
\label{fig:1d P vs weak error kuramoto}

  \end{subfigure}%
  \hfill
  \begin{subfigure}[t]{0.45\textwidth}
\centering

\includegraphics[width=\textwidth]{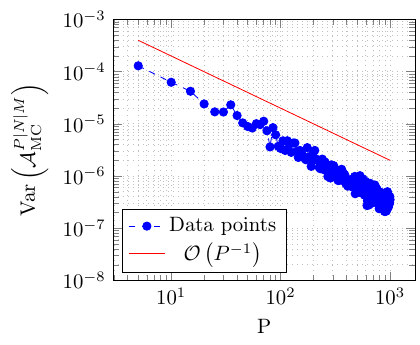}
\caption{Variance convergence for the Kuramoto model in Appendix~\ref{sec:Kuramoto Model}, using Monte Carlo (MC) sampling.}
\label{fig:1d P vs estimator variance kuramoto}

\end{subfigure}
  \caption{Convergence rates for $g(x) = \frac{1}{d} \sum_{i=1}^d \cos(x_i)$ vs. $P$ for the five-dimensional Kuramoto model in Appendix~\ref{sec:Kuramoto Model}. Fourth-order Runge--Kutta (RK4) and Gaussian kernel regression were used, with both rQMC and MC sampling. Parameters: $M=40$, and $N=10$.}
  \label{fig:5 dim convergence rate kuramoto}
\end{figure}

\Cref{fig:5 dim convergence rate gaussian} and \Cref{fig:5 dim convergence rate kuramoto} estimate the convergence rates for the weak error and variance, with computations on the five-dimensional Gaussian and Kuramoto models, respectively, as specified in Appendices~\ref{sec:Gaussian Model} and \ref{sec:Kuramoto Model}. The experiments apply the analytical solution to estimate the weak error for the Gaussian model. For the Kuramoto model, a Richardson extrapolation is applied with double the number of particles $P$ and the same number of time steps $N$. This work applies the classical RK scheme with $N=10$ time steps and $M=40$ independent randomizations of the initial points. The resulting weak-error is plotted. Figure panels \ref{fig:P vs weak error gaussian} and \ref{fig:P vs weak error kuramoto} present a weak error convergence rate of $\mathcal{O}(P^{-1})$, as allowed by \Cref{thm:rQMC-transport}. Moreover, Figure panels \ref{fig:P vs estimator variance gaussian} and \ref{fig:P vs estimator variance kuramoto} reveal a variance convergence rate of $\mathcal{O}(P^{-2})$.

For comparison, the same computation is repeated with plain MC sampling, in Figure panels \ref{fig:1d P vs weak error gaussian}, \ref{fig:1d P vs estimator variance gaussian}, \ref{fig:1d P vs weak error kuramoto}, and \ref{fig:1d P vs estimator variance kuramoto}. Slower rates of $\alpha = 0.5$ and $\beta = 1$ were observed. Therefore, the adaptive algorithm with MC sampling is expected to have a computational complexity of $\mathcal{O}\left( \text{TOL}^{-\left( 4 + \frac{1}{\gamma}\right)} \right)$ by \Cref{complexity: optimal values}. Worse constants are also depicted in \Cref{fig:5 dim convergence rate gaussian} and \Cref{fig:5 dim convergence rate kuramoto}.

\begin{figure}[htbp]
  \centering
  \begin{subfigure}[t]{0.45\textwidth}
    \centering

\includegraphics[width=\textwidth]{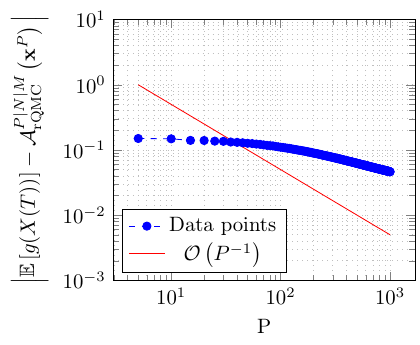}
\caption{Weak error convergence, using Gaussian kernel regression, randomized quasi-Monte Carlo (rQMC) sampling, and $N=10$.}
\label{fig:20 P vs weak error gaussian}

  \end{subfigure}%
  \hfill
  \begin{subfigure}[t]{0.45\textwidth}
\centering

\includegraphics[width=\textwidth]{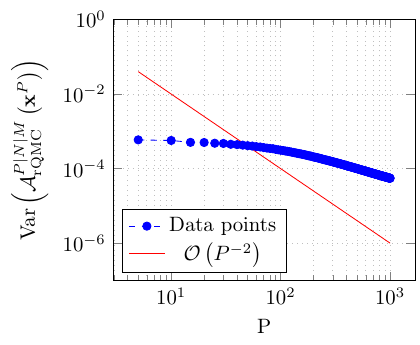}
\caption{Variance convergence, using RKHS kernel regression, randomized quasi-Monte Carlo (rQMC) sampling, and $N=10$.}
\label{fig:20 P vs estimator variance gaussian}

\end{subfigure}

  \vspace{0.1em} 

  \begin{subfigure}[t]{0.45\textwidth}
    \centering

\includegraphics[width=\textwidth]{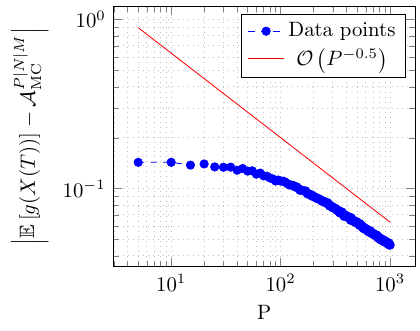}
\caption{Weak error convergence, using Gaussian kernel regression, Monte Carlo (MC) sampling, and $N=10$.}
\label{fig:20d P vs weak error gaussian}

  \end{subfigure}%
  \hfill
  \begin{subfigure}[t]{0.45\textwidth}
\centering

\includegraphics[width=\textwidth]{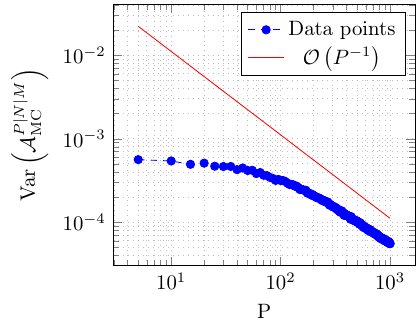}
\caption{Variance convergence, using Gaussian kernel regression, Monte Carlo (MC) sampling, and $N=10$.}
\label{fig:20d P vs estimator variance gaussian}

\end{subfigure}

  \caption{Convergence rates for $g(x) = \frac{1}{d} \sum_{i=1}^d \cos(x_i)$ vs. $P$ for the $20$-dimensional Gaussian model in Appendix~\ref{sec:Gaussian Model}. The fourth-order Runge--Kutta (RK4) method was applied. Gaussian kernel regression with randomized quasi-Monte Carlo (rQMC) and Monte Carlo (MC) sampling were applied. Parameter: $M=40$.}
  \label{fig:20 dim convergence rate kernel score}
\end{figure}

\begin{figure}[htbp]
  \centering

  \begin{subfigure}[t]{0.45\textwidth}
    \centering

\includegraphics[width=\textwidth]{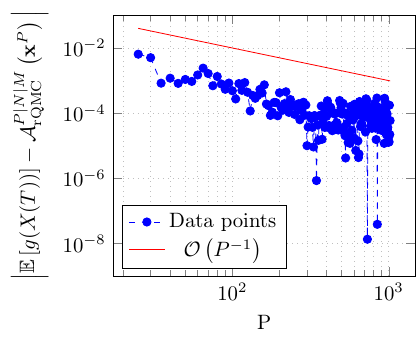}
\caption{Weak error convergence, using plug-in score, randomized quasi-Monte Carlo (rQMC) sampling, and $N=40$.}
\label{fig:P vs weak error gaussian true}

  \end{subfigure}%
  \hfill
  \begin{subfigure}[t]{0.45\textwidth}
\centering

\includegraphics[width=\textwidth]{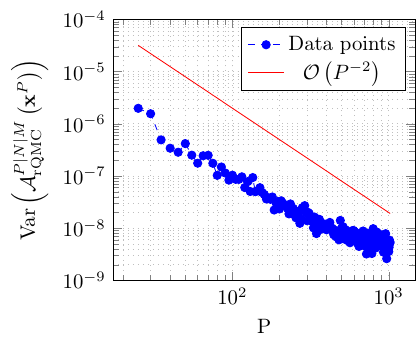}
\caption{Variance convergence, using plug-in score, randomized quasi-Monte Carlo (rQMC) sampling, and $N=40$.}
\label{fig:P vs estimator variance gaussian true}

\end{subfigure}

  \vspace{0.1em} 

  \begin{subfigure}[t]{0.45\textwidth}
    \centering

\includegraphics[width=\textwidth]{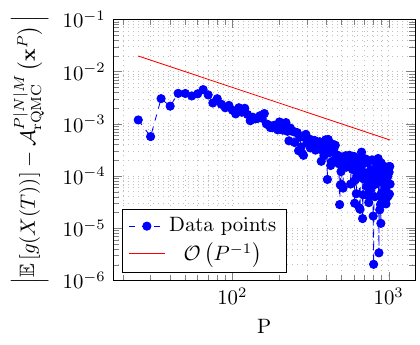}
\caption{Weak error convergence, using linear multivariate polynomial score estimation, randomized quasi-Monte Carlo (rQMC) sampling, and $N=40$.}
\label{fig:P vs weak error gaussian polynomial}

  \end{subfigure}%
  \hfill
  \begin{subfigure}[t]{0.45\textwidth}
\centering

\includegraphics[width=\textwidth]{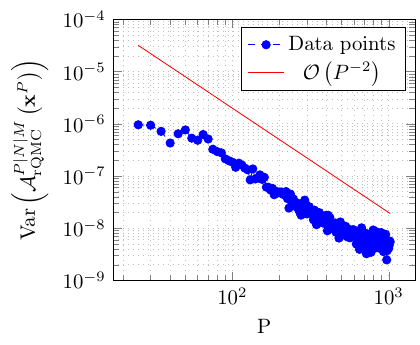}
\caption{Variance convergence, using linear multivariate polynomial score estimation, randomized quasi-Monte Carlo (rQMC) sampling, and $N=40$.}
\label{fig:P vs estimator variance gaussian polynomial}

\end{subfigure}
  \caption{Convergence rates for $g(x) = \frac{1}{d} \sum_{i=1}^d \cos(x_i)$ vs. $P$ for the $20$-dimensional Gaussian model in Appendix~\ref{sec:Gaussian Model}. The fourth-order Runge--Kutta (RK4) method was applied. Plug-in score, and a linear multivariate polynomial score estimation with rQMC sampling were applied. Parameter: $M=40$.}
  \label{fig:20 dim convergence rate non kernel score}
\end{figure}

\Cref{fig:20 dim convergence rate kernel score} and \Cref{fig:20 dim convergence rate non kernel score} display the convergence plots for the $d=20$-dimensional Gaussian model. Figure panels \ref{fig:20 P vs weak error gaussian} and \ref{fig:20 P vs estimator variance gaussian} do not present substantial improvements for the rQMC convergence rates, compared to the plain MC rates. Figure panels \ref{fig:20 P vs weak error gaussian}, \ref{fig:20 P vs estimator variance gaussian}, \ref{fig:20d P vs weak error gaussian}, and \ref{fig:20d P vs estimator variance gaussian}, confirm  that a very large number of particles is necessary to display asymptotic behavior and that further increasing $P$ for \Cref{fig:20 dim convergence rate kernel score} becomes computationally prohibitive.  Figure panels \ref{fig:20d P vs weak error gaussian} and \ref{fig:20d P vs estimator variance gaussian} indicate that the rQMC convergence rates degrade with an increasing value of $d$, an effect investigated in \citep{a46f288a95ff44a589103486f365ea79} for the QoI $g(x) = \frac{1}{d} \sum_{i=1}^d \cos(x_i)$, which considers all $20$ components of the interacting particle system.

To distinguish the moderate-dimensional performance of the rQMC method from the score approximation error, the $d=20$ Gaussian experiment is repeated with (i) the plug-in score in Figure panels \ref{fig:P vs weak error gaussian true} and \Cref{fig:P vs estimator variance gaussian true} and (ii) a polynomial score surrogate in Figure panels \ref{fig:P vs weak error gaussian polynomial} and \ref{fig:P vs estimator variance gaussian polynomial}. The first case employs the plug-in score of a multivariate Gaussian distribution with covariance matrix $\Gamma_t$ and mean $m_t$
\begin{align*}
\nabla \log \rho_t(x) = -\Gamma_t^{-1} (x-m_t),
\end{align*}
where $\Gamma_t$ and $m_t$ are estimated using samples at each time step. The expected rates of $\mathcal{O}\left(P^{-1}\right)$ and $\mathcal{O}\left(P^{-2}\right)$ are observed without a preasymptotic regime, indicating that the appropriate score representations recover the predicted rQMC scaling. Multivariate polynomials of degree at most $1$ are employed (similar to Section~5.1 in \citet{GradientEstimation}) for the second case. Rates near $\mathcal{O}\left(P^{-1}\right)$ and $\mathcal{O}\left(P^{-2}\right)$ are also observed in Figure panels \ref{fig:P vs weak error gaussian polynomial} and \ref{fig:P vs estimator variance gaussian polynomial}. However, at higher $P$ values, the RKHS regression error dominates and masks the rQMC convergence. \Cref{fig:20 dim convergence rate non kernel score} indicates that, in the present implementation and in higher dimensions, the primary preasymptotic limitation arises from the score approximation rather than the rQMC quadrature.


\subsection{Adaptive algorithm}
\label{ssec:Adaptive Algorithm}
In \Cref{ssec:Work Rate}, $\delta = 2$, and in \Cref{ssec:Convergence Rates}, $\alpha = 1$, $\beta = 2$, were observed; therefore, according to \Cref{corollary optimatl} the following is expected:
\begin{align*}
\mathcal{W} = \mathcal{O}\left( \text{TOL}_{\text{eff}}^{-2.25} \right)
\end{align*}
for the two test cases \eqref{app Gaussian model} and \eqref{app Kuramoto model}, using $g(x) = \frac{1}{d} \sum_{i=1}^d \cos(x_i)$. The constants $B_1$, $B_2$, and $V$ for initializing \Cref{adaptive alg. 2} were obtained from 
short pilot runs to estimate the convergence rates. This work selects an error split, dictated by $\theta = 0.7$. Although \eqref{optimal NMP} may suggest exceedingly small values of $M$, this study enforces $M\geq20$ as a practical minimum for estimating the randomization variance and forming an approximate normal confidence interval. This work sets $c_\varepsilon = 1.96$, corresponding to a $95\%$-confidence level. The estimates fall below the given tolerance, and the adaptive algorithm terminates without requiring the while-loop, suggesting that the starting constants $B_1$, $B_2$, and $V$ and, thus, the starting values $P$, $N$, and $M$ were set too conservatively. The outputs of \Cref{adaptive alg. 2} were validated against reference solutions. For the Gaussian model in Appendix~\ref{sec:Gaussian Model}, the analytic solution was employed. By \eqref{gaussian model cosinus},
\begin{align*}
\mathbb{E}\left[ \frac{1}{d}\sum_{i=1}^d \cos(X_i(T)) \right] \approx 6.197\cdot 10^{-1}.
\end{align*}

For the Kuramoto model in Appendix~\ref{sec:Kuramoto Model}, the reference solution  was computed with Euler--Maruyama MC, using $P=2000$ particles, $M=1000$ independent realizations, and $N=3000$ time steps. 

\begin{figure}[htbp]
  \centering
  \begin{subfigure}[t]{0.45\textwidth}
    \centering

\includegraphics[width=\textwidth]{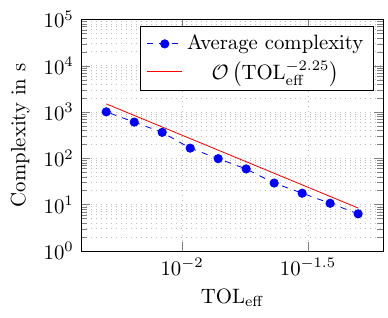}
\caption{Runtime for the Gaussian model in Appendix~\ref{sec:Gaussian Model}.}
\label{fig:Tol vs time gassian}

  \end{subfigure}%
  \hfill
  \begin{subfigure}[t]{0.45\textwidth}
    \centering

\includegraphics[width=\textwidth]{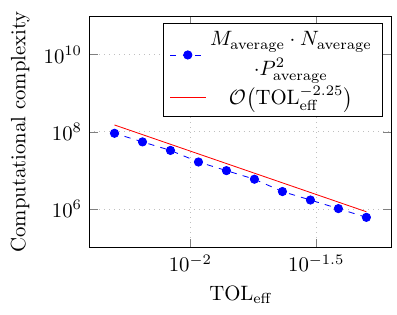}
\caption{Analytic computational cost for the Gaussian model in Appendix~\ref{sec:Gaussian Model}.}
\label{fig:Tol vs analytic work gassian}

  \end{subfigure}

\vspace{0.1em} 

  \begin{subfigure}[t]{0.45\textwidth}
    \centering

\includegraphics[width=\textwidth]{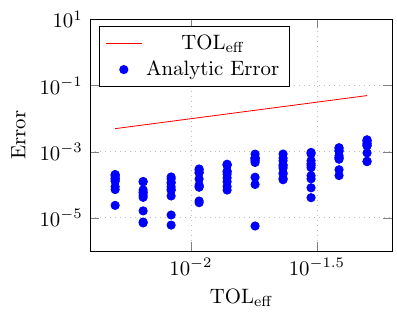}
\caption{Total error for the Gaussian model in Appendix~\ref{sec:Gaussian Model}.}
\label{fig:Tol-error gaussian}

  \end{subfigure}
  
  \caption{Adaptive algorithm output for $g(x) = \frac{1}{d} \sum_{i=1}^d \cos(x_i)$ for the five-dimensional Gaussian model in Appendix~\ref{sec:Gaussian Model}, using fourth-order Runge--Kutta (RK4) time integration, Gaussian kernel-regression score estimation, and randomized quasi-Monte Carlo (rQMC) sampling. Parameters: $\theta = 0.7$, $c_\varepsilon = 1.96$, $\alpha = 1$, $\beta = 2$, and $\gamma = 4$. The starting values $M_{\text{start}}$, $N_{\text{start}}$, and $P_{\text{start}}$ are listed in \Cref{tab. MNP gaussian}.}
  \label{fig:adaptive fine gaussian}
\end{figure}

For the MV-SDE \eqref{app Gaussian model} in Appendix~\ref{sec:Gaussian Model}, \Cref{adaptive alg. 2} was run $10$ times. \Cref{fig:Tol vs time gassian} plots the average computational cost of these outputs, and \Cref{fig:Tol-error gaussian} plots the individual errors of these realizations, using the true solution against the given tolerance $\text{TOL}_{\text{eff}}$. \Cref{fig:Tol vs analytic work gassian} presents the analytic computational cost:
\begin{align*}
\widetilde{\mathcal{W}} = M_{\text{average}} \cdot N_{\text{average}} \cdot P_{\text{average}}^2,
\end{align*}
leading to the same rate as in \Cref{fig:Tol vs time gassian} but with a much larger associated constant. The time to compute the constants $B_1$, $B_2$, and $V$ is not included in \Cref{fig:Tol vs time gassian}. 

\Cref{tab. MNP gaussian} reports the starting values $M_{\text{start}}$, $N_{\text{start}}$, and $P_{\text{start}}$ with the corresponding tolerances $\text{TOL}_{\text{eff}}$ and the corresponding values $M_{\text{average}}$, $N_{\text{average}}$ and $P_{\text{average}}$ for which the algorithm terminates, averaged over $10$ reruns.

\begin{table}[htbp]
  \centering
  \begin{subtable}[t]{1\textwidth}
    \centering

   \begin{tabular}{ c | c c c | c c c}
       $\text{TOL}_{\text{eff}}$ & $M_{\text{start}}$ & $N_{\text{start}}$  & $P_{\text{start}}$ & $M_{\text{average}}$ &	$N_{\text{average}}$  & $P_{\text{average}}$\\
    \hline
      $5.00\cdot 10^{-2}$  & $20$ & $4$ & $87$ &	 $20$ & $4$ & $87$\\ 
      $3.87\cdot 10^{-2}$  & $20$ & $4$ & $113$ &	 $20$ & $4$ & $113$\\ 
      $3.00\cdot 10^{-2}$  & $20$ & $4$ & $146$ &	 $20$ & $4$ & $146$\\
      $2.32\cdot 10^{-2}$  & $20$ & $4$ & $188$ &	 $20$ & $4$ & $188$\\ 
      $1.80\cdot 10^{-2}$  & $20$ & $5$ & $243$ &	 $20$ & $5$ & $243$\\ 
      $1.39\cdot 10^{-2}$  & $20$ & $5$ & $314$ &	 $20$ & $5$ & $314$\\  
      $1.08\cdot 10^{-2}$  & $20$ & $5$ & $406$ &	 $20$ & $5$ & $406$\\ 
      $8.34\cdot 10^{-3}$  & $20$ & $6$ & $524$ &	 $20$ & $6$ & $524$\\ 
      $6.46\cdot 10^{-3}$  & $20$ & $6$ & $676$ &	 $20$ & $6$ & $676$\\ 
      $5.00\cdot 10^{-3}$  & $20$ & $6$ & $874$ &	 $20$ & $6$ & $874$\\          
    \end{tabular}
    \caption{Starting and average terminating values for the corresponding tolerances $\text{TOL}_{\text{eff}}$ for the Gaussian model in Appendix~\ref{sec:Gaussian Model}.}
    \label{tab. MNP gaussian}

  \end{subtable}
  
    \vspace{0.1em} 
    
  \begin{subtable}[t]{1\textwidth}
    \centering

    \begin{tabular}{ c | c c c | c c c}
       $\text{TOL}_{\text{eff}}$ & $M_{\text{start}}$ & $N_{\text{start}}$  & $P_{\text{start}}$ & $M_{\text{average}}$ &	$N_{\text{average}}$  & $P_{\text{average}}$\\
    \hline
      $5.00\cdot 10^{-2}$  & $20$ & $3$ & $42$ &	 $20$ & $3$ & $42$\\ 
      $3.87\cdot 10^{-2}$  & $20$ & $3$ & $55$ &	 $20$ & $3$ & $55$\\ 
      $3.00\cdot 10^{-2}$  & $20$ & $3$ & $70$ &	 $20$ & $3$ & $70$\\
      $2.32\cdot 10^{-2}$  & $20$ & $3$ & $91$ &	 $20$ & $3$ & $91$\\ 
      $1.80\cdot 10^{-2}$  & $20$ & $4$ & $118$ &	 $20$ & $4$ & $118$\\ 
      $1.39\cdot 10^{-2}$  & $20$ & $4$ & $152$ &	 $20$ & $4$ & $152$\\  
      $1.08\cdot 10^{-2}$  & $20$ & $4$ & $196$ &	 $20$ & $4$ & $196$\\ 
      $8.34\cdot 10^{-3}$  & $20$ & $4$ & $253$ &	 $20$ & $4$ & $253$\\ 
      $6.46\cdot 10^{-3}$  & $20$ & $5$ & $327$ &	 $20$ & $5$ & $327$\\ 
      $5.00\cdot 10^{-3}$  & $20$ & $5$ & $422$ &	 $20$ & $5$ & $422$\\            
    \end{tabular}
    \caption{Starting and average terminating values for the corresponding tolerances for the Kuramoto model in Appendix~\ref{sec:Kuramoto Model}.}
    \label{tab. MNP kuramoto}

  \end{subtable} 
   \caption{Starting values $M_{\text{start}}$, $N_{\text{start}}$ and $P_{\text{start}}$ and the average terminating values $M_{\text{average}}$, $N_{\text{average}}$, and $P_{\text{average}}$ for the corresponding tolerances $\text{TOL}_{\text{eff}}$ for the five-dimensional Gaussian model in Appendix~\ref{sec:Gaussian Model} and the Kuramoto model in Appendix~\ref{sec:Kuramoto Model}, using the fourth-order Runge--Kutta (RK4) method and randomized quasi-Monte Carlo (rQMC) sampling. Parameters: $\theta = 0.7$, $c_\varepsilon = 1.96$, $\alpha=1$, $\beta = 2$, $\gamma = 4$. Figures~\ref{fig:adaptive fine gaussian} and \ref{fig:adaptive fine kuramoto} plot the results.}
\end{table}

\begin{figure}[htbp]
  \centering
 \begin{subfigure}[t]{0.45\textwidth}
\centering

\includegraphics[width=\textwidth]{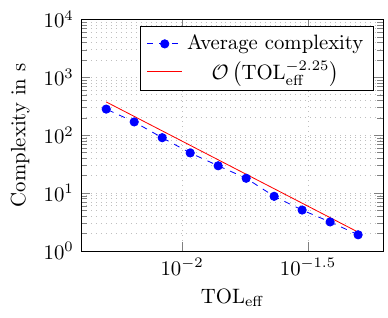}
\caption{Runtime for the Kuramoto model in Appendix~\ref{sec:Kuramoto Model}.}
\label{fig:Tol vs time kuramoto}

\end{subfigure}
  \hfill
 \begin{subfigure}[t]{0.45\textwidth}
\centering

\includegraphics[width=\textwidth]{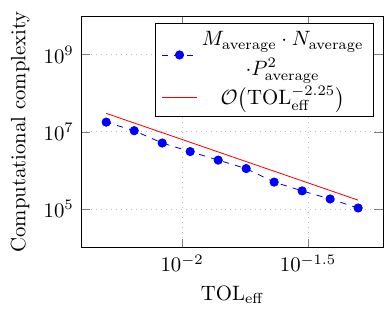}
\caption{Analytic computational cost for the Kuramoto model in Appendix~\ref{sec:Kuramoto Model}.}
\label{fig:Tol vs analytic work kuramoto}

\end{subfigure}

\vspace{0.1em} 

\begin{subfigure}[t]{0.45\textwidth}
\centering

\includegraphics[width=\textwidth]{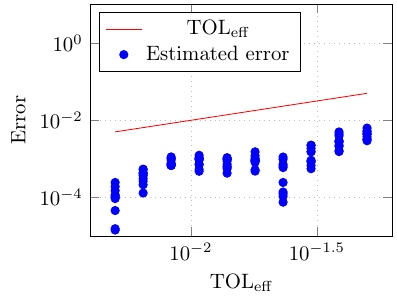}
\caption{Total error for the Kuramoto model in Appendix~\ref{sec:Kuramoto Model}. The Kuramoto reference is computed using the Euler--Maruyama Monte Carlo method; confidence intervals are reported in \Cref{tab. Kuramoto reference}.}
\label{fig:Tol-error kuramoto}

\end{subfigure}
  \caption{Adaptive algorithm output for $g(x) = \frac{1}{d} \sum_{i=1}^d \cos(x_i)$ for the five-dimensional Kuramoto model in Appendix~\ref{sec:Kuramoto Model}, using fourth-order Runge--Kutta (RK4) time integration, Gaussian kernel-regression score estimation, and randomized quasi-Monte Carlo (rQMC) sampling. Parameters: $\theta = 0.7$, $c_\varepsilon = 1.96$, $\alpha = 1$, $\beta = 2$, and $\gamma = 4$. The starting values $M_{\text{start}}$, $N_{\text{start}}$, and $P_{\text{start}}$ are listed in \Cref{tab. MNP kuramoto}.}
  \label{fig:adaptive fine kuramoto}
\end{figure}

The proposed adaptive algorithm terminated after one iteration; hence, the reported behavior is due to pilot-based initialization rather than repeated adaptation. To test the adaptive loop, the algorithm is initialized with the following rough starting values and run with the same tolerances as before:
\begin{align*}
M_{\text{start}} = 20,\quad
N_{\text{start}} = 1,\quad \text{and }
P_{\text{start}} = 10.
\end{align*}
\Cref{tab. MNP gaussian crude} presents the resulting values. The algorithm does not terminate after the first iteration and terminates with far fewer particles $P_{\text{average}}$. The number of time steps $N_{\text{average}}$ is less than the corresponding values in \Cref{tab. MNP gaussian}. 

\begin{figure}[htbp]
  \centering
  \begin{subfigure}[t]{0.45\textwidth}
    \centering

\includegraphics[width=\textwidth]{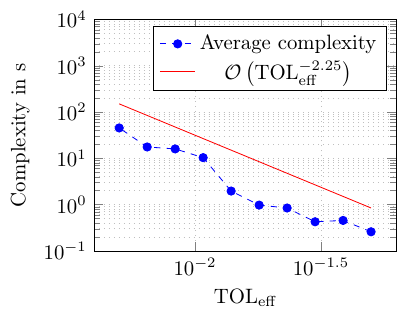}
\caption{Runtime for the Gaussian model in Appendix~\ref{sec:Gaussian Model}.}
\label{fig:Tol vs time gassian crude}

  \end{subfigure}%
  \hfill
  \begin{subfigure}[t]{0.45\textwidth}
    \centering

\includegraphics[width=\textwidth]{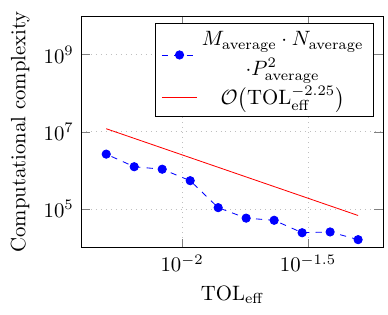}
\caption{Analytic computational cost for the Gaussian model in Appendix~\ref{sec:Gaussian Model}.}
\label{fig:Tol vs analytic work gassian crude}

  \end{subfigure}%

\vspace{0.1em} 

  \begin{subfigure}[t]{0.45\textwidth}
    \centering

\includegraphics[width=\textwidth]{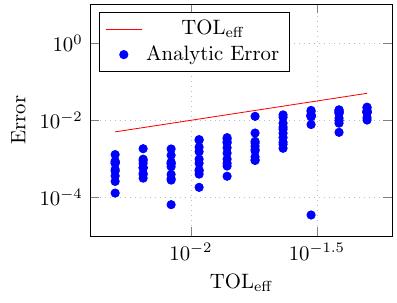}
\caption{Total error for the Gaussian model in Appendix~\ref{sec:Gaussian Model}.}
\label{fig:Tol-error gaussian crude}

  \end{subfigure}%

  \caption{Adaptive algorithm output for $g(x) = \frac{1}{d} \sum_{i=1}^d \cos(x_i)$ for the five-dimensional Gaussian model in Appendix~\ref{sec:Gaussian Model}, using fourth-order Runge--Kutta (RK4) time integration, Gaussian kernel-regression score estimation, and randomized quasi-Monte Carlo (rQMC) sampling. Parameters: $\theta = 0.7$, $c_\varepsilon = 1.96$, $\alpha = 1$, $\beta = 2$, and $\gamma = 4$. The starting values $M_{\text{start}}$, $N_{\text{start}}$, and $P_{\text{start}}$ are listed in \Cref{tab. MNP gaussian crude}.}
  \label{fig:adaptive crude gaussian}
\end{figure}

Figure panels \ref{fig:Tol vs time gassian crude} and \ref{fig:Tol vs analytic work gassian crude} plot the average computational complexity and the analytic computational complexity $\widetilde{\mathcal{W}}$, with a rate of $\mathcal{O}\left( \text{TOL}_{\text{eff}}^{-2.25} \right)$. All outputs are below the prescribed tolerance, but with tighter margins than before (\Cref{fig:Tol-error gaussian crude}).

\begin{figure}[htbp]
  \centering
  \begin{subfigure}[t]{0.45\textwidth}
\centering

\includegraphics[width=\textwidth]{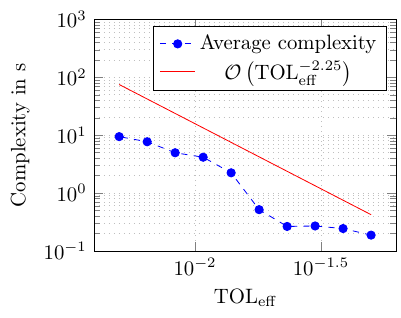}
\caption{Runtime for the Kuramoto model in Appendix~\ref{sec:Kuramoto Model}.}
\label{fig:Tol vs time kuramoto crude}

\end{subfigure}
  \hfill
  \begin{subfigure}[t]{0.45\textwidth}
\centering

\includegraphics[width=\textwidth]{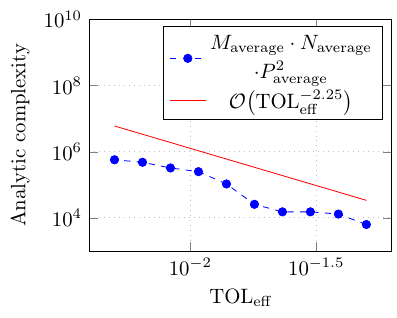}
\caption{Analytic computational cost for the Kuramoto model in Appendix~\ref{sec:Kuramoto Model}.}
\label{fig:Tol vs analytic work kuramoto crude}

\end{subfigure}

\vspace{0.1em} 

  \begin{subfigure}[t]{0.45\textwidth}
\centering

\includegraphics[width=\textwidth]{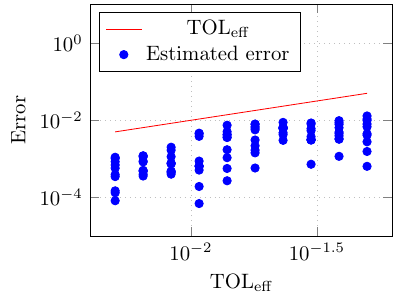}
\caption{Total error for the Kuramoto model in Appendix~\ref{sec:Kuramoto Model}. The Kuramoto reference is computed by using the Euler--Maruyama Monte Carlo method; \Cref{tab. Kuramoto reference} reports confidence intervals.}
\label{fig:Tol-error kuramoto crude}

\end{subfigure}
  \caption{Adaptive algorithm output for $g(x) = \frac{1}{d} \sum_{i=1}^d \cos(x_i)$ for the five-dimensional Kuramoto model in Appendix~\ref{sec:Kuramoto Model}, using fourth-order Runge--Kutta (RK4) time integration, Gaussian kernel-regression score estimation, and randomized quasi-Monte Carlo (rQMC) sampling. Parameters: $\theta = 0.7$, $c_\varepsilon = 1.96$, $\alpha = 1$, $\beta = 2$, and $\gamma = 4$. The starting values $M_{\text{start}}$, $N_{\text{start}}$, and $P_{\text{start}}$ are listed in \Cref{tab. MNP kuramoto crude}.}
  \label{fig:adaptive crude kuramoto}
\end{figure}

The same computations were repeated for the Kuramoto model in Appendix~\ref{sec:Kuramoto Model}, running \Cref{adaptive alg. 2} 10 times. The cost of computing $B_1$, $B_2$, and $V$ is not included. \Cref{fig:Tol vs time kuramoto} plots the average computational complexity of these 10 outputs. \Cref{fig:Tol-error kuramoto} presents the individual errors of these realizations, computed with the reference solution using Euler--Maruyama MC sampling, plotted against the given tolerance $\text{TOL}_{\text{eff}}$. \Cref{tab. MNP kuramoto} lists the starting values $M_{\text{start}}$, $N_{\text{start}}$, and $P_{\text{start}}$ with corresponding tolerances $\text{TOL}_{\text{eff}}$, and the corresponding values $M_{\text{average}}$, $N_{\text{average}}$, and $P_{\text{average}}$ for which the algorithm terminates, averaged over 10 reruns. A computational complexity of order $\mathcal{O}\left(\text{TOL}_{\text{eff}}^{-2.25}\right)$ is again expected, as observed in \Cref{fig:Tol vs time kuramoto}. As in \Cref{tab. MNP gaussian}, \Cref{tab. MNP kuramoto} reveals that, for every rerun, the initial values $M_{\text{start}}$, $N_{\text{start}}$, and $P_{\text{start}}$ are sufficient and that the adaptive algorithm immediately terminates, suggesting overly conservative estimates. \Cref{tab. MNP kuramoto crude} presents the results for the rough starting values, where the algorithm never terminates after the first iteration. As for the Gaussian model, far fewer particles $P_{\text{average}}$ and time steps $N_{\text{average}}$ are necessary compared with the values in \Cref{tab. MNP kuramoto}. \Cref{fig:Tol vs time kuramoto crude} and \Cref{fig:Tol vs analytic work kuramoto crude} plot the average computational complexity and product $M_{\text{average}} N_{\text{average}} P_{\text{average}}^2$. Again, the rate is $\mathcal{O}\left( \text{TOL}_{\text{eff}}^{-2.25} \right)$. In \Cref{fig:Tol-error kuramoto crude}, all outputs are below the tolerance, with several close to it.

\begin{table}[htbp]
  \centering
  \begin{subtable}[t]{1\textwidth}
    \centering

    \begin{tabular}{ c | c c c | c c c}
       $\text{TOL}_{\text{eff}}$ & $M_{\text{start}}$ & $N_{\text{start}}$  & $P_{\text{start}}$ & $M_{\text{average}}$ &	$N_{\text{average}}$  & $P_{\text{average}}$\\
    \hline
      $5.00\cdot 10^{-2}$  & $20$ & $1$ & $10$ &	 $20$ & $2.0$ & $20$\\ 
      $3.87\cdot 10^{-2}$  & $20$ & $1$ & $10$ &	 $20$ & $2.2$ & $24$\\ 
      $3.00\cdot 10^{-2}$  & $20$ & $1$ & $10$ &	 $20$ & $2.1$ & $24$\\
      $2.32\cdot 10^{-2}$  & $20$ & $1$ & $10$ &	 $20$ & $2.2$ & $34$\\ 
      $1.80\cdot 10^{-2}$  & $20$ & $1$ & $10$ &	 $20$ & $2.0$ & $38$\\ 
      $1.39\cdot 10^{-2}$  & $20$ & $1$ & $10$ &	 $20$ & $2.0$ & $52$\\  
      $1.08\cdot 10^{-2}$  & $20$ & $1$ & $10$ &	 $20$ & $2.5$ & $104$\\ 
      $8.34\cdot 10^{-3}$  & $20$ & $1$ & $10$ &	 $20$ & $2.9$ & $136$\\ 
      $6.46\cdot 10^{-3}$  & $20$ & $1$ & $10$ &	 $20$ & $3.0$ & $144$\\ 
      $5.00\cdot 10^{-3}$  & $20$ & $1$ & $10$ &	 $20$ & $3.3$ & $200$\\              
    \end{tabular}
    \caption{Starting and average terminating values for corresponding tolerances $\text{TOL}_{\text{eff}}$ for the Gaussian model in Appendix~\ref{sec:Gaussian Model}.}
    \label{tab. MNP gaussian crude}

  \end{subtable}
  
    \vspace{0.1em} 
    
  \begin{subtable}[t]{1\textwidth}
    \centering

    \begin{tabular}{ c | c c c | c c c}
       $\text{TOL}_{\text{eff}}$ & $M_{\text{start}}$ & $N_{\text{start}}$  & $P_{\text{start}}$ & $M_{\text{average}}$ &	$N_{\text{average}}$  & $P_{\text{average}}$\\
    \hline
      $5.00\cdot 10^{-2}$  & $20$ & $1$ & $10$ &	 $20$ & $1.4$ & $15$\\ 
      $3.87\cdot 10^{-2}$  & $20$ & $1$ & $10$ &	 $20$ & $1.8$ & $19$\\ 
      $3.00\cdot 10^{-2}$  & $20$ & $1$ & $10$ &	 $20$ & $1.9$ & $20$\\
      $2.32\cdot 10^{-2}$  & $20$ & $1$ & $10$ &	 $20$ & $1.9$ & $20$\\ 
      $1.80\cdot 10^{-2}$  & $20$ & $1$ & $10$ &	 $20$ & $1.9$ & $26$\\ 
      $1.39\cdot 10^{-2}$  & $20$ & $1$ & $10$ &	 $20$ & $2.3$ & $48$\\  
      $1.08\cdot 10^{-2}$  & $20$ & $1$ & $10$ &	 $20$ & $2.7$ & $68$\\ 
      $8.34\cdot 10^{-3}$  & $20$ & $1$ & $10$ &	 $20$ & $2.8$ & $76$\\ 
      $6.46\cdot 10^{-3}$  & $20$ & $1$ & $10$ &	 $20$ & $3.1$ & $88$\\ 
      $5.00\cdot 10^{-3}$  & $20$ & $1$ & $10$ &	 $20$ & $3.1$ & $96$\\            
    \end{tabular}
    \caption{Starting and average terminating values for corresponding tolerances for the Kuramoto model in Appendix~\ref{sec:Kuramoto Model}.}
    \label{tab. MNP kuramoto crude}

  \end{subtable} 
   \caption{Starting values $M_{\text{start}}$, $N_{\text{start}}$, and $P_{\text{start}}$ and average terminating values $M_{\text{average}}$, $N_{\text{average}}$, and $P_{\text{average}}$ for corresponding tolerances $\text{TOL}_{\text{eff}}$ for the five-dimensional Gaussian model in Appendix~\ref{sec:Gaussian Model} and the Kuramoto model in Appendix~\ref{sec:Kuramoto Model}, using the fourth-order Runge--Kutta (RK4) method and randomized quasi-Monte Carlo (rQMC) sampling. Parameters: $\theta = 0.7$, $c_\varepsilon = 1.96$, $\alpha=1$, $\beta = 2$, $\gamma = 4$. Figures~\ref{fig:adaptive crude gaussian} and \ref{fig:adaptive crude kuramoto} plot the results.}
\end{table}

\Cref{tab:summary of rates} summarizes the observed rates $\alpha$, $\beta$, $\gamma$, and $\delta$, and \Cref{ssec:Convergence Rates} depicts $\alpha$ and $\beta$ for plain MC and rQMC sampling. The predicted and observed rates agree, aside from the $20$-dimensional Gaussian model in \Cref{fig:20 dim convergence rate kernel score}. However, the expected rates are recovered using a polynomial score-estimation space or the plug-in score in \Cref{fig:20 dim convergence rate non kernel score}. The rate $\gamma = 4$ was observed in \Cref{ssec:Adaptive Algorithm}. For the rate $\delta$, this work must assume that $\delta = 3$ holds asymptotically; however, in \Cref{ssec:Work Rate} and \Cref{ssec:Adaptive Algorithm} $\delta=2$ was observed, due to MATLAB's optimized backslash operator.
\begin{table}[htbp]
\begin{center}
\begin{tabular}{c|c|c|c|c} 
& \makecell{Predicted rate,\\MC sampling}
& \makecell{Observed fit,\\MC sampling}
& \makecell{Predicted rate,\\rQMC sampling}
& \makecell{Observed fit,\\rQMC sampling} \\
\hline
$\alpha$ & $1/2$ & $1/2$ & $1$ & $1$ \\
$\beta$ & $1$ & $1$ & $2$ & $2$\\
$\gamma$ & $-$ & $-$ & $4$ & $4$\\
$\delta$ & $3$ & $2$ & $3$ & $2$\\
\end{tabular}
\caption{Predicted and observed rates $\alpha$, $\beta$, $\gamma$, and $\delta$, (see Section~\ref{sec:Numerical Experiments} for the Gaussian and Kuramoto models), using either Monte Carlo (MC) or randomized quasi-Monte Carlo (rQMC) sampling.}
\label{tab:summary of rates}
\end{center}
\end{table}

\section{Conclusion}
\label{chap:conclusion}

Under a probability-flow assumption, this work represents the original MV-SDE marginals via a deterministic mean-field ODE. The effect of replacing the exact score function by an  RKHS surrogate is controlled through a prescribed score residual, which is quantified in Theorem \ref{thm:vanishing-wasserstein-stability} under \nameref{ass:A}. The weak error of the deterministic particle system is then derived under \nameref{ass:A} and \nameref{ass:B}, while the variance across randomizations, time-discretization bias, and particle-work rate enter as separate quantitative models. In the tested moderate-dimensional regimes, these models are consistent with an effective tolerance-to-work complexity exponent near $2.25$. In the assessed examples, the effective exponents were close to $\alpha\approx 1$, $\beta\approx 2$, $\gamma \approx 4$, and  $\delta\approx 2$, yielding an empirical cost close to  $\mathcal{O} \left( \text{TOL}^{ -2.25 } \right)$. 

Several extensions are possible, including alternative score estimators (e.g., wavelet basis functions or other minimization spaces); iterative solvers, low-rank or fast-kernel methods for efficient linear solves; or multilevel MC/rQMC methods. Due to the deterministic evolution of the mean-field ODE, importance sampling is an interesting application. Although using the Liouville equation preserves the marginal densities, the score estimation introduces a bias that can alter the tails. Hyvärinen's score matching was applied to approximate the score in the weighted $L_{\mu_t}^2$ space. Incorporating a score estimator that better captures tail behavior is conceptually related to importance sampling, warranting further attention. Moreover, applying mean-field ODEs in the context of optimal control suggests the potential application of Pontryagin's maximum principle to derive the necessary conditions, instead of using the Hamilton--Jacobi--Bellman equation, which becomes intractable in high dimensions. Overall, these research directions offer several paths for advancing the theoretical and numerical aspects of the proposed approach.
\clearpage
\appendix
\label{sec:appendix}

\section{Score Estimator}
\label{sec:Score Estimator}
The objective function of the minimization problem \eqref{argmin} is given by
\begin{align*}
J_\mu(\phi) &=
\int_{\mathbb{R}^d}
\left(
\phi^2 + 2 r_i\phi + 2 \partial_{x_i}\phi
\right)\, d\mu
+ \lambda \|\phi\|_{\mathcal{H}}^2.
\end{align*}
Taking an arbitrary perturbation $h\in\mathcal{H}$ with $\epsilon\in\mathbb{R}$ yields
\begin{align*}
J_\mu(\phi + \epsilon h) &=
\int_{\mathbb{R}^d}
\left(
\left( \phi + \epsilon h \right)^2 + 2 r_i\left( \phi + \epsilon h\right) + 2 \partial_{x_i}\left( \phi + \epsilon h\right)
\right)\, d\mu
+ \lambda \| \phi + \epsilon h\|_{\mathcal{H}}^2,
\end{align*}
which leads to the first variation
\begin{align*}
\delta J_\mu(\phi)[h] = \frac{d}{d\epsilon} J_\mu(\phi + \epsilon h) \Big|_{\epsilon=0} = 2 \int_{\mathbb{R}^d} \left( \phi h +   r_i h  +  \partial_{x_i} h \right) \ d\mu + 2 \lambda \langle \phi, h \rangle_\mathcal{H}.
\end{align*}
For the minimizer $\phi_{i,\mu}^*$, $\delta J_\mu(\phi_{i,\mu}^*)[h]=0$ is required for all $h\in\mathcal{H}$, that is,
\begin{align*}
\frac{1}{2} \delta J_\mu(\phi_{i,\mu}^*)[h] = \int_{\mathbb{R}^d} \left( \left( \phi_{i,\mu}^* + r_i \right) h  +  \partial_{x_i} h \right) \ d\mu + \lambda \langle \phi_{i,\mu}^*, h \rangle_\mathcal{H} = 0,\quad \forall h\in\mathcal{H}.
\end{align*}
Using
\begin{align*}
h(x) &= \langle h, K(x,\cdot) \rangle_\mathcal{H}
\\
\partial_{x_i}h(x) &= \langle h,\partial_{x_i} K(x,\cdot) \rangle_\mathcal{H}
\end{align*}
yields
\begin{align*}
&\frac{1}{2} \delta J_\mu(\phi_{i,\mu}^*)[h]
\\
 =&\langle h , \int_{\mathbb{R}^d} \left( \left( \phi_{i,\mu}^*(x) + r_i(x) \right) K(x,\cdot) + \partial_{x_i} K(x,\cdot) \right) \ \mu(dx) + \lambda \phi_{i,\mu}^* \rangle_\mathcal{H} = 0,\quad \forall h\in\mathcal{H}
\\
\Rightarrow &\int_{\mathbb{R}^d} \left( \left( \phi_{i,\mu}^*(x) + r_i(x) \right) K(x,\cdot) + \partial_{x_i} K(x,\cdot) \right) \ \mu(dx) + \lambda \phi_{i,\mu}^* = 0.
\end{align*}
For $f\in\mathcal{H}$ and $i=1,\ldots,d$, the following is defined:
\begin{align*}
&\mathcal{C}_\mu:\mathcal{H} \to \mathcal{H},\quad  \mathcal{C}_\mu f := \int_{\mathbb{R}^d} \left( K(\cdot,y) f(y) \right) \ \mu(dy)
\\
&\xi_{i,\mu} := \int_{\mathbb{R}^d} \left( \partial_{x_i} K(x,\cdot) \right) \ \mu(dx)\in\mathcal{H}.
\end{align*}
Thus,
\begin{align*}
&\int_{\mathbb{R}^d} \left( \left( \phi_{i,\mu}^*(x) + r_i(x) \right) K(x,\cdot) + \partial_{x_i} K(x,\cdot) \right) \ \mu(dx) + \lambda \phi_{i,\mu}^*
\\
=& \mathcal{C}_\mu ( \phi_{i,\mu}^* + r_i) + \xi_{i,\mu} + \lambda \phi_{i,\mu}^* = 0.
\end{align*}
Then, the minimizer solves
\begin{align}\label{score esitmation minimizer}
 \left(\lambda I + \mathcal{C}_\mu\right)\phi_{i,\mu}^* =& - \left( \mathcal{C}_\mu r_i + \xi_{i,\mu} \right),
\end{align}
with the identity $I:\mathcal{H}\to\mathcal{H}$.

\subsection{Properties of the Score Estimator}
\subsubsection*{Boundedness:}
By the boundedness of the kernel, \eqref{A kernel bound} of \nameref{ass:A}, the operator $\mathcal{C}_\mu$ is also bounded, that is,
\begin{align*}
\|\mathcal{C}_\mu\|_{\mathcal{H}\to\mathcal{H}} \leq B,
\end{align*}
for all $\mu\in\mathcal{P}_2(\mathbb{R}^d)$. By definition, $\mathcal{C}_\mu$ is also self-adjoint for all $\mu\in\mathcal{P}_2(\mathbb{R}^d)$. If $f,g\in\mathcal{H}$, then
\begin{equation}\label{Cmu self adjoint}
\begin{split}
\langle \mathcal{C}_\mu f,g \rangle_\mathcal{H} &= \langle \int_{\mathbb{R}^d} \left( K(\cdot,y) f(y) \right) \ \mu(dy), g\rangle_\mathcal{H}
\\
&=  \int_{\mathbb{R}^d} \left( f(y) \langle K(\cdot,y) , g\rangle_\mathcal{H} \right) \ \mu(dy)
\\
&= \int_{\mathbb{R}^d} \left( f(y) g(y) \right) \ \mu(dy)
\\
&= \int_{\mathbb{R}^d} \left(  \langle f, K(\cdot,y)\rangle_\mathcal{H} g(y) \right) \ \mu(dy)
\\
&= \langle  f,\mathcal{C}_\mu g \rangle_\mathcal{H}.
\end{split}
\end{equation}
Setting $g:=f$ in \eqref{Cmu self adjoint} reveals that $\mathcal{C}_\mu$ is positive for all $\mu\in\mathcal{P}_2(\mathbb{R}^d)$,
\begin{align*}
\langle \mathcal{C}_\mu f,f \rangle_\mathcal{H} = \int_{\mathbb{R}^d} |f(y)|^2 \ \mu(dy) \geq 0.
\end{align*}
This result implies injectivity of $\lambda I + \mathcal{C}_\mu$ for all $\lambda >0$ and all $\mu\in\mathcal{P}_2(\mathbb{R}^d)$ by
\begin{equation}\label{operator bounded from below}
\begin{split}
\langle (\lambda I + \mathcal{C}_\mu) f,f \rangle_\mathcal{H} &= \lambda \|f\|_\mathcal{H}^2 + \int_{\mathbb{R}^d} |f(y)|^2 \ \mu(dy) \geq \lambda \|f\|_\mathcal{H}^2
\\
\Rightarrow \|(\lambda I + \mathcal{C}_\mu)f\|_\mathcal{H}&\geq \lambda \|f\|_\mathcal{H}.
\end{split}
\end{equation}
In \eqref{operator bounded from below}, this work establishes that $(\lambda I + \mathcal{C}_\mu)$ is bounded from below; therefore $\text{Ran}(\lambda I + \mathcal{C}_\mu)$ is closed and $(\lambda I + \mathcal{C}_\mu)$ is injective. Because $(\lambda I + \mathcal{C}_\mu)$ is self-adjoint, 
\begin{align*}
\overline{\text{Ran}(\lambda I + \mathcal{C}_\mu)} = \left( \text{ker}(\lambda I + \mathcal{C}_\mu) \right)^\perp = \mathcal{H}
\end{align*}
and due to closedness,
\begin{align*}
\text{Ran}(\lambda I + \mathcal{C}_\mu) = \mathcal{H}.
\end{align*}
As a direct consequence $(\lambda I + \mathcal{C}_\mu)^{-1}$ is also self-adjoint for all $\lambda >0$ and all $\mu\in\mathcal{P}_2(\mathbb{R}^d)$.

This result implies that, for all $\lambda >0$ and all $\mu\in\mathcal{P}_2(\mathbb{R}^d)$, 
\begin{align*}
\| \left(\lambda I + \mathcal{C}_\mu\right)^{-1} \|_{\mathcal{H} \to \mathcal{H}}\leq \frac{1}{\lambda}.
\end{align*}
Therefore, for all $i=1,\ldots,d$ and all $\mu\in\mathcal{P}_2(\mathbb{R}^d)$
\begin{align*}
\|\phi_{i,\mu}^* \|_\mathcal{H} \leq \frac{1}{\lambda} \left( \|  \mathcal{C}_\mu r_i \|_\mathcal{H} + \|\xi_{i,\mu}\|_\mathcal{H} \right),
\end{align*}
with the further bound
\begin{align*}
\|\phi_{i,\mu}^* \|_\mathcal{H}  \leq \frac{B^2 + B}{\lambda},
\end{align*}
and that, by the RKHS embedding, $\|h\|_{L^\infty(\mathbb{R}^d)} \leq \sqrt{B} \|h\|_\mathcal{H}$ for all $h\in\mathcal{H}$. Hence
\begin{align*}
\|\phi_{i,\mu}^* \|_{L^\infty(\mathbb{R}^d)}   \leq\sqrt{B} \frac{B^2 + B}{\lambda}.
\end{align*}
Every component $i=1,\ldots,d$ is approximated; thus,
\begin{align*}
\partial_{x_i} \log \rho \approx r_i + \phi_{i,\mu}^*,
\end{align*}
\nameref{ass:A} yields
\begin{align}\label{score estimate Linfinity bound}
\| r_i + \phi_{i,\mu}^* \|_{L^\infty(\mathbb{R}^d)} \leq \sqrt{B}\frac{B^2 + B}{\lambda} + B.
\end{align}

\subsubsection*{Pointwise Difference for Different Measures:}
Let $\mu$, and $\nu$ be two probability measures and $\phi_{i,\mu}^*$, and $\phi_{i,\nu}^*$ be the corresponding minimizers, where $i=1,\ldots,d$, given by \eqref{score esitmation minimizer}.
Then,
\begin{align*}
& \left(\lambda I + \mathcal{C}_\mu\right)\phi_{i,\mu}^* -  \left(\lambda I + \mathcal{C}_\nu\right)\phi_{i,\nu}^* = - \left( \mathcal{C}_\mu r_i + \xi_{i,\mu} \right) + \left( \mathcal{C}_\nu r_i + \xi_{i,\nu} \right)
 \\
 \Leftrightarrow & \left( \lambda I + \mathcal{C}_\mu \right) \left( \phi_{i,\mu}^* -\phi_{i,\nu}^* \right)   = \left( \mathcal{C}_\nu - \mathcal{C}_\mu \right) r_i + \left( \xi_{i,\nu} - \xi_{i,\mu} \right)  + \left( \mathcal{C}_\nu - \mathcal{C}_\mu \right) \phi_{i,\nu}^*.
\end{align*}
For the rQMC analysis in \Cref{sssec:Randomized Quasi-Monte Carlo Sampling}, the pointwise difference is of interest:
\begin{align*}
 \left( \phi_{i,\mu}^* -\phi_{i,\nu}^* \right)(x)   = \langle \left( \lambda I + \mathcal{C}_\mu \right)^{-1} \left(\left( \mathcal{C}_\nu - \mathcal{C}_\mu \right) r_i + \left( \xi_{i,\nu} - \xi_{i,\mu} \right)  + \left( \mathcal{C}_\nu - \mathcal{C}_\mu \right) \phi_{i,\nu}^*\right),K(x,\cdot) \rangle_\mathcal{H},
\end{align*}
for any $x\in\mathbb{R}^d$. Examining the individual terms and combining the fact that $\left( \lambda I + \mathcal{C}_\mu \right)^{-1}$ is a bounded linear operator with Bochner integrability, the integration order is interchangeable, yielding:
\begin{align*}
 &\left( \lambda I + \mathcal{C}_\mu \right)^{-1} \left( \mathcal{C}_\nu - \mathcal{C}_\mu \right) r_i
 \\
 =& \left( \lambda I + \mathcal{C}_\mu \right)^{-1} \int_{\mathbb{R}^d} \left( K(\cdot,y) r_i(y) \right) \ \left( \nu - \mu \right)(dy)
 \\
  =& \int_{\mathbb{R}^d} \left( \left( \lambda I + \mathcal{C}_\mu \right)^{-1} K(\cdot,y) r_i(y) \right) \ \left( \nu - \mu \right)(dy),
\end{align*}
and
\begin{align*}
 &\left( \lambda I + \mathcal{C}_\mu \right)^{-1} \left( \xi_{i,\nu} - \xi_{i,\mu} \right)
 \\
 =& \left( \lambda I + \mathcal{C}_\mu \right)^{-1} \int_{\mathbb{R}^d} \left( \partial_{x_i} K(\cdot,y) \right) \ \left( \nu - \mu \right)(dy)
 \\
  =& \int_{\mathbb{R}^d} \left( \left( \lambda I + \mathcal{C}_\mu \right)^{-1} \partial_{x_i} K(\cdot,y) \right) \ \left( \nu - \mu \right)(dy).
\end{align*}
For the last term, this work applies the reproducing property and the fact that $\left( \lambda I + \mathcal{C}_\mu \right)^{-1}$ is self-adjoint:
\begin{align*}
& \langle \left( \lambda I + \mathcal{C}_\mu \right)^{-1} \left( \mathcal{C}_\nu - \mathcal{C}_\mu \right) \phi_{i,\nu}^*, K(x,\cdot) \rangle_\mathcal{H}
\\
=& \langle  \left( \mathcal{C}_\nu - \mathcal{C}_\mu \right) \phi_{i,\nu}^*, \left( \lambda I + \mathcal{C}_\mu \right)^{-1} K(x,\cdot) \rangle_\mathcal{H}
\\
=& \int_{\mathbb{R}^d} \left( \phi_{i,\nu}^*(z) \langle  K(z,\cdot), \left( \lambda I + \mathcal{C}_\mu \right)^{-1} K(x,\cdot) \rangle_\mathcal{H} \right) \ (\nu - \mu)(dz)
\\
=&\int_{\mathbb{R}^d} \left(  \phi_{i,\nu}^*(z) \psi_{\mu,x}(z) \right) \ (\nu - \mu)(dz),
\end{align*}
where 
\begin{align*}
\psi_{\mu,x}(z) := \langle K(z,\cdot), \left( \lambda I + \mathcal{C}_\mu \right)^{-1} K(x,\cdot)\rangle_\mathcal{H},\quad \forall x,z\in\mathbb{R}^d,\mu\in\mathcal{P}_2(\mathbb{R}^d),
\end{align*}
and
\begin{align}\label{new function space}
\mathcal{F}:= \{ z\mapsto h_1(z) h_2(z) \mid h_1\in\mathcal{H}, h_2\in\mathcal{H}, \|h_1\|_\mathcal{H}\leq 1, \|h_2\|_\mathcal{H}\leq 1 \}.
\end{align}
Then, under the trivial assumption $\phi_{i,\nu}^* \neq 0 \neq \psi_{\mu,x}$, the following can be bounded:
\begin{align*}
&|\langle \left( \lambda I + \mathcal{C}_\mu \right)^{-1} \left( \mathcal{C}_\nu - \mathcal{C}_\mu \right) \phi_{i,\nu}^*, K(x,\cdot) \rangle_\mathcal{H}|
\\
=&  \Bigg|\int_{\mathbb{R}^d} \left(  \phi_{i,\nu}^*(z) \psi_{\mu,x}(z) \right) \ (\nu - \mu)(dz) \Bigg|
\\
\leq  & \|\phi_{i,\nu}^*\|_{\mathcal{H}} \|\psi_{\mu,x}\|_{\mathcal{H}} \sup_{f\in\mathcal{F}} \Bigg| \int_{\mathbb{R}^d} f(z) \ (\nu - \mu)(dz) \Bigg|.
\end{align*}
Because $\|\phi_{i,\nu}^*\|_{\mathcal{H}} \leq \frac{B^2 + B}{\lambda}$ and $\|\psi_{\mu,x}\|_{\mathcal{H}} \leq \|\left( \lambda I + \mathcal{C}_\mu \right)^{-1}\|_{\mathcal{H}\to \mathcal{H}} \|K(x,\cdot)\|_\mathcal{H} \leq \frac{1}{\lambda}\sqrt{B}$,
\begin{align*}
&|\langle \left( \lambda I + \mathcal{C}_\mu \right)^{-1} \left( \mathcal{C}_\nu - \mathcal{C}_\mu \right) \phi_{i,\nu}^*, K(x,\cdot) \rangle_\mathcal{H}|
\\
\leq  & \sqrt{B}\frac{B^2 + B}{\lambda^2}  \sup_{f\in\mathcal{F}} \Bigg| \int_{\mathbb{R}^d} f(z) \ (\nu - \mu)(dz) \Bigg|.
\end{align*}
In total, using the reproducing property,
\begin{equation}\label{difference score estimator}
\begin{split}
&| \left( \phi_{i,\mu}^* -\phi_{i,\nu}^* \right)(x)| 
\\
\leq& \Bigg| \int_{\mathbb{R}^d} \left( \left[ \left( \lambda I + \mathcal{C}_\mu \right)^{-1} K(\cdot,y) \right](x) r_i(y) \right) \ \left( \nu - \mu \right)(dy) \Bigg|
\\
& +\Bigg|\int_{\mathbb{R}^d} \left( \left[ \left( \lambda I + \mathcal{C}_\mu \right)^{-1} \partial_{x_i} K(\cdot,y) \right](x) \right)  \ \left( \nu - \mu \right)(dy)\Bigg|
\\
&+ \sqrt{B}\frac{B^2 + B}{\lambda^2}  \sup_{f\in\mathcal{F}} \Bigg| \int_{\mathbb{R}^d} f(z) \ (\nu - \mu)(dz) \Bigg|.
\end{split}
\end{equation}
All integrals in the bound \eqref{difference score estimator} are well defined. By the reproducing property, the Cauchy--Schwarz inequality, and \eqref{A kernel bound}, every $h\in\mathcal{H}$ satisfies
\begin{align}\label{RKHS sup bound}
|h(x)| = |\langle h, K(x,\cdot)\rangle_\mathcal{H}| \leq \|h\|_\mathcal{H}\, K(x,x)^{1/2} \leq \sqrt{B}\,\|h\|_\mathcal{H},\quad x\in\mathbb{R}^d.
\end{align}
Hence the mappings
\begin{align*}
y &\mapsto \left[ \left( \lambda I + \mathcal{C}_\mu \right)^{-1} K(\cdot,y) \right](x) r_i(y), \text{ and}
\\
y  &\mapsto \left[ \left( \lambda I + \mathcal{C}_\mu \right)^{-1} \partial_{x_i} K(\cdot,y) \right](x)
\end{align*}
are bounded for all $x\in\mathbb{R}^d$ under \nameref{ass:A}, and are therefore integrable with respect to the finite signed measure $\nu-\mu$. For the last integral, let $f = h_1 h_2\in\mathcal{F}$. By \eqref{RKHS sup bound}, $|f(z)| \leq B\|h_1\|_\mathcal{H}\|h_2\|_\mathcal{H} \leq B$ for all $z\in\mathbb{R}^d$, so
\begin{align*}
\Bigg| \int_{\mathbb{R}^d} f(z) \ (\nu - \mu)(dz) \Bigg| \leq \int_{\mathbb{R}^d} |f(z)| \ |\nu - \mu|(dz) \leq 2B
\end{align*}
for any probability measures $\mu$ and $\nu$, including empirical measures $\bar{\mu}^P$ of $P\in\mathbb{N}$ particles. In particular, \eqref{RKHS sup bound} yields $\|h\|_{L^2_\nu(\mathbb{R}^d)} \leq \sqrt{B}\|h\|_\mathcal{H}$ for every probability measure $\nu$, so $\mathcal{H}\hookrightarrow L^2_\nu(\mathbb{R}^d)$ holds without further assumptions. This enables the application of the bound \eqref{difference score estimator} in the context of low discrepancies in \Cref{sssec:Randomized Quasi-Monte Carlo Sampling}.

\section{Proof of Theorem \ref{thm:vanishing-wasserstein-stability}}
\label{sec:Wasserstein}
\begin{proof}
Set
\[
e_t(x):=\nabla \log \rho_t - \left( r(x) +  \phi_{\mu_t^\infty}^*(x) \right),
\qquad
u(t):=\E\left[\|\hat{X}(t)-\hat{X}^\infty(t)\|_2^2\right].
\]
Because $\hat{X}(0)=\hat{X}^\infty(0)$ almost surely and both trajectories are absolutely continuous, \(u\) is also absolutely continuous on \([0,T]\). In addition,
\[
\pi_t:=\Law\left(\hat{X}(t),\hat{X}^\infty(t)\right)
\]
is a coupling of \(\mu_t\) and \(\mu_t^\infty\); thus, by the definition of the \(2\)-Wasserstein distance,
\[
W_2^2(\mu_t,\mu_t^\infty)\le \int \|x-y\|_2^2 \ d\pi_t(x,y)=u(t).
\]
Differentiating \(u\) for almost every \(t\) yields
\[
\frac12 u'(t)
=
\E\Big[\big\langle \hat{X}(t)-\hat{X}^\infty(t),\,
\hat{b}(\hat{X}(t),\mu_t)-b^\mathcal{H}(\hat{X}^\infty(t),\mu_t^\infty)\big\rangle\Big].
\]
Add and subtract \(b^\mathcal{H}(\hat{X}(t),\mu_t)\):
\[
\frac12 u'(t)=A_t+B_t,
\]
where
\[
A_t:=\E\Big[\big\langle \hat{X}(t)-\hat{X}^\infty(t),\,
b^\mathcal{H}(\hat{X}(t),\mu_t)-b^\mathcal{H}(\hat{X}^\infty(t),\mu_t^\infty)\big\rangle\Big],
\]
\[
B_t:=\E\Big[\big\langle \hat{X}(t)-\hat{X}^\infty(t),\,
\hat{b}(\hat{X}(t),\mu_t)-b^\mathcal{H}(\hat{X}(t),\mu_t)\big\rangle\Big].
\]
For \(A_t\), Lipschitz continuity \eqref{lipschitz of velocity general} yields
\[
\|b^\mathcal{H}(\hat{X}(t),\mu_t) - b^\mathcal{H}(\hat{X}^\infty(t),\mu_t^\infty)\|_2 \leq C(\|\hat{X}(t)-\hat{X}^\infty(t)\|_2 + W_2(\mu_t,\mu_t^\infty))
\]
hence,
\[
A_t
\le
C\,\E\left[\|\hat{X}(t)-\hat{X}^\infty(t)\|_2^2\right]
+
C\,\E\left[ \|\hat{X}(t)-\hat{X}^\infty(t)\|_2 \right]\,W_2(\mu_t,\mu_t^\infty).
\]
As \(W_2(\mu_t,\mu_t^\infty)\le u(t)^{1/2}\) and \(\E\left[\|\hat{X}(t)-\hat{X}^\infty(t)\|_2\right]\le u(t)^{1/2}\),
\[
A_t\le 2C u(t).
\]
For \(B_t\), due to the Lipschitz continuity in \eqref{Lipschitz MVSDE}, \eqref{Lipschitz_divergence},  \eqref{lipschitz of velocity general}, the bounds \eqref{D bound} and \eqref{score estimate Linfinity bound}, and since $\|D(x,\mu)-D(x,\nu)\|_2 \leq \sqrt{2B}\,C\,W_2(\mu,\nu)$ by \eqref{Lipschitz MVSDE} and \eqref{D bound},
\begin{align*}
&\|\hat{b}(x,\mu_t)-b^\mathcal{H}(x,\mu_t)\|_2 
\\
\leq &\|\hat{b}(x,\mu_t)- b^\mathcal{H}(x,\mu_t^\infty)\|_2 + \|b^\mathcal{H}(x,\mu_t^\infty) -b^\mathcal{H}(x,\mu_t)\|_2
\\
\leq&
\| D(x,\mu_t)e_t(x)\|_2 + \underbrace{C\left( 3 + \sqrt{2Bd}\left( \sqrt{B}\frac{B^2+B}{\lambda} + B \right) \right)}_{=:\tilde{C}} W_2(\mu_t,\mu_t^\infty).
\end{align*}
Therefore, using the bound \eqref{D bound} and \(W_2(\mu_t,\mu_t^\infty)\le u(t)^{1/2}\),
\begin{align*}
|B_t|
&\leq
\E\Big[\|\hat{X}(t)-\hat{X}^\infty(t)\|_2\, \left( \|D(\hat{X}(t),\mu_t)\|_{2}\,\|e_t(\hat{X}(t))\|_2 + \tilde{C} W_2(\mu_t,\mu_t^\infty)\right)\Big]
\\
&\leq
B\,u(t)^{1/2}\,\|e_t\|_{L^2(\mu_t)} + \tilde{C} u(t).
\end{align*}
Using the residual assumption componentwise,
\[
|B_t|\le B\,u(t)^{1/2} d^{1/2} \epsilon(t) + \tilde{C} u(t).
\]
Finally, Young's inequality \(2ab\le a^2+b^2\) is applied to the first term on the right-hand side:
\[
2|B_t|
\le
(2\tilde{C} + 1)u(t)+B^2 d \epsilon(t)^2.
\]
Combining the estimates for \(A_t\) and \(B_t\), yields
\[
u'(t)
\le
(4C + 2\tilde{C} + 1)u(t)+B^2 d \epsilon(t)^2.
\]
Because \(u(0)=0\), Gr\"onwall's lemma results in
\[
u(t)
\le
B^2d\int_0^t e^{(4C + 2\tilde{C} + 1)(t-s)}\epsilon(s)^2 \ ds.
\]
Using \(W_2(\mu_t,\mu_t^\infty)^2\le u(t)\), yields
\[
W_2(\mu_t,\mu_t^\infty)^2
\le
B^2d\int_0^t e^{(4C + 2\tilde{C} + 1)(t-s)}\epsilon(s)^2 \ ds,
\]
which is the desired estimate.
\end{proof}

\section{Proof of Theorem \ref{thm:rQMC-transport}}
\label{sec:RQMCTransport}
Next, \Cref{thm:rQMC-transport} is proved by expressing the error using the Kolmogorov backward equation (KBE). Then, \nameref{ass:B}  is employed to linearize the resulting difference.
\begin{proof}
Using the cost-to-go function from the KBE \eqref{rQMC: constant cost to go} results in
\begin{align*}
&\mathbf{g}_P\left( \hat{\mathbf{X}}^{\infty|0|\boldsymbol{\mathcal{T}}(\boldsymbol{\xi}(\mathbf{x}))}(T) \right) - \mathbf{g}_P\left( \bar{\mathbf{X}}^{P|0|\boldsymbol{\mathcal{T}}(\boldsymbol{\xi}(\mathbf{x}))}(T) \right)\\
=&  u\left(T,\hat{\mathbf{X}}^{\infty|0|\boldsymbol{\mathcal{T}}(\boldsymbol{\xi}(\mathbf{x}))}(T)\right) - u\left(0,\hat{\mathbf{X}}^{\infty|0|\boldsymbol{\mathcal{T}}(\boldsymbol{\xi}(\mathbf{x}))}(0)\right)\\
=&\int_0^T \frac{d}{dt} u\left(t,\hat{\mathbf{X}}^{\infty|0|\boldsymbol{\mathcal{T}}(\boldsymbol{\xi}(\mathbf{x}))}(t)\right) \ dt\\
=&\int_0^T \frac{\partial u}{\partial t} \left(t,\hat{\mathbf{X}}^{\infty|0|\boldsymbol{\mathcal{T}}(\boldsymbol{\xi}(\mathbf{x}))}(t)\right)  + \sum_{p=1}^P \nabla_{\zeta_p}u \left(t,\hat{\mathbf{X}}^{\infty|0|\boldsymbol{\mathcal{T}}(\boldsymbol{\xi}(\mathbf{x}))}(t)\right) \cdot \frac{d}{dt}\hat{X}^{\infty|0|\mathcal{T}_0(\xi(x_p))}(t)  \ dt.
\end{align*}
By the KBE \eqref{rQMC: constant cost to go},
\begin{align*}
&\frac{\partial u}{\partial t} \left(t,\hat{\mathbf{X}}^{\infty|0|\boldsymbol{\mathcal{T}}(\boldsymbol{\xi}(\mathbf{x}))}(t)\right)
\\
 = &-\sum_{p=1}^P  b^{\mathcal{H}}\left(\hat{X}^{\infty|0|\mathcal{T}_0(\xi(x_p))}(t),\bar{\mu}^{\hat{\mathbf{X}}^{\infty|0|\boldsymbol{\mathcal{T}}(\boldsymbol{\xi}(\mathbf{x}))}(t)} \right) \cdot \nabla_{\zeta_p}u\left(t,\hat{\mathbf{X}}^{\infty|0|\boldsymbol{\mathcal{T}}(\boldsymbol{\xi}(\mathbf{x}))}(t)\right),
\end{align*}
where $\bar{\mu}^{\hat{\mathbf{X}}^{\infty|0|\boldsymbol{\mathcal{T}}(\boldsymbol{\xi}(\mathbf{x}))}(t)}$ defines the empirical measure of the independent particle system $\hat{\mathbf{X}}^{\infty|0|\boldsymbol{\mathcal{T}}(\boldsymbol{\xi}(\mathbf{x}))}(t)$:
\begin{align*}
\bar{\mu}^{\hat{\mathbf{X}}^{\infty|0|\boldsymbol{\mathcal{T}}(\boldsymbol{\xi}(\mathbf{x}))}(t)} := \frac{1}{P}\sum_{p=1}^P \delta\left(\cdot-\hat{X}^{\infty|0|\mathcal{T}_0(\xi(x_p))}(t)\right).
\end{align*}
In contrast, the term $\bar{\mu}^{\bar{\mathbf{X}}^{P|0|\boldsymbol{\mathcal{T}}(\boldsymbol{\xi}(\mathbf{x}))}(t)}$ denotes the empirical measure of the coupled particle system $\bar{\mathbf{X}}^{P|0|\boldsymbol{\mathcal{T}}(\boldsymbol{\xi}(\mathbf{x}))}(t)$, defined in \eqref{empirical measure interacting particle system}.
In summary,
\begin{align*}
&\mathbf{g}_P\left( \hat{\mathbf{X}}^{\infty|0|\boldsymbol{\mathcal{T}}(\boldsymbol{\xi}(\mathbf{x}))}(T) \right) - \mathbf{g}_P\left( \bar{\mathbf{X}}^{P|0|\boldsymbol{\mathcal{T}}(\boldsymbol{\xi}(\mathbf{x}))}(T) \right)\\
=&  \int_0^T - \sum_{p=1}^P b^{\mathcal{H}}\left(\hat{X}^{\infty|0|\mathcal{T}_0(\xi(x_p))}(t),\bar{\mu}^{\hat{\mathbf{X}}^{\infty|0|\boldsymbol{\mathcal{T}}(\boldsymbol{\xi}(\mathbf{x}))}(t)}\right) \cdot \nabla_{\zeta_p}u\left(t,\hat{\mathbf{X}}^{\infty|0|\boldsymbol{\mathcal{T}}(\boldsymbol{\xi}(\mathbf{x}))}(t)\right)\\
& + \sum_{p=1}^P b^{\mathcal{H}}\left(\hat{X}^{\infty|0|\mathcal{T}_0(\xi(x_p))}(t),\mu_t^\infty\right) \cdot \nabla_{\zeta_p}u \left(t,\hat{\mathbf{X}}^{\infty|0|\boldsymbol{\mathcal{T}}(\boldsymbol{\xi}(\mathbf{x}))}(t)\right)  \ dt
\\
=& \int_0^T \sum_{p=1}^P  \left( b^{\mathcal{H}}\left(\hat{X}^{\infty|0|\mathcal{T}_0(\xi(x_p))}(t),\mu_t^\infty\right) - b^{\mathcal{H}}\left(\hat{X}^{\infty|0|\mathcal{T}_0(\xi(x_p))}(t),\bar{\mu}^{\hat{\mathbf{X}}^{\infty|0|\boldsymbol{\mathcal{T}}(\boldsymbol{\xi}(\mathbf{x}))}(t)}\right) \right)
\\
&\hspace{9cm}\cdot \nabla_{\zeta_p}u \left(t,\hat{\mathbf{X}}^{\infty|0|\boldsymbol{\mathcal{T}}(\boldsymbol{\xi}(\mathbf{x}))}(t)\right) \ dt.
\end{align*}
By \nameref{ass:B}, the original drift and diffusion functions can be linearized:
\begin{align*}
&b\left(\hat{X}^{\infty|0|\mathcal{T}_0(\xi(x_p))}(t),\mu_t^\infty\right) - b\left(\hat{X}^{\infty|0|\mathcal{T}_0(\xi(x_p))}(t),\bar{\mu}^{\hat{\mathbf{X}}^{\infty|0|\boldsymbol{\mathcal{T}}(\boldsymbol{\xi}(\mathbf{x}))}(t)}\right)
\\
=& \int_{\mathbb{R}^d} A_b\left(\hat{X}^{\infty|0|\mathcal{T}_0(\xi(x_p))}(t),y;\mu_t^\infty\right) \ \left(\mu_t^\infty -\bar{\mu}^{\hat{\mathbf{X}}^{\infty|0|\boldsymbol{\mathcal{T}}(\boldsymbol{\xi}(\mathbf{x}))}(t)}\right)(dy)
\\
& + R_b\left(\hat{X}^{\infty|0|\mathcal{T}_0(\xi(x_p))}(t);  \mu_t^\infty ,\bar{\mu}^{\hat{\mathbf{X}}^{\infty|0|\boldsymbol{\mathcal{T}}(\boldsymbol{\xi}(\mathbf{x}))}(t)}\right),
\\
&\nabla \cdot D\left(\hat{X}^{\infty|0|\mathcal{T}_0(\xi(x_p))}(t),\mu_t^\infty\right) - \nabla \cdot D\left(\hat{X}^{\infty|0|\mathcal{T}_0(\xi(x_p))}(t),\bar{\mu}^{\hat{\mathbf{X}}^{\infty|0|\boldsymbol{\mathcal{T}}(\boldsymbol{\xi}(\mathbf{x}))}(t)}\right)
\\
=& \int_{\mathbb{R}^d} A_{\nabla\cdot D}\left(\hat{X}^{\infty|0|\mathcal{T}_0(\xi(x_p))}(t),y;\mu_t^\infty\right) \ \left(\mu_t^\infty -\bar{\mu}^{\hat{\mathbf{X}}^{\infty|0|\boldsymbol{\mathcal{T}}(\boldsymbol{\xi}(\mathbf{x}))}(t)}\right)(dy)
\\
& + R_{\nabla\cdot D}\left(\hat{X}^{\infty|0|\mathcal{T}_0(\xi(x_p))}(t);  \mu_t^\infty ,\bar{\mu}^{\hat{\mathbf{X}}^{\infty|0|\boldsymbol{\mathcal{T}}(\boldsymbol{\xi}(\mathbf{x}))}(t)}\right),
\\
&D\left(\hat{X}^{\infty|0|\mathcal{T}_0(\xi(x_p))}(t),\mu_t^\infty\right) - D\left(\hat{X}^{\infty|0|\mathcal{T}_0(\xi(x_p))}(t),\bar{\mu}^{\hat{\mathbf{X}}^{\infty|0|\boldsymbol{\mathcal{T}}(\boldsymbol{\xi}(\mathbf{x}))}(t)}\right)
\\
=& \int_{\mathbb{R}^d} A_D\left(\hat{X}^{\infty|0|\mathcal{T}_0(\xi(x_p))}(t),y;\mu_t^\infty\right) \ \left(\mu_t^\infty -\bar{\mu}^{\hat{\mathbf{X}}^{\infty|0|\boldsymbol{\mathcal{T}}(\boldsymbol{\xi}(\mathbf{x}))}(t)}\right)(dy)
\\
& + R_D\left(\hat{X}^{\infty|0|\mathcal{T}_0(\xi(x_p))}(t);  \mu_t^\infty,\bar{\mu}^{\hat{\mathbf{X}}^{\infty|0|\boldsymbol{\mathcal{T}}(\boldsymbol{\xi}(\mathbf{x}))}(t)}\right).
\end{align*} 
Substituting and taking the expectation with regard to $\boldsymbol{\xi}$ yields a decomposition into the following terms:
\begin{align}
&\Big| \mathbb{E}_{\boldsymbol{\xi}} \left[ \mathbf{g}_P\left( \hat{\mathbf{X}}^{\infty|0|\boldsymbol{\mathcal{T}}(\boldsymbol{\xi}(\mathbf{x}))}(T) \right) \right]- \mathbb{E}_{\boldsymbol{\xi}} \left[\mathbf{g}_P\left( \bar{\mathbf{X}}^{P|0|\boldsymbol{\mathcal{T}}(\boldsymbol{\xi}(\mathbf{x}))}(T) \right) \right]\Big|
\notag\\
\leq & \mathbb{E}_{\boldsymbol{\xi}} \Bigg[ \Bigg| \int_0^T \sum_{p=1}^P  \Bigg( \int_{\mathbb{R}^d} A_b\left(\hat{X}^{\infty|0|\mathcal{T}_0(\xi(x_p))}(t),y;\mu_t^\infty\right) \ \mu_t^\infty(dy)
\notag \\
&\hspace{1cm} -\int_{\mathbb{R}^d} A_b\left(\hat{X}^{\infty|0|\mathcal{T}_0(\xi(x_p))}(t),y;\mu_t^\infty\right) \ \bar{\mu}^{\hat{\mathbf{X}}^{\infty|0|\boldsymbol{\mathcal{T}}(\boldsymbol{\xi}(\mathbf{x}))}(t)}(dy)
\notag\\
&\hspace{1cm} + R_b\left(\hat{X}^{\infty|0|\mathcal{T}_0(\xi(x_p))}(t);  \mu_t^\infty , \bar{\mu}^{\hat{\mathbf{X}}^{\infty|0|\boldsymbol{\mathcal{T}}(\boldsymbol{\xi}(\mathbf{x}))}(t)}\right) \Bigg)
\notag\\
&\hspace{9cm}\cdot \nabla_{\zeta_p}u \left(t,\hat{\mathbf{X}}^{\infty|0|\boldsymbol{\mathcal{T}}(\boldsymbol{\xi}(\mathbf{x}))}(t)\right) \ dt  \Bigg| \Bigg]
\notag\\
&+ \mathbb{E}_{\boldsymbol{\xi}} \Bigg[ \Bigg| \int_0^T \sum_{p=1}^P  \Bigg( \int_{\mathbb{R}^d} A_{\nabla\cdot D}\left(\hat{X}^{\infty|0|\mathcal{T}_0(\xi(x_p))}(t),y;\mu_t^\infty\right) \ \mu_t^\infty(dy) 
\notag\\
&\hspace{1cm} -\int_{\mathbb{R}^d} A_{\nabla\cdot D}\left(\hat{X}^{\infty|0|\mathcal{T}_0(\xi(x_p))}(t),y;\mu_t^\infty\right) \ \bar{\mu}^{\hat{\mathbf{X}}^{\infty|0|\boldsymbol{\mathcal{T}}(\boldsymbol{\xi}(\mathbf{x}))}(t)}(dy)
\notag\\
&\hspace{1cm} + R_{\nabla\cdot D}\left(\hat{X}^{\infty|0|\mathcal{T}_0(\xi(x_p))}(t);  \mu_t^\infty ,\bar{\mu}^{\hat{\mathbf{X}}^{\infty|0|\boldsymbol{\mathcal{T}}(\boldsymbol{\xi}(\mathbf{x}))}(t)}\right) \Bigg)
\notag\\
&\hspace{9cm}\cdot \nabla_{\zeta_p}u \left(t,\hat{\mathbf{X}}^{\infty|0|\boldsymbol{\mathcal{T}}(\boldsymbol{\xi}(\mathbf{x}))}(t)\right) \ dt \Bigg| \Bigg]
\notag\\
&+ \mathbb{E}_{\boldsymbol{\xi}} \Bigg[ \Bigg| \int_0^T \sum_{p=1}^P \Bigg( \int_{\mathbb{R}^d} A_D\left(\hat{X}^{\infty|0|\mathcal{T}_0(\xi(x_p))}(t),y;\mu_t^\infty\right) \ \mu_t^\infty(dy)
\notag \\
&\hspace{1cm} -\int_{\mathbb{R}^d} A_D\left(\hat{X}^{\infty|0|\mathcal{T}_0(\xi(x_p))}(t),y;\mu_t^\infty\right) \ \bar{\mu}^{\hat{\mathbf{X}}^{\infty|0|\boldsymbol{\mathcal{T}}(\boldsymbol{\xi}(\mathbf{x}))}(t)}(dy) 
\notag\\
&\hspace{1cm} + R_D\left(\hat{X}^{\infty|0|\mathcal{T}_0(\xi(x_p))}(t);  \mu_t^\infty ,\bar{\mu}^{\hat{\mathbf{X}}^{\infty|0|\boldsymbol{\mathcal{T}}(\boldsymbol{\xi}(\mathbf{x}))}(t)}\right) \Bigg) 
\notag\\
&\hspace{1cm}\cdot\Bigg(  r\left(\hat{X}^{\infty|0|\mathcal{T}_0(\xi(x_p))}(t)\right) + \phi_{\mu_t^\infty}^*\left(\hat{X}^{\infty|0|\mathcal{T}_0(\xi(x_p))}(t)\right) \Bigg)
\notag\\
&\hspace{9cm} \cdot \nabla_{\zeta_p}u \left(t,\hat{\mathbf{X}}^{\infty|0|\boldsymbol{\mathcal{T}}(\boldsymbol{\xi}(\mathbf{x}))}(t)\right) \ dt \Bigg| \Bigg]
\notag\\
&+ \mathbb{E}_{\boldsymbol{\xi}} \Bigg[ \Bigg| \int_0^T \sum_{p=1}^P D\left(\hat{X}^{\infty|0|\mathcal{T}_0(\xi(x_p))}(t),\bar{\mu}^{\hat{\mathbf{X}}^{\infty|0|\boldsymbol{\mathcal{T}}(\boldsymbol{\xi}(\mathbf{x}))}(t)}\right)
\notag\\
&\hspace{1cm} \cdot\Bigg( \phi_{\mu_t^\infty}^*\left(\hat{X}^{\infty|0|\mathcal{T}_0(\xi(x_p))}(t) \right)- \phi_{\bar{\mu}^{\hat{\mathbf{X}}^{\infty|0|\boldsymbol{\mathcal{T}}(\boldsymbol{\xi}(\mathbf{x}))}(t)}}^*\left(\hat{X}^{\infty|0|\mathcal{T}_0(\xi(x_p))}(t) \right) \Bigg)
\notag\\
&\hspace{9cm}\cdot \nabla_{\zeta_p}u \left(t,\hat{\mathbf{X}}^{\infty|0|\boldsymbol{\mathcal{T}}(\boldsymbol{\xi}(\mathbf{x}))}(t)\right)  \ dt \Bigg| \Bigg].\label{the differences}
\end{align}
By the inequality \eqref{difference score estimator}, 
\begin{align*}
&\Bigg| \phi_{\mu_t^\infty}^*\left(\hat{X}^{\infty|0|\mathcal{T}_0(\xi(x_p))}(t) \right)- \phi_{\bar{\mu}^{\hat{\mathbf{X}}^{\infty|0|\boldsymbol{\mathcal{T}}(\boldsymbol{\xi}(\mathbf{x}))}(t)}}^*\left(\hat{X}^{\infty|0|\mathcal{T}_0(\xi(x_p))}(t) \right)\Bigg| 
\\
\leq& \Bigg| \int_{\mathbb{R}^d}  \left( \lambda I + \mathcal{C}_{\mu_t^\infty} \right)^{-1} K\left(\hat{X}^{\infty|0|\mathcal{T}_0(\xi(x_p))}(t),y\right) r_i(y) \ \left( \mu_t^\infty - \bar{\mu}^{\hat{\mathbf{X}}^{\infty|0|\boldsymbol{\mathcal{T}}(\boldsymbol{\xi}(\mathbf{x}))}(t)} \right)(dy) \Bigg|
\\
& +\Bigg|\int_{\mathbb{R}^d}  \left( \lambda I + \mathcal{C}_{\mu_t^\infty} \right)^{-1} \partial_{x_i} K\left(\hat{X}^{\infty|0|\mathcal{T}_0(\xi(x_p))}(t),y\right)  \ \left( \mu_t^\infty - \bar{\mu}^{\hat{\mathbf{X}}^{\infty|0|\boldsymbol{\mathcal{T}}(\boldsymbol{\xi}(\mathbf{x}))}(t)} \right)(dy)\Bigg|
\\
&+ \sqrt{B}\frac{B^2 + B}{\lambda^2}  \sup_{f\in\mathcal{F}} \Bigg| \int_{\mathbb{R}^d} f(z) \ (\mu_t^\infty - \bar{\mu}^{\hat{\mathbf{X}}^{\infty|0|\boldsymbol{\mathcal{T}}(\boldsymbol{\xi}(\mathbf{x}))}(t)})(dz) \Bigg|,
\end{align*}
where $\mathcal{F}$ is defined in \eqref{new function space}. 

Define
\begin{align*}
\Delta_t^P := \mu_t^\infty - \bar{\mu}^{\hat{\mathbf{X}}^{\infty|0|\boldsymbol{\mathcal{T}}(\boldsymbol{\xi}(\mathbf{x}))}(t)}
\end{align*}
\begin{align*}
Q_t^P = \sup_{f\in\mathcal{G}_t} \Big| \int_{\mathbb{R}^d} f(z) \ \Delta_t^P(dz) \Big|.
\end{align*}
This gives
\begin{align*}
\Bigg\| \int_{\mathbb{R}^d} A_b(x,z;\mu_t^\infty) \ \Delta_t^P(dz) \Bigg\|_2 &\leq \sqrt{d}Q_t^P
\\
\Bigg\| \int_{\mathbb{R}^d} A_ {\nabla\cdot D}(x,z;\mu_t^\infty) \ \Delta_t^P(dz) \Bigg\|_2 &\leq \sqrt{d}Q_t^P
\end{align*}
and for the matrix-valued diffusion linearization
\begin{align*}
\Bigg\| \int_{\mathbb{R}^d} A_ {D}(x,z;\mu_t^\infty) \ \Delta_t^P(dz) \Bigg\|_2 &\leq dQ_t^P.
\end{align*}
These inequalities hold pathwise and uniformly in $x$. Also define
\begin{align*}
\mathscr{R}_t^P = \sup_{x\in\mathbb{R}^{d}} \max \Bigg\{ &\Big\|R_{b}\left(x;\mu_t^\infty,\bar{\mu}^{\hat{\mathbf{X}}^{\infty|0|\boldsymbol{\mathcal{T}}(\boldsymbol{\xi}(\mathbf{x}))}(t)}\right) \Big\|_2,  
\Big\|R_{D}\left(x;\mu_t^\infty,\bar{\mu}^{\hat{\mathbf{X}}^{\infty|0|\boldsymbol{\mathcal{T}}(\boldsymbol{\xi}(\mathbf{x}))}(t)}\right)\Big\|_2, 
\\ 
&\Big\|R_{\nabla\cdot D}\left(x;\mu_t^\infty,\bar{\mu}^{\hat{\mathbf{X}}^{\infty|0|\boldsymbol{\mathcal{T}}(\boldsymbol{\xi}(\mathbf{x}))}(t)}\right)\Big\|_2 \Bigg\}.
\end{align*}
Together, this gives
\begin{align*}
\Big\|b(x,\mu_t^\infty) - b\left(x,\bar{\mu}^{\hat{\mathbf{X}}^{\infty|0|\boldsymbol{\mathcal{T}}(\boldsymbol{\xi}(\mathbf{x}))}(t)}\right)\Big\|_2 &\leq \sqrt{d} Q_t^P + \mathscr{R}_t^P
\\
\Big\|\nabla \cdot D(x,\mu_t^\infty) - \nabla \cdot D\left(x,\bar{\mu}^{\hat{\mathbf{X}}^{\infty|0|\boldsymbol{\mathcal{T}}(\boldsymbol{\xi}(\mathbf{x}))}(t)}\right)\Big\|_2 &\leq \sqrt{d} Q_t^P + \mathscr{R}_t^P
\\
\Big\|D(x,\mu_t^\infty) - D\left(x,\bar{\mu}^{\hat{\mathbf{X}}^{\infty|0|\boldsymbol{\mathcal{T}}(\boldsymbol{\xi}(\mathbf{x}))}(t)}\right)\Big\|_2 &\leq d Q_t^P + \mathscr{R}_t^P
\end{align*}
As a result
\begin{align*}
&\Bigg| \phi_{\mu_t^\infty}^*\left(\hat{X}^{\infty|0|\mathcal{T}_0(\xi(x_p))}(t) \right)- \phi_{\bar{\mu}^{\hat{\mathbf{X}}^{\infty|0|\boldsymbol{\mathcal{T}}(\boldsymbol{\xi}(\mathbf{x}))}(t)}}^*\left(\hat{X}^{\infty|0|\mathcal{T}_0(\xi(x_p))}(t) \right)\Bigg| 
\\
\leq & \left( 2 + \sqrt{B}\frac{B^2 + B}{\lambda^2}  \right) Q_t^P
\\
\Rightarrow &\Bigg\| \phi_{\mu_t^\infty}^*\left(\hat{X}^{\infty|0|\mathcal{T}_0(\xi(x_p))}(t) \right)- \phi_{\bar{\mu}^{\hat{\mathbf{X}}^{\infty|0|\boldsymbol{\mathcal{T}}(\boldsymbol{\xi}(\mathbf{x}))}(t)}}^*\left(\hat{X}^{\infty|0|\mathcal{T}_0(\xi(x_p))}(t) \right)\Bigg\|_2 
\\
\leq & \sqrt{d}\left( 2 + \sqrt{B}\frac{B^2 + B}{\lambda^2}  \right) Q_t^P.
\end{align*}
Combining 
 $$\sup_{t\in[0,T], x\in\mathbb{R}^d} \| r(x) + \phi_{\mu_t^\infty}^*\|_2 \leq \sqrt{d}\left( \sqrt{B} \frac{B^2 + B}{\lambda} + B \right)$$ 
with 
$$\sup_{x\in\mathbb{R}^d, \mu\in\mathcal{P}(\mathbb{R}^d)} \|D(x,\mu)\|_2\leq B$$
yields
\begin{align*}
&\Big\| b^\mathcal{H}\left(\hat{X}^{\infty|0|\mathcal{T}_0(\xi(x_p))}(t),\mu_t^\infty\right) - b^\mathcal{H}\left(\hat{X}^{\infty|0|\mathcal{T}_0(\xi(x_p))}(t),\bar{\mu}^{\hat{\mathbf{X}}^{\infty|0|\boldsymbol{\mathcal{T}}(\boldsymbol{\xi}(\mathbf{x}))}(t)}\right) \Big\|_2
\\
\leq & 2 \left( \sqrt{d} Q_t^P + \mathscr{R}_t^P \right)+ \left( \sqrt{d}\left( \sqrt{B} \frac{B^2 + B}{\lambda} + B\right) \right) \left( d Q_t^P + \mathscr{R}_t^P \right)
\\
& + B \sqrt{d}\left( 2 + \sqrt{B}\frac{B^2 + B}{\lambda^2}  \right) Q_t^P.
\end{align*}
With the constants
\begin{align*}
C_0 =& 2\sqrt{d} + d\sqrt{d}\left( \sqrt{B} \frac{B^2 + B}{\lambda} + B \right) + \sqrt{d}B \left( 2 + \sqrt{B} \frac{B^2 + B}{\lambda^2} \right),
\\
C_1 =& \sqrt{d}\left( \sqrt{B}\frac{B^2 + B}{\lambda} + B \right)
\end{align*}
this simplifies to
\begin{align*}
&\Big\| b^\mathcal{H}\left(\hat{X}^{\infty|0|\mathcal{T}_0(\xi(x_p))}(t),\mu_t^\infty\right) - b^\mathcal{H}\left(\hat{X}^{\infty|0|\mathcal{T}_0(\xi(x_p))}(t),\bar{\mu}^{\hat{\mathbf{X}}^{\infty|0|\boldsymbol{\mathcal{T}}(\boldsymbol{\xi}(\mathbf{x}))}(t)}\right) \Big\|_2
\\
\leq & C_0 Q_t^P  + (2+C_1) \mathscr{R}_t^P,
\end{align*}
hence
\begin{align*}
&\Big|  \mathbf{g}_P\left( \hat{\mathbf{X}}^{\infty|0|\boldsymbol{\mathcal{T}}(\boldsymbol{\xi}(\mathbf{x}))}(T) \right)- \mathbf{g}_P\left( \bar{\mathbf{X}}^{P|0|\boldsymbol{\mathcal{T}}(\boldsymbol{\xi}(\mathbf{x}))}(T) \right) \Big|
\\
\leq & C_u \int_0^T \left[  C_0 Q_t^P  + (2+C_1) \mathscr{R}_t^P \right] \ dt.
\end{align*}
Taking the expectation with respect to $\boldsymbol{\xi}$ gives
\begin{align*}
&\Big| \mathbb{E}_{\boldsymbol{\xi}} \left[ \mathbf{g}_P\left( \hat{\mathbf{X}}^{\infty|0|\boldsymbol{\mathcal{T}}(\boldsymbol{\xi}(\mathbf{x}))}(T) \right) \right]- \mathbb{E}_{\boldsymbol{\xi}} \left[\mathbf{g}_P\left( \bar{\mathbf{X}}^{P|0|\boldsymbol{\mathcal{T}}(\boldsymbol{\xi}(\mathbf{x}))}(T) \right) \right]\Big|
\\
\leq & C_u C_0 \mathbb{E}_{\boldsymbol{\xi}} \left[\int_0^T  Q_t^P  \ dt \right]+ C_u(2 + C_1) \mathbb{E}_{\boldsymbol{\xi}} \left[\int_0^T  \mathscr{R}_t^P  \ dt \right]
\end{align*}
In \nameref{ass:B} 
\begin{align*}
\mathbb{E}_{\boldsymbol{\xi}} \left[\int_0^T  Q_t^P \ dt \right] \leq T \mathbb{E}_{\boldsymbol{\xi}} \left[ \sup_{t\in[0,T]} Q_t^P \right]  \leq T C_{\mathcal{G}} P^{-\alpha_Q}
\end{align*}
and
\begin{align*}
\mathbb{E}_{\boldsymbol{\xi}} \left[ \int_0^T  \mathscr{R}_t^P  \ dt \right] \leq T C_R P^{-\alpha_R}
\end{align*}
was assumed, with $\alpha_R>1$. In conclusion
\begin{align*}
&\Big| \mathbb{E}_{\boldsymbol{\xi}}\left[  \mathbf{g}_P\left( \hat{\mathbf{X}}^{\infty|0|\boldsymbol{\mathcal{T}}(\boldsymbol{\xi}(\mathbf{x}))}(T) \right) \right]- \mathbb{E}_{\boldsymbol{\xi}} \left[\mathbf{g}_P\left( \bar{\mathbf{X}}^{P|0|\boldsymbol{\mathcal{T}}(\boldsymbol{\xi}(\mathbf{x}))}(T) \right) \right] \Big|
\\
 \leq &T C_u C_0 C_{\mathcal{G}}P^{-\alpha_Q} + T C_u (2+C_1) C_R P^{-\alpha_R}.
\end{align*}
\end{proof}

\section{Proof of Theorem \ref{complexity: optimal values}}
\label{sec:Optimal Work Complexity}
\begin{proof}
First, $\text{TOL}_{\text{eff}} := (\text{TOL} - e_1)$ is defined, and the Lagrangian \eqref{Lagrangian ODE} results in the following constraints:
\begin{align}
\frac{\partial \mathcal{L}_{\text{ODE}}}{\partial M} &= NP^\delta - \lambda_2 \frac{c_\varepsilon^2 V}{M^2P^{\beta}} = 0\label{dM}
\\
\frac{\partial \mathcal{L}_{\text{ODE}}}{\partial N} &= MP^\delta - \frac{\lambda_1 B_2\gamma}{N^{\gamma + 1}} = 0\label{dN}
\\
\frac{\partial \mathcal{L}_{\text{ODE}}}{\partial P} &= \delta MNP^{\delta -1} - \frac{\lambda_1 B_1 \alpha}{P^{\alpha + 1}} - \frac{\lambda_2 c_\varepsilon^2 V \beta}{M P^{\beta + 1}} = 0 \label{dP}
\\
\frac{\partial \mathcal{L}_{\text{ODE}}}{\partial \lambda_1} &= \frac{B_1}{P^\alpha} + \frac{B_2}{N^\gamma} - \theta \text{TOL}_{\text{eff}} = 0\label{dlambda1}
\\
\frac{\partial \mathcal{L}_{\text{ODE}}}{\partial \lambda_2} &= \frac{c_\varepsilon^2 V}{MP^{\beta}} - (1-\theta)^2 \text{TOL}_{\text{eff}}^2 = 0\label{dlambda2}.
\end{align}
From the constraints in \eqref{dM} and \eqref{dN},
\begin{align*}
\lambda_2 = \frac{NP^{\beta+\delta} M^2}{V c_\varepsilon^2}, \quad \text{and} \quad
\lambda_1 = \frac{MP^\delta N^{\gamma + 1}}{B_2 \gamma}.
\end{align*}
Substituting this into the constraint in \eqref{dP} yields the following:
\begin{align*}
0 =& \delta MNP^{\delta-1} - \frac{B_1 \alpha}{B_2\gamma} MP^{\delta-\alpha-1}N^{\gamma + 1}  - \beta NMP^{\delta-1} \\
\Leftrightarrow 0 =& M P^{\delta-1} (\delta N - \frac{B_1 \alpha}{B_2\gamma} P^{-\alpha}N^{\gamma + 1}  - \beta N)  \\
\Rightarrow  P=&\left( \frac{B_1 \alpha}{B_2\gamma}\frac{N^\gamma}{\delta - \beta }\right)^{1/\alpha},
\end{align*}
which requires $\beta\neq \delta$. This result is applied in the fourth constraint \eqref{dlambda1}, to obtain the following:
\begin{align*}
\frac{B_2\gamma}{ \alpha}\frac{\delta - \beta }{N^\gamma} + \frac{B_2}{N^\gamma} &= \theta \text{TOL}_{\text{eff}}\\
\Leftrightarrow  N &= \left( B_2 \left( \frac{\gamma(\delta-\beta)}{\alpha} + 1 \right) \frac{1}{\theta \text{TOL}_{\text{eff}}} \right)^{1/\gamma} \\
\Rightarrow P& = \left( \frac{B_1(\delta-\beta + \frac{\alpha}{\gamma} )}{\delta-\beta} \frac{1}{\theta \text{TOL}_{\text{eff}}} \right)^{1/\alpha}.
\end{align*}
Finally, the last constraint, \eqref{dlambda2}, results in
\begin{align*}
M &=  \frac{c_\varepsilon^2V}{(1-\theta)^2} \left( \frac{(\delta-\beta)\theta}{B_1(\delta-\beta + \frac{\alpha}{\gamma})} \right)^{\beta/\alpha} \text{TOL}_{\text{eff}}^{\beta/\alpha - 2},
\end{align*}
which uses $\delta-\beta + \frac{\alpha}{\gamma}\neq 0$. The total optimal computational complexity is
\begin{equation*}
\begin{split}
\mathcal{W}=& \frac{c_\varepsilon^2V}{(1-\theta)^2} \left( \frac{(\delta-\beta)\theta}{B_1(\delta-\beta + \frac{\alpha}{\gamma})} \right)^{\beta/\alpha} \left( B_2 \left( \frac{\gamma(\delta-\beta)}{\alpha} + 1 \right) \frac{1}{\theta } \right)^{1/\gamma}\\
&\left( \frac{B_1(\delta-\beta + \frac{\alpha}{\gamma} )}{\delta-\beta} \frac{1}{\theta } \right)^{\delta/\alpha} \text{TOL}_{\text{eff}}^{-\left(2 - \frac{\beta}{\alpha} + \frac{1}{\gamma} + \frac{\delta}{\alpha} \right)}.
\end{split}
\end{equation*}
\end{proof} 

\section{Gaussian Model}
\label{sec:Gaussian Model}
Let $A$, $B$, $C_{\text{diff}}\in\mathbb{R}^{d\times d}$, and let $W(t)$ be a $d$-dimensional Brownian motion. The $d$-dimensional Gaussian MV-SDE is defined as follows:
\begin{equation}\label{app Gaussian model}
\begin{split}
dX(t) &= (A X(t) + B\mathbb{E}[X(t)]) dt + C_{\text{diff}}dW(t),\quad 0 < t \leq T\\
X(0)&\sim \mathcal{N}(m,\Gamma),
\end{split}
\end{equation}
with a mean of $m\in\mathbb{R}^d$ and a covariance of $\Gamma \in \mathbb{R}^{d\times d}$. Let $m_t := \mathbb{E}[X(t)]$, and $\Gamma_t = \text{Cov}(X(t))$. For the MV-SDE \eqref{app Gaussian model},
\begin{align*}
m_t' &= (A+B)m_t
\\
\Gamma_t' &= A \Gamma_t + \Gamma_t A^T + C_{\text{diff}}C_{\text{diff}}^T.
\end{align*}
Hence, $m_t = \exp( (A+B)t) m$. By variation of constants,
\begin{align*}
\Gamma_t = e^{At} \Gamma e^{A^T t} + \int_0^t e^{A(t-s)}  C_{\text{diff}} C_{\text{diff}}^T e^{A^T(t-s)} \ ds.
\end{align*}
Thus, $X(t) \sim \mathcal{N}(m_t,\Gamma_t)$ and
\begin{align}\label{gaussian model cosinus}
\mathbb{E}\left[ \cos(X_i(T)) \right] = \cos((m_T)_i) \exp \left( -\frac{1}{2}(\Gamma_T)_{ii} \right).
\end{align}

\section{Kuramoto Model} 
\label{sec:Kuramoto Model}
First introduced by Yoshiki Kuramoto, \citep{Kuramoto75,Kuramoto84}, the Kuramoto model originated in statistical physics and has been applied in physics \citep{KuramotoPhysics}, chemistry \citep{KuramotoChemistry}, medicine \citep{KuramotoMedizin}, and sociology \citep{KuramotoSociology}. In its interacting particle form, this model is given by 
\begin{equation}\label{app Kuramoto model particle system}
\begin{split}
dX_p(t) &= \left( \nu_p + \frac{1}{P}\sum_{q=1}^P \sin(X_p(t) - X_q(t)) \right) dt + \sigma dW_p(t),\quad t>0,p=1,\ldots,P\\
X_p(0) & \sim \mu,
\end{split}
\end{equation}
where $\nu_p$ denotes the intrinsic natural frequencies of particles $p=1,\ldots,P$. The system is coupled due to the sinusoidal term, which is applied componentwise, and every particle experiences random noise in the form of the independent Brownian motion $W_p(t)$. In the mean-field limit $P\to\infty$, the MV-SDE is
\begin{equation}\label{app Kuramoto model}
\begin{split}
dX(t) &= \left( \nu + \int_{\mathbb{R}^d} \sin(X(t) - x) \ \mu_t(dx) \right) dt + \sigma dW(t),\quad t>0\\
X(0) &\sim \mu.
\end{split}
\end{equation}
Reference solutions computed with MC and Euler--Maruyama methods are reported in \Cref{tab. Kuramoto reference}. Two finer estimates were computed for validation of the reference value.
\begin{table}
\centering
   \begin{tabular}{ c | c | c | c | c}
   $P$  &  $N$ & $M$  & Reference value & $95\%$-confidence interval\\
    \hline
    $2000$ & $3000$ & $1000$ & $0.4205$ &  $[0.420345, 0.420655]$
    \\
    $4000$ & $3000$ & $1000$ & $0.4205$ &   $[0.420392, 0.420608]$
    \\
    $2000$ & $6000$ & $1000$ & $0.4202$ &  $[0.420047, 0.420353]$
    \end{tabular}
    \caption{Reference values of the Kuramoto solution with $95\%$ confidence intervals. Computed using Euler--Maruyama and Monte Carlo (MC) sampling.}
    \label{tab. Kuramoto reference}
\end{table}

\section{Influence of the Hyperparameter $\lambda$ on the Gaussian Model}
\label{Influence of the Hyperparameter Gaussian}
\begin{figure}[htbp]
  \centering
  \begin{subfigure}[t]{0.31\textwidth}
    \centering

\includegraphics[width=\textwidth]{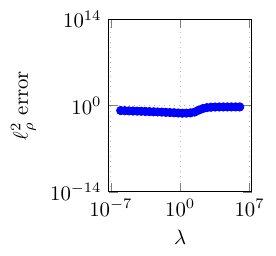}
\caption{$\ell^2_\rho$ score error.}
\label{fig:lambda vs L_rho}
  \end{subfigure}%
  \hfill
  \begin{subfigure}[t]{0.31\textwidth}
    \centering

\includegraphics[width=\textwidth]{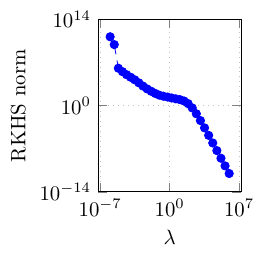}
\caption{Reproducing kernel Hilbert space (RKHS) norm of the score approximation.}
\label{fig:lambda vs rkhs}
  \end{subfigure}%
  \hfill
  \begin{subfigure}[t]{0.31\textwidth}
    \centering

\includegraphics[width=\textwidth]{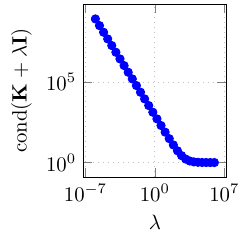}
\caption{Condition number of $\mathbf{K}+\lambda \mathbf{I}$.}
\label{fig:lambda vs condition}
  \end{subfigure}
  \caption{Gaussian kernel approximation errors of a one-dimensional Gaussian density and the condition number of $\mathbf{K}+\lambda \mathbf{I}$, for values of $\lambda$. The reproducing kernel Hilbert space (RKHS) norm increases rapidly as $\lambda\to 0$, reflecting an ill-conditioned $\mathbf{K}+\lambda \mathbf{I}$.}
  \label{fig:lambda}
\end{figure}

To quantify the influence of the regularization parameter $\lambda>0$, $P=1000$ independent, 1D points $x_p$, $p=1,\ldots,P$, are generated from the initial density of the Gaussian model in Appendix~\ref{sec:Gaussian Model}, and the approximation of the score is determined using Gaussian kernels with varying $\lambda$ values. Specifically,
\begin{align*}
K(x,y) = \exp\left(-\frac{\|x-y\|_2^2}{2\ell^2} \right),
\end{align*}
where the kernel length scale $\ell$ is set to twice the standard deviation of the particles. The average results, repeated over $M=100$ independent realizations, are reported. First, \Cref{fig:lambda vs L_rho} presents the discrete error norm for vector-valued functions $f$,
\begin{align*}
\|f\|_{\ell_\rho^2}^2 &= \frac{1}{P} \sum_{p=1}^P  \|f(x_p)\|_2^2,
\end{align*}
where $\rho(x)$ represents the probability density function of the initial density. Second, \Cref{fig:lambda vs rkhs} plots the RKHS norm of the score estimator. By the approximation \eqref{kernel regression particle system}, this approach requires computing the inverses of $(\mathbf{K}+\lambda \mathbf{I})$, which become nearly singular for exceedingly small values of $\lambda$.

\begin{figure}[htbp]
  \centering
  \begin{subfigure}[t]{0.33\textwidth}
    \centering

\includegraphics[width=\textwidth]{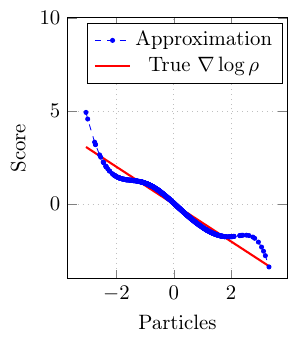}
\caption{$\lambda = 10^{-3}$.}
\label{fig:lambda m3}
  \end{subfigure}%
  \hfill
  \begin{subfigure}[t]{0.30\textwidth}
    \centering

\includegraphics[width=\textwidth]{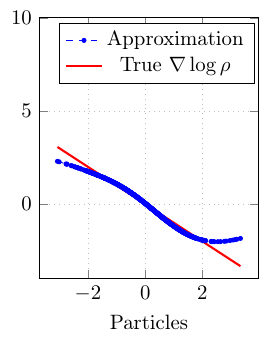}
\caption{$\lambda = 10^{0}$.}
\label{fig:lambda 0}
  \end{subfigure}%
  \hfill
  \begin{subfigure}[t]{0.30\textwidth}
    \centering

\includegraphics[width=\textwidth]{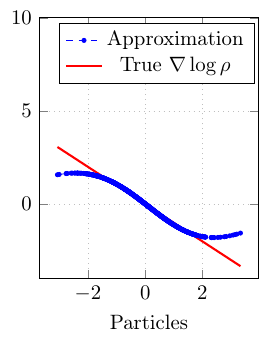}
\caption{$\lambda = 10^{1}$.}
\label{fig:lambda 1}
  \end{subfigure}

  \caption{Plot of the score for values of $\lambda$ for the one-dimensional Gaussian model with $P=1000$ particles.}
  \label{fig:lambda and gradlogs}
\end{figure}

The particles are also evolved in time using the same score approximation, with one forward Euler step of size $dt = 1/100$, and the resulting average was compared with the analytic solution for the QoI $g$. The error was constant, with values of $(6.54 \pm 0.01) \cdot 10^{-3}$. Finally, the condition number of the matrix $\mathbf{K}+\lambda \mathbf{I}$ is plotted, and \Cref{fig:lambda} shows the three quantities. \Cref{fig:lambda and gradlogs} also visualizes the score estimation for $\lambda\in\{10^{-3},10^{0},10^{1}\}$.

\section{Influence of the Hyperparameter $\lambda$ on a Log-Normal Model}
\label{Influence of the Hyperparameter LogNormal}
Similar to Figure~A2 in the work by \citet{GradientEstimation}, this work estimates the 1D score:
\begin{align*}
\nabla \log \rho(x) = \frac{m - \sigma^2 - \log(x)}{\sigma^2 x}
\end{align*}
of a log-normal density with parameters $m = 0$ and $\sigma = 0.5$ using Gaussian kernels with a length scale $\ell = 1$. Particle systems of size $P=1000$ with $M=100$ independent realizations were used.

\begin{figure}[htbp]
  \centering
  \begin{subfigure}[t]{0.31\textwidth}
    \centering

\includegraphics[width=\textwidth]{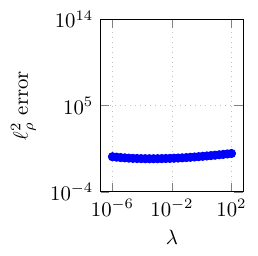}
\caption{$\ell^2_\rho$ score error.}
\label{fig:lambda vs L_rho ln}
  \end{subfigure}%
  \hfill
  \begin{subfigure}[t]{0.31\textwidth}
    \centering

\includegraphics[width=\textwidth]{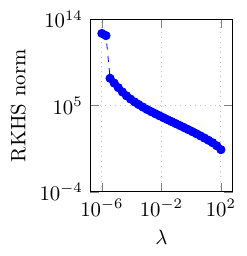}
\caption{Reproducing kernel Hilbert space (RKHS) norm of the score approximation.}

\label{fig:lambda vs rkhs ln}
  \end{subfigure}%
  \hfill
  \begin{subfigure}[t]{0.31\textwidth}
    \centering

\includegraphics[width=\textwidth]{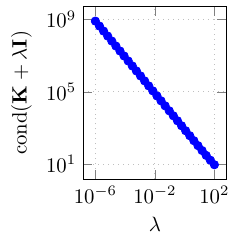}
\caption{Condition number of $\mathbf{K}+\lambda \mathbf{I}$.}
\label{fig:lambda vs condition ln}
  \end{subfigure}
  \caption{Gaussian kernel approximation errors of a one-dimensional log-normal density and the condition number of $\mathbf{K}+\lambda \mathbf{I}$, for values of $\lambda$. The reproducing kernel Hilbert space (RKHS) norm increases rapidly as $\lambda\to 0$, reflecting an ill-conditioned $\mathbf{K}+\lambda \mathbf{I}$.}
  \label{fig:lambda ln}
\end{figure}

\begin{figure}[htbp]
  \centering
  \begin{subfigure}[t]{0.33\textwidth}
    \centering

\includegraphics[width=\textwidth]{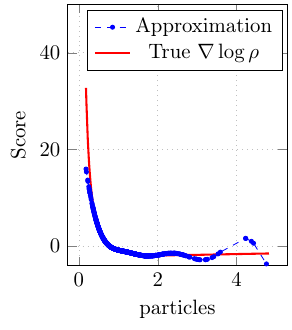}
\caption{$\lambda = 10^{-2}$.}
\label{fig:lambda 2 ln}
  \end{subfigure}%
  \hfill
  \begin{subfigure}[t]{0.30\textwidth}
    \centering

\includegraphics[width=\textwidth]{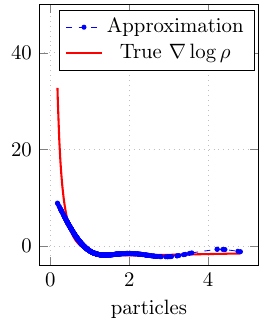}
\caption{$\lambda = 10^{0}$.}
\label{fig:lambda 0 ln}
  \end{subfigure}%
  \hfill
  \begin{subfigure}[t]{0.30\textwidth}
    \centering

\includegraphics[width=\textwidth]{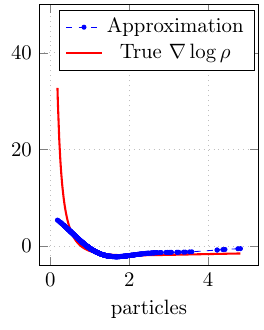}
\caption{$\lambda = 10^{1}$.}
\label{fig:lambda 1 ln}
  \end{subfigure}

  \caption{Plot of the score for values of $\lambda$ for the one-dimensional log-normal model with $P=1000$ particles.}
  \label{fig:lambda and gradlogs ln}
\end{figure}
In summary, a large $\lambda$ value imposes a strong prior on smoothness, as observed in Figures~\ref{fig:lambda vs rkhs} and \ref{fig:lambda vs rkhs ln}, which is relatively harmless in the Gaussian case, whose true score is smooth, but is detrimental in the log-normal case, where the score function contains higher-frequency components. In the log-normal case, reducing $\lambda$ decreases the bias in \Cref{fig:lambda vs L_rho ln}, allowing the estimator to capture sharp changes in $\nabla\log\rho(x)$ at a higher variance cost with worse conditioning, as demonstrated in \Cref{fig:lambda vs condition ln}. In contrast, for the Gaussian case in \Cref{fig:lambda vs L_rho}, the bias was of order $10^{-3}$, even with a moderate $\lambda$ value. Thus, decreasing $\lambda$ offers no significant benefit but leads to a more irregular score estimate in the 
RKHS. This behavior is also consistent with that in Figures~\ref{fig:lambda and gradlogs} and \ref{fig:lambda and gradlogs ln}.

\clearpage
\bibliography{References}
\bibliographystyle{plainnat}
\setcitestyle{authoryear,open={[},close={]}}
\end{document}